\documentclass[12pt]{article}
\usepackage[english]{babel}
\usepackage{epsfig}
\usepackage{amsfonts}
\usepackage{amssymb}
\usepackage{amsmath}
\usepackage{amsthm}
\usepackage{latexsym}
\usepackage{graphicx}
\usepackage{bm}
\usepackage{tabularx}
\usepackage{booktabs} 
\usepackage{enumerate}
\usepackage{caption}
\usepackage{wrapfig}  
\usepackage{geometry}
\usepackage{mathtools}  
\usepackage{eqnarray}
\usepackage{enumitem}

\usepackage[titletoc]{appendix}

\usepackage[numbers]{natbib}
\usepackage{comment}

\usepackage{bbm} 
\usepackage{url}
\usepackage{subfig}
\usepackage{hyperref}
\usepackage{color}
\usepackage{float}

\usepackage{hyperref}
\hypersetup{colorlinks=true, citecolor=blue, linkcolor=red, urlcolor=green}
 
\usepackage{bbm}
\usepackage{url}
\usepackage{color}
\usepackage{xcolor}
\usepackage[section]{placeins}

\graphicspath{{New_Graphs/}}

\theoremstyle{plain}
\newtheorem{theorem}{Theorem}[section]
\newtheorem{lemma}[theorem]{Lemma}
\newtheorem{corollary}[theorem]{Corollary}
\newtheorem{proposition}[theorem]{Proposition}

\theoremstyle{definition}

\theoremstyle{remark}
\newtheorem{remark}{Remark}

\newcommand{\shortdot}[1]{\raisebox{-0.4pt}{$\stackrel{\bullet}{#1}$}}
\newcommand{\updot}[1]{\raisebox{0.9pt}{$\stackrel{\bullet}{#1}$}} 

\title{The Number of Overlapping Customers in Erlang-A Queues: An Asymptotic Approach}
\author{ 
Young Myoung Ko
\\
Department of Industrial and Management Engineering
\\
Pohang University of Science Technology
\\
{youngko@postech.ac.kr} \\ \\
Jamol Pender \footnote{Corresponding Author}
\\
School of Operations Research and Information Engineering
\\
Cornell University
\\
{jjp274@cornell.edu} \\ \\
Jin Xu 
\\
School of Management
\\
Huazhong University of Science and Technology 
\\
{$\text{xu\_jin@hust.edu.cn}$}
}

\usepackage{natbib}
\usepackage{graphicx}
\usepackage{amsmath}
\date{}
\begin{document}

\maketitle

\begin{abstract}
In this paper, we investigate the number of customers that overlap or coincide with a virtual customer in an Erlang-A queue. Our study provides a novel approach that exploits fluid and diffusion limits for the queue to approximate the mean and variance of the number of overlapping customers. We conduct a detailed analysis of the fluid and diffusion limit differential equations to derive these approximations. We also construct new accurate approximations for the mean and variance of the waiting time in the Erlang-A queue by combining fluid limits with the polygamma function. Our findings have important implications for queueing theory and evaluating the overlap risk of more complicated service systems.   
\end{abstract}

\section{Introduction}

Have you ever wondered how many people you come in contact with when you visit a store or a retail shop? While this is an interesting question on its own, it has taken on greater significance in light of the COVID-19 pandemic. Most research on COVID-19 has utilized deterministic compartmentalized models to estimate infection rates and spread dynamics, as demonstrated in works such as \citet{dandekar2021safe, nguemdjo2020simulating, kaplan2020om}. However, we must acknowledge the importance of stochastic effects in determining the spread, as evidenced in works like \citet{drakopoulos2017network, palomo2020flattening, pang2020functional, forien2020epidemic, moein2021inefficiency}. As shoppers crowd stores to purchase essential goods such as water and non-perishable items, the risk of virus transmission increases. Service facilities and systems have implemented several measures, including the installation of air filters, transparent barriers, and mandatory mask-wearing for patrons. In addition, these systems have adopted various forms of social and physical distancing to minimize customer proximity, as outlined in \citet{bove2020restrict}.

Nonetheless, there are certain places where maintaining distance between people is not always feasible, such as in workplaces or crowded retail environments. In these situations, it becomes crucial to understand the extent to which customers overlap with one another. Recent research by \citet{kang2021queueing, palomooverlap, pender2021overlap, palomo2023overlap} has shed light on how to calculate customer overlap times in both single-server and multi-server queues. Specifically, \citet{kang2021queueing} demonstrates how overlaps can be used to calculate a new $R_0$ value for analyzing infection rates in compartmentalized epidemic models. Meanwhile, \citet{palomooverlap} proves that the overlap distribution is exponential for the M/M/1 queue and depends explicitly on the service distribution, and demonstrates via simulation that a similar result holds true in the non-exponential service setting as well.  Recent studies by \citet{pender2021overlap, palomo2023overlap} have brought fresh insights into the analysis of overlap times. \citet{pender2021overlap} analyzes overlap times in the infinite server queue, while \citet{palomo2023overlap} extends the analysis to the case of batch arrivals, highlighting the practical applications of this analysis for transportation systems like trains and buses.

Unfortunately, the existing literature on the number of overlaps has only focused on infinite server systems, as demonstrated by \citet{palomooverlap}. In an infinite server queue, the number of overlaps is equivalent to the queue length upon arrival and the number of additional arrivals during service, as there is no waiting. However, in the context of a finite server system, customer waiting can increase the number of overlaps. Consequently, the number of overlaps in an infinite server system can serve as a lower bound for any finite server system without abandonment. Given the prevalence of abandonment in real-world queueing systems, it is crucial to explore how abandonment and a finite number of servers impact the number of overlaps experienced by customers. While recent research by \citet{ko2022overlapping} has analyzed the total amount of overlapping time in a time-varying multi-server queue, their work has not investigated the number of overlaps or the impact of abandonment on the system. Thus, our work aims to address this gap in the literature.

Our paper presents a new analysis that provides critical insights into finite server systems with customer abandonment. Specifically, our analysis enables us to determine the number of overlaps that occur accurately and sheds light on the distribution of potential overlaps, allowing us to establish prediction intervals. This analysis has practical applications in preventing large overlaps, and can also serve as a valuable design tool for constructing appropriate overlaps by regulating arrival rates, service distributions, and the number of servers. Overall, our analysis offers a valuable contribution to the field of queueing systems and has the potential to inform decision-making in a variety of service system applications.


 \subsection{Main Contributions of the Paper}

\noindent The main contributions of this work can be summarized as follows:
\begin{itemize}
\item We derive the exact expressions for the Erlang-A fluid and diffusion differential equations.
\item We provide new approximations for the mean and variance of the number of overlapping customers.
\item We provide new approximations for the mean and variance of the waiting time in the Erlang-A queue using the digamma and trigamma functions.  
\end{itemize}


\subsection{Organization of the Paper}

The remainder of this paper is organized as follows. Section~\ref{secQMod} introduces the Erlang-A queueing model and its importance in stochastic network theory.  We also derive exact expressions for the fluid and diffusion variance differential equations.  In Section \ref{num_overlaps}, we introduce novel approximations for the mean and variance of the number of overlaps for a virtual customer. Our work also includes new approximations for the mean and variance of the wait time for the Erlang-A queue. We also validate the effectiveness of these approximations through simulation experiments. Lastly, in Section \ref{secConc}, we summarize our findings and suggest promising areas for future research in this area.

\section{The Erlang-A Queueing Model} \label{secQMod}

The Erlang-A queueing model is a fundamental queueing model in the applied probability literature.  The seminal work of \citet{mandelbaum1998strong} shows that the queue length process for an $M(t)/M/c+M$ queueing system $ Q\equiv \{ Q(t) | t \geq 0 \} $ is represented by the following stochastic, time changed integral equation:
\begin{eqnarray*} \label{queue_length}
Q(t) &=& Q(0) + \mathit{\Pi}_1 \left(\int^{t}_{0} \lambda(s) ds \right) -  \mathit{\Pi}_2 \left(\int^{t}_{0} \mu \cdot (Q(s) \wedge c )ds \right) \nonumber \\
&& - \mathit{\Pi}_3 \left(\int^{t}_{0} \theta \cdot (Q(s) - c )^+ds \right) ,
\end{eqnarray*}
where $  \mathit{\Pi}_i \equiv \{  \mathit{\Pi}_i(t) | t \geq 0 \} $
for $ i = {1,2,3} $ are i.i.d.\  standard (rate 1) Poisson processes and $(x \wedge y) = \min\{x,y\}$.  Thus, the sample path of the Erlang-A queue can be written in terms of three independent unit rate Poisson processes.  A deterministic time change for $\mathit{\Pi}_1 $ transforms it into a non-homogeneous Poisson arrival process with rate $ \lambda(t)$ that counts the customer arrivals that occurred in the time interval [0,t).  A random time change for the Poisson process $\mathit{\Pi}_2$, gives us a departure process that counts the number of customers that have received service from any one of the $c$ servers. We implicitly assume that the number of servers is $c \in \mathbb{Z}^+$ and that each server works at rate $\mu$. Finally, a random time change of $ \mathit{\Pi}_3$ counts the number of customers that abandon the system before beginning service.  We also assume that the abandonment distribution is exponential and the abandonment rate is equal to $\theta$.

The Erlang-A queueing model and its variants are very well studied in the queueing literature, see for example \citet{zeltyn2005call, whitt2006sensitivity, mandelbaum2007service, gurvich2014excursion,  engblom2014approximations,  pender2014gram, niyirora2016optimal, aktekin2016stochastic, braverman2017stein, pender2017sampling, pender2017approximating, bitton2019joint, azriel2019erlang, van2019economies}.  The Erlang-A queue has been so extensively analyzed since several important queueing models are special cases of it.  One special case is the infinite server queue.  The infinite server queue can be constructed from the Erlang-A queue in two ways.  The first way is to set the number of servers to infinity i.e., $c = \infty$.  The second way to construct the infinite server queue is to set the service rate $\mu$ equal to the abandonment rate $\theta$.  A second special case of the Erlang-A is when one sets the abandonment parameter $\theta =0$.  In this case, we shut off the abandonment process and we obtain the Erlang-C queue.  Finally, when we let $\theta$ get large i.e.,  $\theta \to \infty$, we obtain an Erlang-B queue, see for example \citet{hampshire2009time}. Thus, the Erlang-A queueing model encapsulates three special cases of queueing models and this reflects its importance as a model in the applied probability and queueing literature.

A key insight from \cite{halfin1981heavy} is that for multi-server queueing systems, it is natural to scale up the arrival rate and the number of servers simultaneously.  This scaling -- known as the \emph{Halfin-Whitt} scaling -- has been an important technique for modeling call centers in the queueing literature, see for example \citet{pang2007martingale}.  Since the $M(t)/M/c+M$ queueing process is a special case of a single node \emph{Markovian service network},
we can also construct an associated, \emph{uniformly accelerated} queueing process where both the new arrival rate
$\eta\cdot\lambda(t)$ and the new number of servers $\eta\cdot c$ are scaled by the same factor $\eta>0$.  Thus, using the \emph{Halfin-Whitt} scaling for the Erlang-A model, we arrive at the following sample path representation for the queue length process as
 \begin{eqnarray*}
Q^{\eta}(t) &=& Q^{\eta}(0) + \mathit{\Pi}_1 \left(\int^{t}_{0}  \eta \cdot \lambda(s) ds \right) -  \mathit{\Pi}_2 \left(\int^{t}_{0} \mu \cdot (Q^{\eta}(s) \wedge  \eta \cdot c )ds \right) \\
&&- \mathit{\Pi}_3 \left(\int^{t}_{0} \theta \cdot (Q^{\eta}(s) - \eta \cdot c )^+ds \right) \\
&=& Q^{\eta}(0) + \mathit{\Pi}_1 \left( \int^{t}_{0}  \eta \cdot \lambda(s) ds \right) -  \mathit{\Pi}_2 \left(\int^{t}_{0}  \eta \cdot \mu \cdot \left( \frac{Q^{\eta}(s)}{\eta} \wedge c \right)ds \right) \\
&&- \mathit{\Pi}_3 \left(\int^{t}_{0} \eta \cdot \theta \cdot \left( \frac{Q^{\eta}(s)}{\eta} - c \right)^+ds \right) .
\end{eqnarray*}

Applying the Halfin-Whitt scaling and taking the following limits gives us the \emph{fluid} models of \cite{mandelbaum1998strong}, i.e.,
 \begin{equation*}
\lim_{\eta\to\infty} \frac{1}{\eta} Q^{\eta}(t) =   q(t) \hspace{3mm}
\mathrm{a.s.}
\end{equation*}
where the deterministic process $q(t)$, the \emph{fluid mean}, is
governed by the one dimensional ordinary differential equation (ODE),
\begin{equation}
\label{fldmean}
 \shortdot{q}(t) = \lambda(t) - \mu \cdot (q(t) \wedge c) - \theta \cdot (q(t) - c)^+ .
\end{equation}
Moreover, if one takes a diffusion limit i.e.,
 \begin{equation*}
\lim_{\eta\to\infty} \sqrt{\eta} \left( \frac{1}{\eta} Q^{\eta}(t) -   q(t) \right) \Rightarrow \tilde{Q}(t), \hspace{3mm}
\end{equation*}
one obtains a diffusion process where the variance of the diffusion is given by the following ODE,
\begin{eqnarray}
\label{diffvar}
 \updot{\mathrm{Var}}\left[ \tilde{Q}(t) \right]
 & =& \lambda(t) + \mu \cdot (q(t) \wedge c) + \theta \cdot (q(t) - c)^+  \nonumber \\
&&- 2 \cdot  \mathrm{Var}\left[\tilde{Q}(t) \right] \cdot \left( \mu \cdot \{ q(t) < c \} +   \theta \cdot \{ q(t) \geq  c \} \right) .
\end{eqnarray}

In what follows, we provide an extensive analysis of the differential equations that are obtained from the fluid and diffusion limits.  Despite the equations being ubiquitous in the literature, this paper provides the first detailed analysis of the fluid and diffusion differential equation dynamics in the constant arrival rate case. We begin with the analysis of the differential equations of the fluid limits.   

\subsection{Fluid Analysis}

In this section, we provide results regarding the fluid limit dynamics when the arrival rate is constant.  Before we describe our main results regarding the fluid analysis, we provide a standard result for linear differential equations.  

\begin{lemma} \label{ode_soln}
Let $q(t)$ be the solution to the following differential equation
\begin{equation*}
\shortdot{q} = \lambda(t) - \mu(t) q(t) 
\end{equation*}
where $q(0) = q_0$.  Then the solution for any value of $t$ is given by 
\begin{eqnarray*}
q(t) &=& q_0  \exp\left\{ - \int^{t}_{0} \mu(s) ds \right\}  \\
&&+ \left( \exp\left\{ - \int^{t}_{0} \mu(s) ds \right\}  \cdot \left( \int^{t}_{0} \lambda(s) \exp\left\{  \int^{s}_{0} \mu(r) dr \right\}  ds \right) \right) \nonumber  .
\end{eqnarray*}
\begin{proof}
This follows from standard results on ordinary differential equations by varying parameters (see \cite{tenenbaum1985ordinary}).  
\end{proof}
\end{lemma}

Now that we characterize the dynamics of $q(t)$ with time-varying arrival rate $\lambda(t)$ and service rate $\mu(t)$, our next result provides the explicit solution to $q(t)$ in the steady state, i.e., when the arrival and service rates are both constant.  

\begin{proposition}\label{prop2.2}When $\lambda(t)$ and $\mu(t)$ are constant, the solution $q(t)$ to Equation \eqref{fldmean} is given as:
\begin{eqnarray*}
q(t) & = & \begin{cases}
\frac{\lambda-\mu c+\theta c}{\theta}+\left(q(0)-\frac{\lambda-\mu c+\theta c}{\theta}\right)e^{-\theta t} & \mbox{if }q(0)>c,\lambda>\mu c\\
\frac{\lambda-\mu c+\theta c}{\theta}+\left(q(0)-\frac{\lambda-\mu c+\theta c}{\theta}\right)e^{-\theta t}, & \mbox{if }q(0)>c,\lambda\leq\mu c,t\leq t_{1}^{*}\\
\frac{\lambda}{\mu}+\left(c-\frac{\lambda}{\mu}\right)e^{-\mu(t-t^{*})}, & \mbox{if }q(0)>c,\lambda\leq\mu c,t>t_{1}^{*}\\
\frac{\lambda}{\mu}+\left(q(0)-\frac{\lambda}{\mu}\right)e^{-\mu t} & \mbox{if }q(0)\leq c,\lambda\leq\mu c\\
\frac{\lambda}{\mu}+\left(q(0)-\frac{\lambda}{\mu}\right)e^{-\mu t}, & \mbox{if }q(0)\leq c,\lambda>\mu c,t\leq t_{2}^{*}\\
\frac{\lambda-\mu c+\theta c}{\theta}+\left(\frac{-\lambda+\mu c}{\theta}\right)e^{-\theta(t-t^{*})}, & \mbox{if }q(0)\leq c,\lambda>\mu c,t>t_{2}^{*}
\end{cases}
\end{eqnarray*}
where $t_{1}^{*}=\frac{\log\left(\frac{\theta q(0)-\lambda+\mu c-\theta c}{\mu c-\lambda}\right)}{\theta}$
and $t_{2}^{*}=\frac{\log\left(\frac{q(0)-\frac{\lambda}{\mu}}{c-\frac{\lambda}{\mu}}\right)}{\mu}.$

\begin{proof} When $q(0)>c$ and $\lambda>\mu c$, we apply Lemma
\ref{ode_soln} to the differential equation given in Equation \eqref{fldmean}.
When $q(0)>c$ and $\lambda\leq\mu c$, from time 0 to time $t_{1}^{*}$
(which is the time that the differential equation hits the value $c$),
we can apply Lemma \ref{ode_soln} to the differential equation to
obtain the solution. After $t_{1}^{*}$, we know the solution is larger
than $c$ so that it follows a new differential equation and the solution
is also given by Lemma \ref{ode_soln}. The discussion of the other
cases are similar, so we omit the details here. This completes the
proof. \end{proof} 
\end{proposition}

Now that we have completely described the dynamics of the fluid queue length as a function of time in the previous proposition. We can easily obtain the next corollary following Proposition \ref{prop2.2}. 

\begin{corollary}\label{cCor}
Suppose that $\lambda(t)$ and $\mu(t)$ are constant, and $q(\infty) = \lim_{t \to \infty} q(t)$. Then we have \begin{eqnarray*}
q(\infty) & = & \begin{cases}
\frac{\lambda-\mu c+\theta c}{\theta} & \mbox{if }\lambda>\mu c\\
\frac{\lambda}{\mu} & \mbox{if }\lambda\leq\mu c
\end{cases}
\end{eqnarray*}

\begin{proof}
The results follows from Proposition \ref{prop2.2} by letting $t\to \infty$.
\end{proof}
\end{corollary}

We then turn our attention to analyzing the diffusion variance differential equations.  Analyzing the fluid equations first is essential since the diffusion variance differential equations depend on the fluid dynamics in an explicit way.  

\subsection{Diffusion Variance Analysis}

In this section, we provide an analysis of the diffusion variance differential equations given in Equation \eqref{diffvar}.  We start with explicit solutions to the diffusion variance equations.  

\begin{proposition}\label{prop:var}
When $\lambda(t)$ and $\mu(t)$ are constant, the solution $v(t)$  to Equation \eqref{diffvar}, is given as:

\begin{align}\label{eq:variance}
v(t) =  \begin{cases}
\left(v(0)-q(0)-\frac{\mu c-\theta c}{\theta}\right)e^{-2\theta t}+\frac{\lambda}{\theta}+\left(q(0)-\frac{\lambda-\mu c+\theta c}{\theta}\right)e^{-\theta t}, & \mathrm{if}\ q(0)>c,\lambda>\mu c\\
\left(v(0)-q(0)-\frac{\mu c-\theta c}{\theta}\right)e^{-2\theta t}+\frac{\lambda}{\theta}+\left(q(0)-\frac{\lambda-\mu c+\theta c}{\theta}\right)e^{-\theta t}, & \mathrm{if}\ q(0)>c,\lambda\leq\mu c,t\leq t_{1}^{*}\\
(v(t_{1}^{*})-c)e^{-2\mu(t-t_{1}^{*})}+\frac{\lambda}{\mu}+\left(c-\frac{\lambda}{\mu}\right)e^{-\mu(t-t_{1}^{*})}, & \mathrm{if}\ q(0)>c,\lambda\leq\mu c,t>t_{1}^{*}\\
\left(v(0)-q(0)\right)e^{-2\mu t}+\left(1-e^{-\mu t}\right)\cdot\frac{\lambda}{\mu}+q(0)e^{-\mu t}, & \mathrm{if}\ q(0)\leq c,\lambda\leq\mu c\\
\left(v(0)-q(0)\right)e^{-2\mu t}+\left(1-e^{-\mu t}\right)\cdot\frac{\lambda}{\mu}+q(0)e^{-\mu t}, & \mathrm{if}\ q(0)\leq c,\lambda>\mu c,t\leq t_{2}^{*}\\
(v(t_{2}^{*})-c)e^{-2\theta(t-t_{2}^{*})}+\frac{\lambda}{\theta}+\left(\frac{-\lambda+\mu c}{\theta}\right)e^{-\theta(t-t_{2}^{*})}, & \mathrm{if}\ q(0)\leq c,\lambda>\mu c,t>t_{2}^{*}
\end{cases}
\end{align}
where
\begin{equation*}
    t^*_1 = \frac{\log\left( \frac{\theta q(0) - \lambda + \mu c - \theta c}{\mu c - \lambda} \right)}{\theta} \quad \mathrm{and} \quad  t^*_2 = \frac{\log\left( \frac{q(0) - \frac{\lambda}{\mu}}{c - \frac{\lambda}{\mu}} \right)}{\mu}.
\end{equation*}
\begin{proof}
We prove the results in Equation \eqref{eq:variance} by discussing the following four cases.

\noindent\textbf{CASE 1}:
For the first case with $q(0) > c$ and $\lambda > \mu c$, we have
\begin{eqnarray*}
 \updot{v} &=& \lambda + \mu \cdot  c + \theta \cdot (q(t) - c) -   2 \cdot \theta \cdot v(t).
\end{eqnarray*}
Using the theory of linear odes, this implies that 
\begin{eqnarray*}
v(t) &=& v(0) e^{-2\theta t} + e^{-2\theta t} \int^{t}_{0} e^{2\theta s} \left( \lambda + \mu c - \theta c  + \theta q(s) \right) ds \\
&=& v(0) e^{-2\theta t} + (1 - e^{-2\theta t} )  \cdot \left( \frac{\lambda + \mu c - \theta c }{2\theta} \right)  + \theta e^{-2\theta t} \int^{t}_{0} e^{2\theta s}   q(s) ds \\
&=& v(0) e^{-2\theta t} + (1 - e^{-2\theta t} )  \cdot \left( \frac{\lambda + \mu c - \theta c }{2\theta} \right)  \\
&+& \theta e^{-2\theta t} \int^{t}_{0} e^{2\theta s}   \left(  \frac{\lambda - \mu c + \theta c}{\theta} + e^{-\theta s} \left( q(0) - \frac{\lambda - \mu c + \theta c}{\theta} \right) \right)ds \\
&=& v(0) e^{-2\theta t} + (1 - e^{-2\theta t} )  \cdot \left( \frac{\lambda + \mu c - \theta c }{2\theta} \right)  \\
&+& \left( \frac{\lambda - \mu c + \theta c}{2 \theta} ( 1 -  e^{-2\theta t} ) + \left(q(0) - \frac{\lambda - \mu c + \theta c }{\theta} \right) \cdot \left(e^{-\theta t} - e^{-2\theta t} \right)  \right) \\
&=& \left(v(0) -q(0)+\frac{-\mu c + \theta c}{\theta}\right) e^{-2\theta t} +  \frac{\lambda}{\theta} 
+ \left(q(0) - \frac{\lambda - \mu c + \theta c }{\theta} \right) e^{-\theta t}.
\end{eqnarray*}

We thus showed the first item of Equation \eqref{eq:variance}.

\noindent\textbf{CASE 2}:
For the second case where $ q(0) > c$ and  $\lambda < \mu c$, we have 
\begin{eqnarray*}
 \updot{v} =  \begin{cases}
\lambda + \mu c + \theta q(t) - \theta c - 2 \theta v(t) , & \mathrm{if} \ t \leq t^* \\
\lambda + \mu q(t) - 2 \mu v(t), & \mathrm{if} \ t > t^* 
\end{cases} 
\end{eqnarray*}
where $t^*$ is equal to 
\begin{equation*}
    t^*_1 = \frac{\log\left( \frac{\theta q(0) + \lambda - \mu c + \theta c}{\mu c - \lambda} \right)}{\theta}.
\end{equation*}
Using the theory of linear odes, this implies for $t \leq t^*_1$ that
\begin{eqnarray*}
v(t) = \left(v(0) -q(0)+\frac{-\mu c + \theta c}{\theta}\right) e^{-2\theta t} +  \frac{\lambda}{\theta} 
+ \left(q(0) - \frac{\lambda - \mu c + \theta c }{\theta} \right) e^{-\theta t}.
\end{eqnarray*}
Lastly for $t > t^*_1$ we have that
\begin{eqnarray*}
v(t) &=& (v(t^*_1) - c) e^{-2\mu (t - t^*_1)} +  \frac{\lambda }{\mu} 
+ \left(c - \frac{\lambda }{\mu} \right) e^{-\mu (t - t^*_1)} .
\end{eqnarray*}

We hence proved the second and third items in Equation \eqref{eq:variance}.

\noindent\textbf{CASE 3}: 
For the third case where $ q(0) \leq c$ and $\lambda < \mu c$, we have
\begin{eqnarray*}
 \updot{v} &=& \lambda + \mu \cdot  q(t) -  2 \cdot \mu \cdot v(t).
\end{eqnarray*}
Using the theory of linear odes, this implies that 
\begin{eqnarray*}
v(t) &=& v(0) e^{-2\mu t} + e^{-2\mu t} \int^{t}_{0} e^{2\mu s} \left( \lambda + \mu q(s) \right) ds \\
&=& v(0) e^{-2\mu t} + \left( 1 - e^{-2\mu t} \right) \cdot \frac{\lambda}{2\mu} + \mu e^{-2\mu t} \int^{t}_{0} e^{2\mu s} q(s)ds \\
&=& v(0) e^{-2\mu t} + \left( 1 - e^{-2\mu t} \right) \cdot \frac{\lambda}{2\mu} + \mu e^{-2\mu t} \int^{t}_{0} e^{2\mu s} \left( \frac{\lambda}{\mu} + \left( q(0) - \frac{\lambda}{\mu} \right) e^{-\mu s} \right) ds \\
&=& v(0) e^{-2\mu t} + \left( 1 - e^{-2\mu t} \right) \cdot \frac{\lambda}{\mu} + \mu e^{-2\mu t} \int^{t}_{0}  \left( q(0) - \frac{\lambda}{\mu} \right) e^{\mu s}  ds \\
&=& v(0) e^{-2\mu t} + \left( 1 - e^{-2\mu t} \right) \cdot \frac{\lambda}{\mu} + \left( q(0) - \frac{\lambda}{\mu} \right) \left( e^{-\mu t} - e^{-2\mu t} \right)  \\
&=& \left( v(0) - q(0) \right) e^{-2\mu t} + \left( 1 - e^{-\mu t} \right) \cdot \frac{\lambda}{\mu} + q(0) e^{-\mu t} .
\end{eqnarray*}

We hence have the fourth item of Equation \eqref{eq:variance}.

\noindent\textbf{CASE 4}:
For the fourth case where $\lambda > \mu c$ and $ q(0) \leq c$, we have 
\begin{eqnarray*}
 \updot{v} &=&  \begin{cases}
\lambda + \mu q(t) - 2 \mu v(t), & \mathrm{if} \ t \leq t^*_2 \\
\lambda + \mu c + \theta q(t) - \theta c - 2 \theta v(t) , & \mathrm{if} \ t > t^*_2 
\end{cases} 
\end{eqnarray*}
where $t^*$ is equal to 
\begin{equation*}
    t^*_2 = \frac{\log\left( \frac{q(0) - \frac{\lambda}{\mu}}{c - \frac{\lambda}{\mu}} \right)}{\mu}.
\end{equation*}
Using the theory of linear odes, this implies for $t \leq t^*_2$ that
\begin{eqnarray*}
v(t) &=& v(0) e^{-2\mu t} + e^{-2\mu t} \int^{t}_{0} e^{2\mu s} \left( \lambda + \mu q(s) \right) ds \\
&=& v(0) e^{-2\mu t} + \frac{\lambda}{2\mu} \left(1 - e^{-2\mu t} \right) \\
&+& \mu  e^{-2\mu t}  \int^{t}_{0} e^{2\mu s}  \left( \frac{\lambda}{\mu} + \left( q(0) - \frac{\lambda}{\mu} \right) e^{-\mu s}\right) ds \\
&=& v(0) e^{-2\mu t} + \left( \frac{\lambda}{2\mu} + \frac{\lambda }{2\mu} \right) \left(1 - e^{-2\mu t} \right) + \mu  e^{-2\mu t}  \int^{t}_{0} e^{\mu s}  \left( q(0) - \frac{\lambda }{\mu} \right)  ds \\
&=& v(0) e^{-2\mu t} + \left( \frac{\lambda}{2\mu} + \frac{\lambda }{2\mu} \right) \left(1 - e^{-2\mu t} \right) + \left( q(0) - \frac{\lambda }{\mu} \right) \left(  e^{-\mu t}  - e^{-2\mu t} \right)  \\
&=& (v(0) - q(0)) e^{-2\mu t} + \frac{\lambda}{\mu}  + \left( q(0) - \frac{\lambda }{\mu} \right) e^{-\mu t}  .
\end{eqnarray*}
Lastly for $t > t^*_2$ we have that
\begin{eqnarray*}
v(t) = (v(t^*_2) - c) e^{-2\theta (t - t^*_2)} +  \frac{\lambda  }{\theta}  
+ \left( \frac{-\lambda + \mu c  }{\theta}  \right) e^{-\theta (t - t^*_2)} .
\end{eqnarray*}

We thus proved the fifth and sixth items of Equation \eqref{eq:variance}.
\end{proof}
\end{proposition}

Now that we have provided an analysis of the transient dynamics of the diffusion variance differential equations, we analyze the steady state behavior of the diffusion variance differential equation in the following result.

\begin{corollary}
Let $v(t)$ be the solution to Equation \eqref{diffvar}, then in steady state we have 
\begin{eqnarray*}
v(\infty) &=&  \begin{cases}
\frac{\lambda}{\mu}, \mathrm{if} \ \lambda \leq \mu c \\
\frac{\lambda  }{\theta} , \mathrm{if} \ \lambda > \mu c .
\end{cases} 
\end{eqnarray*}
\begin{proof}
    The result follows directly from Proposition \ref{prop:var} by letting $v \to \infty$ in Equation \eqref{eq:variance}.
\end{proof}

\end{corollary}

What is apparent from the steady state variance is that, it is the same as the steady state mean queue length when $\lambda \leq \mu c$.  This is because the system is effectively an infinite server queue in this regime.  We also see that the when the arrival rate $\lambda$ is larger than the maximum service rate $\mu c$, then the steady state variance is $\lambda/\theta$.  This is interesting because when $\theta < \mu$, the mean is smaller than the variance and when $\theta \geq \mu$, the steady state variance is smaller than the steady state mean queue length.  This describes over-dispersion and under-dispersion respectively.  For more on the over-dispersion and under-dispersion relationships between the mean and variance in the Erlang-A queue see \citet{daw2019new}.

Now that we have a good understanding of the fluid mean and diffusion variance of the Erlang-A queue, we will show how to use these results in the context of studying overlaps in the Erlang-A queue.  

\section{The Number of Overlaps} \label{num_overlaps}

In this section, we introduce the virtual overlap process, which counts the number of customers that the virtual customer will overlap with during their time in the queue.  This process is important from a epidemiology perspective since it would be the number of people that would need to be contact traced for an exposure if the virtual customer had an infectious disease, see for example \citet{kang2021queueing, palomo2020flattening}.  

Like the infinite server setting, the virtual customer will also overlap with the number present in the queue and the number of customers that arrive during the service time of the virtual customer.  However, unlike the infinite server setting, the virtual customer, who we assume does not abandon from the queue, must also overlap with customers that also arrive during their wait for service.  Thus, the overlap process can be written as the following expression in terms of the queue length $Q(t)$, the virtual waiting time $W(t)$, and the service time of the virtual customer $\mathcal{S}$ i.e. 
\begin{eqnarray} \label{overlap_process}
O(t) &=& \underbrace{N\left(t+\mathcal{S} + W(t) \right) - N(t)}_{\text{Arrivals During Service and Wait}}  + \underbrace{Q(t)}_{\text{Queue Upon Arrival}}.
\end{eqnarray}
Using the above representation, we can compute the mean number of overlaps at time $t$ by taking the expectation of the overlap process. Thus, the mean number of overlaps can be written as 
\begin{eqnarray*}
\mathbb{E}\left[ O(t) \right] &=&\mathbb{E}\left[  N\left(t+\mathcal{S} + W(t) \right) - N(t)  + Q(t) \right] \\
&=&\mathbb{E}\left[  N\left(t+\mathcal{S} + W(t) \right) - N(t)\right]  + \mathbb{E}\left[ Q(t) \right] \\
&=& \lambda \mathbb{E}\left[ \mathcal{S}\right]  + \lambda \mathbb{E}\left[  W(t) \right]  + \mathbb{E}\left[ Q(t) \right] .
\end{eqnarray*}
Unfortunately, the transient mean queue length and the transient mean wait time are not known in closed form for the Erlang-A queue, except in the case where $\mu = \theta$.  This is a major difference between the Erlang-A and the infinite server queue.  In the infinite server queue, the mean wait time is zero and mean queue can be written as an explicit integral with respect to the service time distribution, see for example \citet{eick1993mt, eick1993physics}.  Consequently, we will use limit theory to approximate the mean number of overlaps.

\subsection{The Transient Mean Number of Overlaps}

In this section, we show how to approximate the transient mean number of overlap using asymptotic analysis.  We first leverage the results of \citet{mandelbaum1998strong}, which proves almost sure limit theorems for the queue length process in the Halfin-Whitt regime.  However, we also need results for the virtual waiting time in order to fully analyze the overlap process.  To do this, we will exploit a recent result by \citet{massey2018dynamic}, which proves the following theorem.  

\begin{theorem} \label{v_wait_time}
Let $W^{\eta}(t)$ be the virtual wait time of a customer at time t who is not going to abandon, in the scaled process.  Then we have 
\begin{eqnarray*}
\lim_{\eta \to \infty} W^{\eta}(t) &{\buildrel a.s \over =} & w(t),
\end{eqnarray*}
where $w(t)$ satisfies the following equation
\begin{eqnarray} \label{fluid_wait}
w(t) &=& \frac{1}{\theta} \log \left( 1 + \frac{\theta \cdot ( q(t) - c)^+}{\mu c }\right) .
\end{eqnarray}
Moreover, as $t \to \infty$ we have that when $\lambda > \mu c$
\begin{eqnarray*} \label{fluid_wait_2}
\lim_{t \to \infty} w(t) &=&  \frac{1}{\theta} \log \left( \frac{\lambda}{\mu c }\right) .
\end{eqnarray*}
\begin{proof}
See Theorem 6 of \citet{massey2018dynamic}. 
\end{proof}
\end{theorem}

This result shows that the limiting virtual waiting time is a function of the fluid queue length function $q(t)$.  When the fluid queue length is less than the number of servers $c$, then the virtual waiting time is equal to zero.  Moreover, when the fluid queue length is greater than the number of servers $c$, then the virtual waiting time is positive. What is more important is that the fluid limit for the virtual waiting time yields a deterministic function of time.  Thus, we are able to approximate the virtual waiting time with a non-random function of time. 
Now that we have a limiting expression for the virtual waiting time, we can define the scaled overlap process as follows
\begin{eqnarray}\label{eq:3-3}
O^{\eta}(t) = N^{\eta} \left(t+\mathcal{S} + W^{\eta}(t)  \right) - N^{\eta}(t)  + Q^{\eta}(t).
\end{eqnarray}

We will use the scaled overlap process to prove our main result of the paper, which gives us a deterministic approximation for the transient mean number of overlaps in the Erlang-A queue.  

\begin{theorem} \label{overlap_fluid}
Let $O^{\eta}(t)$ be the scaled number of people that a customer who arrives at time $t$ will overlap with.  Then, we have
\begin{eqnarray*}
\lim_{\eta \to \infty} \frac{1}{\eta} \mathbb{E} \left[ O^{\eta}(t) \right] &=& \frac{\lambda}{\mu} + \frac{\lambda }{\theta} \log \left( 1 + \frac{\theta \cdot ( q(t) - c)^+}{\mu c }\right) + q(t) .
\end{eqnarray*}
where $q(t)$ is the solution to Equation \eqref{fldmean}.  Moreover, when $t \to \infty$ we have that
\begin{eqnarray}\label{eq:number-steady-state}
\lim_{t \to \infty} \lim_{\eta \to \infty} \frac{1}{\eta} \mathbb{E} \left[ O^{\eta}(t) \right] &=& \begin{cases}
\frac{2\lambda}{\mu} , \mathrm{if} \ \lambda \leq \mu c \\
\frac{\lambda}{\mu} + \frac{\lambda }{\theta} \log \left(  \frac{\lambda }{\mu c } \right) + c + \frac{\lambda - \mu c}{\theta} , \mathrm{if} \ \lambda > \mu c. 
\end{cases}
\end{eqnarray}
\begin{proof}
Based on Equation \eqref{eq:3-3} and Theorem \ref{v_wait_time}, we have
\begin{eqnarray*}
\lim_{\eta \to \infty} \frac{1}{\eta} \mathbb{E} \left[ O^{\eta}(t) \right] &=& \frac{1}{\eta} \mathbb{E} \left[\left( N^{\eta} \left(t+\mathcal{S}+ W^{\eta}(t) \right) - N^{\eta}(t) \right) \right] + \frac{1}{\eta} \mathbb{E} \left[Q^{\eta}(t) \right]\\
&=& \frac{1}{\eta} \left( \eta \lambda \mathbb{E} \left[ \mathcal{S} \right] + \eta \lambda \mathbb{E} \left[  W^{\eta}(t) \right] \right)  + q(t) \\
&=& \frac{\lambda}{\mu} + \lambda w(t)  + q(t) \\
&=& \frac{\lambda}{\mu} + \frac{\lambda }{\theta} \log \left( 1 + \frac{\theta \cdot ( q(t) - c)^+}{\mu c }\right) + q(t) .
\end{eqnarray*}
In particular, when we look at steady state, we have that when $\lambda \leq \mu c$,
\begin{eqnarray*}
\lim_{t \to \infty} \lim_{\eta \to \infty} \frac{1}{\eta} \mathbb{E} \left[ O^{\eta}(t) \right] &=& \frac{\lambda}{\mu} + \frac{\lambda }{\theta} \log \left( 1 + \frac{\theta \cdot ( q(\infty) - c)^+}{\mu c }\right) + q(\infty)  \\
&=& \frac{\lambda}{\mu} + \frac{\lambda }{\theta} \log \left( 1 + \frac{\theta \cdot \left( \frac{\lambda}{\mu}  - c \right)^+}{\mu c }\right) + \frac{\lambda}{\mu}   \\
&=& \frac{2\lambda}{\mu}.  
\end{eqnarray*}
When $\lambda > \mu c$, we have
\begin{eqnarray*}
\lim_{t \to \infty} \lim_{\eta \to \infty} \frac{1}{\eta} \mathbb{E} \left[ O^{\eta}(t) \right] &=& \frac{\lambda}{\mu} + \frac{\lambda }{\theta} \log \left( 1 + \frac{\theta \cdot ( q(\infty) - c)^+}{\mu c }\right) + q(\infty)  \\
&=& \frac{\lambda}{\mu} + \frac{\lambda }{\theta} \log \left(  \frac{\lambda }{\mu c } \right) + c + \frac{\lambda - \mu c}{\theta} .
\end{eqnarray*}
This completes the proof.
\end{proof}
\end{theorem}

It is important to note that the expectation for Theorem \ref{overlap_fluid} is necessary.  Without the expectation the overlap process would explicitly depend on $\mathcal{S}$ and would be random.  It is also worth noting that in the case $\lambda \leq \mu c$, the steady state mean number of overlaps is equal to the steady state mean number of overlaps in the infinite server queue as well.  It also does not depend on the staffing level $c$ since it behaves like an infinite server in the steady state setting.

As we stated earlier, the Erlang-A model is a generalization of multiple queueing system. We now discuss several special cases of the Erlang-A model in the following based on our results in Theorem \ref{overlap_fluid}.

\begin{remark}(M/M/$\infty$ system): When $c\rightarrow\infty$,
all the customers will enter the service directly upon arrival. The
system with constant arrival and service rates becomes the M/M/$\infty$
system. We find from Equation \eqref{eq:number-steady-state} that the
number of overlapped customers will converge to $\frac{2\lambda}{\mu}$.
This result matches our analysis in \cite{xu2023queueing}.  It is important to distinguish between letting $c \to \infty$ and letting $\theta \to \mu$. Notably, as $c$ approaches infinity, the waiting time converges to zero. Conversely, when $\theta$ tends towards $\mu$, the waiting time does not approach zero. Consequently, although the queue lengths are identical and share the same sample path structure, the expected number of overlaps differs between the systems. This observation of the two different systems emphasizes that the number of overlaps is dependent on the customer's experience instead of only the queue length process.  In particular, the customers in the case where $c \to \infty$ all have the same experience in the system, however, in the case where $\theta \to \mu$, abandoning customers have a different overlap experience than those who wait in the queue and then get served.  Finally, it is also worth noting that the diffusion limits of the two systems are different, see for example page 167 of \citet{mandelbaum1998strong}. This further highlights the difference in the number of overlaps experienced by customers in the two systems.  
\end{remark}

\begin{remark}(Erlang-B system): When letting the abandonment rate
$\theta\rightarrow\infty$, we have the Erlang-B system. According
to Theorem \ref{overlap_fluid}, the expected number of people that
the virtual customer overlaps in the steady state is $\frac{2\lambda}{\mu}$
when $\lambda\leq\mu c$, and is $\frac{\lambda}{\mu}+c$ when $\lambda>\mu c$.
The reason is that no customer will wait in the queue, so the virtual customer will always enter the service directly. According to Corollary \ref{cCor}, the virtual customer will overlap with $\frac{\lambda}{\mu}$ customers in the system when $\lambda\leq\mu c$, and $c$ customers when $\lambda > \mu c$. Moreover, during the service time, the virtual customer is expected to overlap with $\frac{\lambda}{\mu}$ newly arrived customers. Note that although the abandonment rate is infinity, we still assume that every customer will arrive at the system first before abandon, which explains why the virtual customer will overlap with new arrivals.
\end{remark}

\begin{remark} (Erlang-C system): When shutting off the abandonment process, we have the Erlang-C system. Based on Equation \eqref{eq:number-steady-state}, the number of overlapped customers in the fluid limit becomes $2\frac{\lambda}{\mu}$ for $\lambda\leq \mu c$. The reason is that the fluid limit of the virtual waiting time becomes 0 in the steady state, according to Theorem \ref{v_wait_time}. In the fluid limit, the virtual customer will overlap with $\frac{\lambda}{\mu}$ customers upon arrival, and other $\frac{\lambda}{\mu}$ customers during service.  
    
\end{remark}

The transient mean queue length and wait time are not known in closed form.  However, if we condition on the queue length, the waiting time distribution is known in closed form for the Erlang-A model.  By conditioning on the number of customers ahead of you given that you are waiting, it is easily seen that the waiting time has a hypoexponential distribution i.e.,
\begin{eqnarray*}
    W_k = \sum^{k}_{j=0} Y_j
\end{eqnarray*}
where $Y_j \sim \mathrm{Exp}(\mu c + \theta \cdot j )$ and $k$ is the number of customers that are ahead of the customer upon the arrival.  Note that $j$ is allowed to be zero when the queue length is identical to the number of servers since in this case, the customer needs to wait an exponential amount of time with rate $\mu c$.  With this hypoexponential representation, we can compute the conditional mean and variance of the waiting time.  The conditional mean waiting time is
\begin{eqnarray} \label{cond_wait}
  \mathbb{E}\left[   W_k \right] &=& \sum^{k}_{j=0} \frac{1}{\mu c + \theta \cdot j }\nonumber \\
  &=& \frac{1}{\theta} \sum^{k}_{j=0} \frac{1}{\frac{\mu c}{\theta} + j }\nonumber \\
  &=& \frac{1}{\theta} \left( \psi\left( \frac{\mu c}{\theta} + k+1 \right) - \psi\left( \frac{\mu c}{\theta}  \right) \right) \label{cond_mean_wait}
\end{eqnarray}
where $\psi(x)$ is the digamma function.  Moreover, the variance is
\begin{eqnarray}
  \mathrm{Var}\left[   W_k \right] &=& \sum^{k}_{j=0} \frac{1}{(\mu c + \theta \cdot j )^2 } \nonumber \\
  &=& \frac{1}{\theta^2} \sum^{k}_{j=0} \frac{1}{ \left( \frac{\mu c}{\theta} + j \right)^2 }\nonumber \\
  &=& \frac{1}{\theta^2} \left( \psi^{(1)}\left( \frac{\mu c}{\theta}  \right) - \psi^{(1)}\left( \frac{\mu c}{\theta} + k+1 \right) \right) \label{cond_var_wait}
\end{eqnarray}
where $\psi^{(1)}(x)$ is the trigamma function.  One should note that both the digamma and the trigamma functions are special cases of the the Hurwitz-Riemann zeta function  defined as
\begin{eqnarray*}
  \zeta(s,\alpha) = \sum^{\infty}_{j=0} \frac{1}{(j + \alpha)^s}.
\end{eqnarray*}

Now that we have explicit expressions for the mean and variance of the waiting time given that there are $k$ customers in front of the current arrival, we should be able to leverage the fluid limits for the queue length process to approximate the mean and variance and of the waiting time at any time $t$.  

We know that a first order Taylor expansion around the mean is given by the following expression
\begin{eqnarray*}
    f(Q(t)) &\approx& f \left(\mathbb{E}\left[ Q(t) \right] \right) + f' \left(\mathbb{E}\left[ Q(t) \right] \right) \cdot \left( Q(t) - \mathbb{E}\left[ Q(t) \right] \right).
\end{eqnarray*}
Thus, the mean of the function can be approximated by 
\begin{eqnarray*}
    \mathbb{E}\left[ f(Q(t)) \right] &\approx& f \left(\mathbb{E}\left[ Q(t) \right] \right) \approx f \left(q(t)\right)
\end{eqnarray*}
 and the variance of the function can be approximated by 
 \begin{eqnarray*}
    \mathrm{Var} \left[f(Q(t)) \right] &\approx&  f' \left(\mathbb{E}\left[ Q(t) \right] \right)^2  \cdot \mathrm{Var} \left[ Q(t) \right] \approx f' \left(q(t)  \right)^2  \cdot v(t).
\end{eqnarray*}
Finally, it is also important to compute an approximation for the covariance as well since we will need it later for computing the variance of the number of overlaps in the Erlang-A queue.
\begin{eqnarray*} \label{linear_noise}
\mathrm{Cov}\left[ f(Q(t)), Q(t) \right] &\approx& \mathrm{Cov}\left[ f(\mathbb{E}\left[ Q(t) \right]) + f'(\mathbb{E}\left[ Q(t) \right]) ( Q(t) - \mathbb{E}\left[ Q(t) \right]), Q(t) \right]  \\ 
&=&  f'(\mathbb{E}\left[ Q(t) \right]) \cdot \mathrm{Var}\left[ Q(t) \right]\\
&\approx& f'(q(t)) \cdot v(t).
\end{eqnarray*}
It is worth mentioning that the covariance result is well known in physics as the linear noise approximation.  We will show in the sequel how to use these approximations for estimating the mean and variance of the waiting time for the Erlang-A queue.  We will also find that these approximations of the waiting time are essential for computing the mean and variance of the number of overlaps in the Erlang-A queue as well.

By using a first order Taylor expansion around the mean queue we obtain the following expression as an approximation for the mean waiting time
\begin{eqnarray} \label{mean_wait_digamma}
  \mathbb{E}\left[ W(t) \right] &=& \mathbb{E} \left[ \frac{1}{\theta} \left( \psi\left( \frac{\mu c}{\theta} + (Q(t) - c)^+ \right) - \psi\left( \frac{\mu c}{\theta}  \right) \right) \right] \nonumber \\
  &=&   \frac{1}{\theta} \mathbb{E} \left[ \left( \psi\left( \frac{\mu c}{\theta} + (Q(t) - c)^+ \right)  \right) \right] - \frac{1}{\theta} \psi\left( \frac{\mu c}{\theta}  \right)\nonumber \\
  &\approx&  \frac{1}{\theta} \left( \psi\left( \frac{\mu c}{\theta} + (q(t) - c )^+ \right) - \psi\left( \frac{\mu c}{\theta}  \right) \right) \label{mean_wait_digamma2} .
\end{eqnarray}

Moreover, a first order Taylor expansion also yields the following approximation for the variance of the waiting time 
\begin{eqnarray} \label{var_wait_trigamma}
    \mathrm{Var}\left[   W(t) \right] &=&  \mathbb{E} \left[ \mathrm{Var} \left[ W(t) \mid Q(t)  \right] \right] +\mathrm{Var} \left[ \mathbb{E} \left[ W(t) \mid Q(t) \right] \right] \nonumber \\
  &=&  \mathbb{E} \left[  \frac{1}{\theta^2} \left( \psi^{(1)}\left( \frac{\mu c}{\theta}  \right) - \psi^{(1)}\left( \frac{\mu c}{\theta} + (Q(t) - c)^+  \right) \right) \right]\nonumber \\
  &+& \mathrm{Var} \left[  \frac{1}{\theta} \left( \psi\left( \frac{\mu c}{\theta} + (Q(t) - c)^+ \right) - \psi\left( \frac{\mu c}{\theta}  \right) \right)\right] \nonumber \\
  &=&  \frac{1}{\theta^2} \left( \psi^{(1)}\left( \frac{\mu c}{\theta}  \right)  - \mathbb{E} \left[ \psi^{(1)}\left( \frac{\mu c}{\theta} + (Q(t) - c)^+  \right)  \right] \right)\nonumber \\
  &+&  \frac{1}{\theta^2} \mathrm{Var} \left[   \psi\left( \frac{\mu c}{\theta} + (Q(t) - c)^+ \right)  \right] \nonumber \\
    &\approx&  \frac{1}{\theta^2} \left( \psi^{(1)}\left( \frac{\mu c}{\theta}  \right)  - \mathbb{E} \left[ \psi^{(1)}\left( \frac{\mu c}{\theta} + (Q(t) - c)^+  \right)  \right] \right)\nonumber \\
  &+&  \frac{1}{\theta^2} \cdot \psi^{(1)}\left( \frac{\mu c}{\theta} + (q(t) - c)^+ \right)^2  \cdot \mathrm{Var} \left[ Q(t) \right]\nonumber  \\
  &\approx&  \frac{1}{\theta^2} \left( \psi^{(1)}\left( \frac{\mu c}{\theta}  \right)  - \psi^{(1)}\left( \frac{\mu c}{\theta} + (q(t) - c)^+  \right) \right) \nonumber \\
  &+&  \frac{1}{\theta^2} \cdot \psi^{(1)}\left( \frac{\mu c}{\theta} + (q(t) - c)^+ \right)^2 \cdot \{ q(t) > c \}  \cdot v(t) .\label{var_wait_trigamma2}
\end{eqnarray}
Note that we have replaced the value $k$ in the conditional mean and variance formulas with the fluid limit queue length at time $t$.


\subsection{The Transient Variance of the Number of Overlaps}

In addition to the mean, we are also interested in approximating the variance of the number of overlaps.  With the variance, we are able to understand the variation around our approximations of the mean.  This implies that we can construct prediction intervals for the number of overlaps one might expect at any time $t$.  We know from Equation \eqref{overlap_process} that the overlap process satisfies the following equation
\begin{eqnarray*}
O(t) = N\left(t+\mathcal{S} + W(t) \right) - N(t)  + Q(t).
\end{eqnarray*}
Thus, the variance of the number overlap is equal to 
\begin{eqnarray*}
\mathrm{Var}[O(t)] &=& \mathrm{Var}\left[ N\left(t+\mathcal{S} + W(t) \right) - N(t)  + Q(t) \right] \\
&=& \mathrm{Var}\left[ N\left(t+\mathcal{S} \right) - N(t) + N\left(t+\mathcal{S} + W(t) \right) - N\left(t+\mathcal{S} \right)  + Q(t) \right] \\
&=& \mathrm{Var}\left[ N\left(t+\mathcal{S} \right) - N(t) \right] + \mathrm{Var}\left[ N\left(t+\mathcal{S} + W(t) \right) - N\left(t+\mathcal{S} \right)  + Q(t) \right] \\
&=& \mathrm{Var}\left[ N\left(t+\mathcal{S} \right) - N(t) \right] + \mathrm{Var}\left[ N\left( W(t) \right)  + Q(t) \right] \\
&=& \mathrm{Var}\left[ N\left(t+\mathcal{S} \right) - N(t) \right] + \mathrm{Var}\left[ N\left( W(t) \right) \right] + \mathrm{Var}\left[ Q(t) \right] \\
&+& 2\mathrm{Cov}\left[ N\left( W(t) \right)  , Q(t) \right] \\
&=&  \lambda \mathbb{E}\left[\mathcal{S}\right] + \lambda \mathbb{E}\left[ W(t) \right] + \lambda^2 \mathrm{Var}[ \mathcal{S}] + \lambda^2 \mathrm{Var}[ W(t) ] + \mathrm{Var}\left[ Q(t) \right] \\
&+& 2 \lambda \mathrm{Cov}\left[  W(t)  , Q(t) \right] \\
&=&  \frac{\lambda}{\mu} + \frac{\lambda^2}{\mu^2} + \lambda \mathbb{E}\left[ W(t) \right] + \lambda^2 \mathrm{Var}[ W(t) ] + \mathrm{Var}\left[ Q(t) \right] \\
&+& 2 \lambda \mathrm{Cov}\left[  W(t)  , Q(t) \right].
\end{eqnarray*}

The exact variance of both $W(t)$ and $Q(t)$ and their covariance are unknown.  Fortunately, because of the fluid and diffusion limits, we have approximations for the mean and variance of the queue length and waiting times.  Thus, it remains for us to derive a transient approximation for the covariance between the waiting time and the queue length at time $t$.  In order to compute the covariance of the two processes, we use a conditioning argument based on the queue length, which is natural given the conditional mean waiting time formula of Equation \eqref{cond_wait}.  We outline this argument below.   

\begin{eqnarray*}
\mathrm{Cov}\left[ W(t), Q(t) \right] &=& \mathbb{E}\left[  W(t) \cdot Q(t) \right] - \mathbb{E}\left[  W(t) \right] \cdot \mathbb{E}\left[ Q(t) \right] \\
&=& \sum^{\infty}_{k=0} k \cdot \mathbb{E}\left[ W(t)   | Q(t) = k  \right] \cdot \mathbb{P} \left( Q(t) = k \right) -  \mathbb{E}\left[  W(t) \right] \cdot \mathbb{E}\left[ Q(t) \right] \\
&=& \sum^{\infty}_{k=c} k \cdot \left( \sum^{k-c}_{j=0} \frac{1}{\mu c + \theta \cdot j }  \right) \cdot \mathbb{P}\left( Q(t) = k \right) - \mathbb{E}\left[  W(t) \right] \cdot \mathbb{E}\left[ Q(t) \right] \\
&=&\sum^{\infty}_{k=c} k \cdot \left( \psi\left(\frac{\mu c}{\theta} + k-c \right) - \psi \left( \frac{\mu c}{\theta} \right) \right) \cdot \mathbb{P}\left( Q(t) = k \right) - \mathbb{E}\left[  W(t) \right] \cdot \mathbb{E}\left[ Q(t) \right] \\
&=& \mathbb{E}\left[ Q(t) \left( \psi\left(\frac{\mu c}{\theta} + (Q(t)-c)^+ \right) - \psi \left( \frac{\mu c}{\theta} \right)  \right) \right] - \mathbb{E}\left[  W(t) \right] \cdot \mathbb{E}\left[ Q(t) \right] \\
&=& \mathrm{Cov}\left[ \psi\left(\frac{\mu c}{\theta} + (Q(t)-c)^+ \right) - \psi \left( \frac{\mu c}{\theta} \right)  , Q(t) \right] \\
&=& \mathrm{Cov}\left[ \psi\left(\frac{\mu c}{\theta} + (Q(t)-c)^+ \right)   , Q(t) \right] \\
&\approx& \mathrm{Var}\left[  Q(t) \right] \cdot \psi^{(1)} \left(\frac{\mu c}{\theta} + (q(t)-c)^+ \right)  \cdot \{ q(t) > c \} \\
&\approx& v(t) \cdot \psi^{(1)} \left(\frac{\mu c}{\theta} + (q(t)-c)^+ \right)  \cdot \{ q(t) > c \} .
\end{eqnarray*}

Using this result, we can approximate the covariance between the waiting time and the queue length as the product of the variance and the trigamma function, which is evaluated at the fluid limit. By using this approximation, we arrive at the following formula for approximating the transient variance of the number of overlaps.

\begin{eqnarray} 
\mathrm{Var}[O(t)] &=&  \frac{\lambda}{\mu} + \frac{\lambda^2}{\mu^2} + \lambda \mathbb{E}\left[ W(t) \right] + \lambda^2 \mathrm{Var}[ W(t) ] + \mathrm{Var}\left[ Q(t) \right] + 2 \lambda \mathrm{Cov}\left[  W(t)  , Q(t) \right]\nonumber \\
&\approx&  \frac{\lambda}{\mu} + \frac{\lambda^2}{\mu^2} + \frac{\lambda}{\theta} \left( \psi\left( \frac{\mu c}{\theta} + (q(t) - c)^+ \right) - \psi\left( \frac{\mu c}{\theta}  \right) \right)  + v(t) \nonumber \\
&+& \frac{\lambda^2 }{\theta^2} \left( \psi^{(1)}\left( \frac{\mu c}{\theta}  \right) - \psi^{(1)}\left( \frac{\mu c}{\theta} + (q(t) - c)^+ \right) \right) \nonumber   \\
&+& \frac{\lambda^2}{\theta^2} \cdot \psi^{(1)}\left( \frac{\mu c}{\theta} + (q(t) - c)^+ \right)^2 \cdot \{ q(t) > c \}  \cdot v(t)\nonumber  \\
&+& 2 \lambda v(t) \cdot \psi^{(1)} \left(\frac{\mu c}{\theta} + (q(t)-c)^+ \right) \cdot \{ q(t) > c \}  . \label{var_wait}
\end{eqnarray}
We can approximate the steady variance by setting $ t \to \infty$.  Thus, when $\lambda > \mu c $
\begin{eqnarray*}
\lim_{t \to \infty} \mathrm{Var}[O(t)] &\approx&  \frac{\lambda}{\mu} + \frac{\lambda^2}{\mu^2} + \frac{\lambda}{\theta} \left( \psi\left( \frac{\lambda}{\theta} \right) - \psi\left( \frac{\mu c}{\theta}  \right) \right)  + \frac{\lambda}{\theta} + \frac{\lambda^3 }{\theta^3} \psi^{(1)}\left( \frac{\lambda}{\theta}  \right)^2 \\
&+& \frac{\lambda^2 }{\theta^2} \left( \psi^{(1)}\left( \frac{\mu c}{\theta}  \right) - \psi^{(1)}\left( \frac{\lambda}{\theta}  \right) \right)  + 2 \frac{\lambda^2}{\theta}  \cdot \psi^{(1)} \left( \frac{\lambda}{\theta}  \right) 
\end{eqnarray*}
 and when $\lambda < \mu c $
 \begin{eqnarray*}
\lim_{t \to \infty} \mathrm{Var}[O(t)] &\approx&  \frac{2\lambda}{\mu} + \frac{\lambda^2}{\mu^2}. 
\end{eqnarray*}

It is important to note that the last equation is precisely the same variance of the infinite server queue setting.  Thus, when the queue is underloaded, the number of overlaps behaves similarly to an infinite server queue.

\subsection{Numerical Experiments}

To understand how our approximations for the mean queue length and the mean number of overlapping customers perform, we present eight different numerical examples below.  Before we provide the examples, we list the parameter values for each of the examples in Table \ref{Table_1}.  In all simulation examples, we simulate the sample paths 10,000 times to produce each curve.

\begin{table}[!htbp]
\centering%
\begin{tabular}{|c|c|c|c|c|c|}
\hline 
Parameters & $\lambda$ & $\mu$ & $\theta$ & $C$ & $Q(0)$\tabularnewline
\hline 
\hline 
Figure (a) & 10 & 1 & 0.5 & 30 & 10\tabularnewline
\hline 
Figure (b) & 10 & 1 & 0.5 & 30 & 50\tabularnewline
\hline 
Figure (c) & 10 & 1 & 2 & 30 & 10\tabularnewline
\hline 
Figure (d) & 10 & 1 & 2 & 30 & 50\tabularnewline
\hline 
Figure (e) & 40 & 1 & 0.5 & 30 & 10\tabularnewline
\hline 
Figure (f) & 40 & 1 & 0.5 & 30 & 50\tabularnewline
\hline 
Figure (g) & 40 & 1 & 2 & 30 & 10\tabularnewline
\hline 
Figure (h) & 40 & 1 & 2 & 30 & 50\tabularnewline
\hline 
\end{tabular}

\caption{Parameters for Examples}
\label{Table_1}
\end{table}

\begin{figure}[!htbp]
\centering
\subfloat[]{\includegraphics[scale=.31]{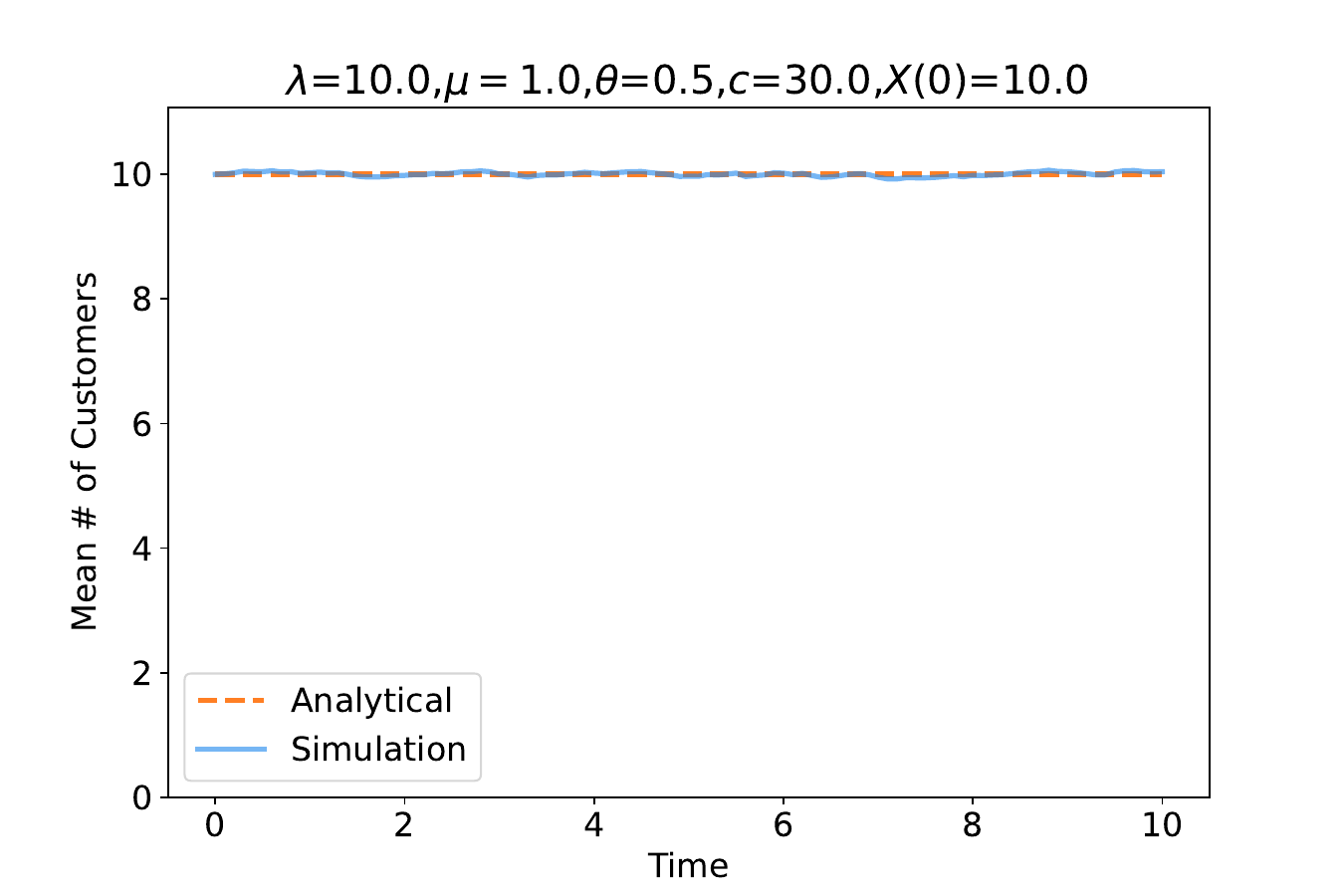}}
\subfloat[]{\includegraphics[scale=.31]{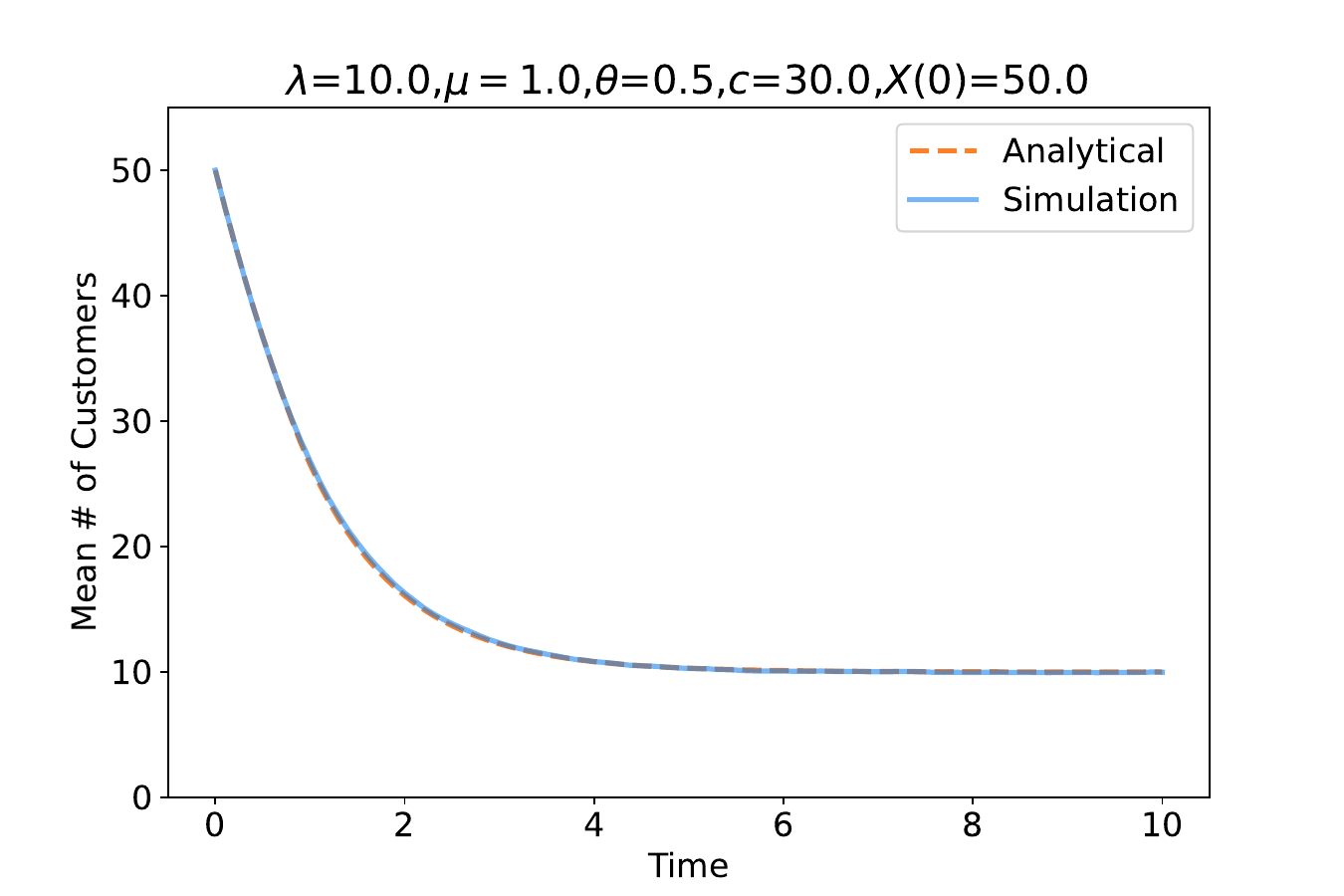}}
\\
\subfloat[]{\includegraphics[scale=.31]{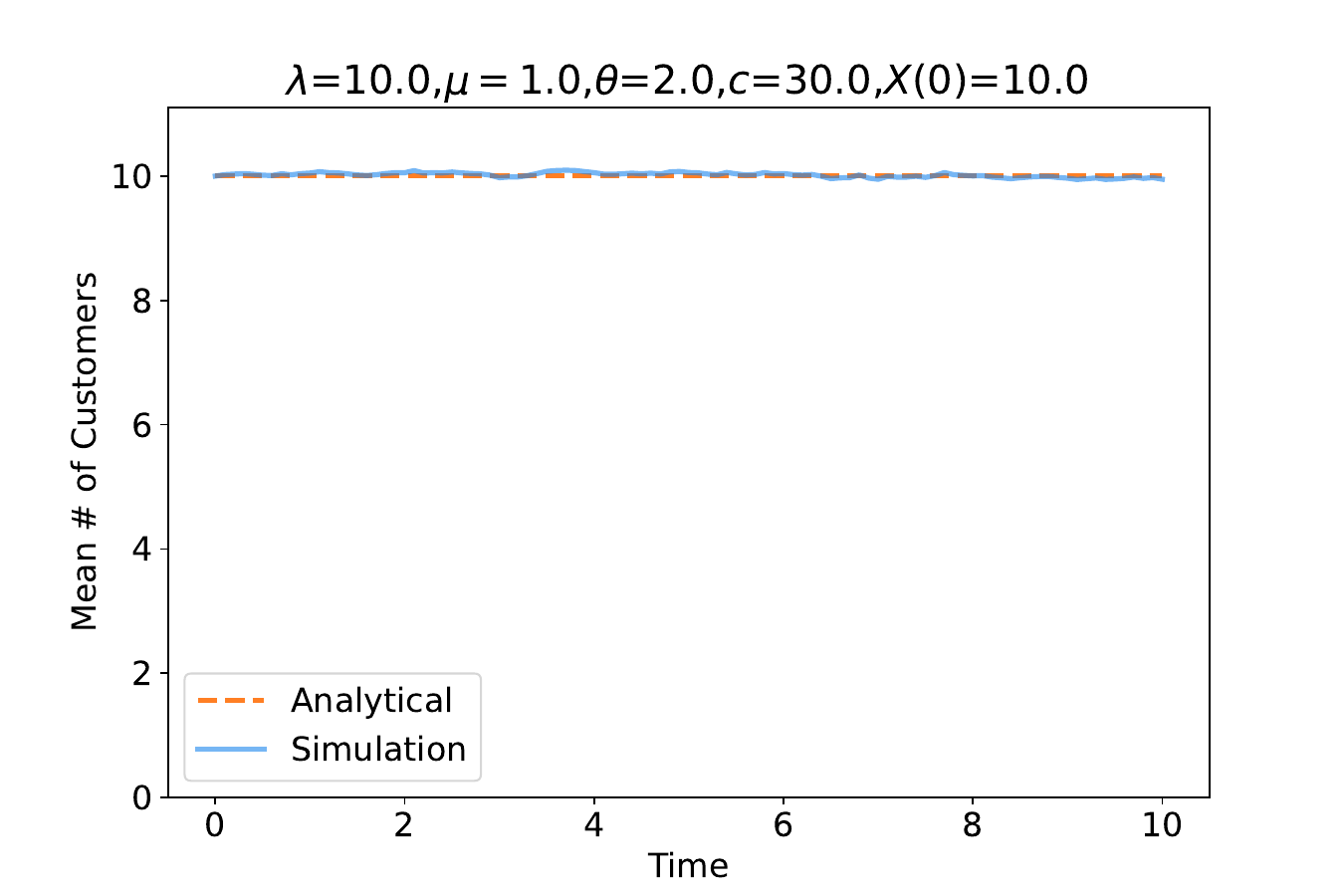}}
\subfloat[]{\includegraphics[scale=.31]{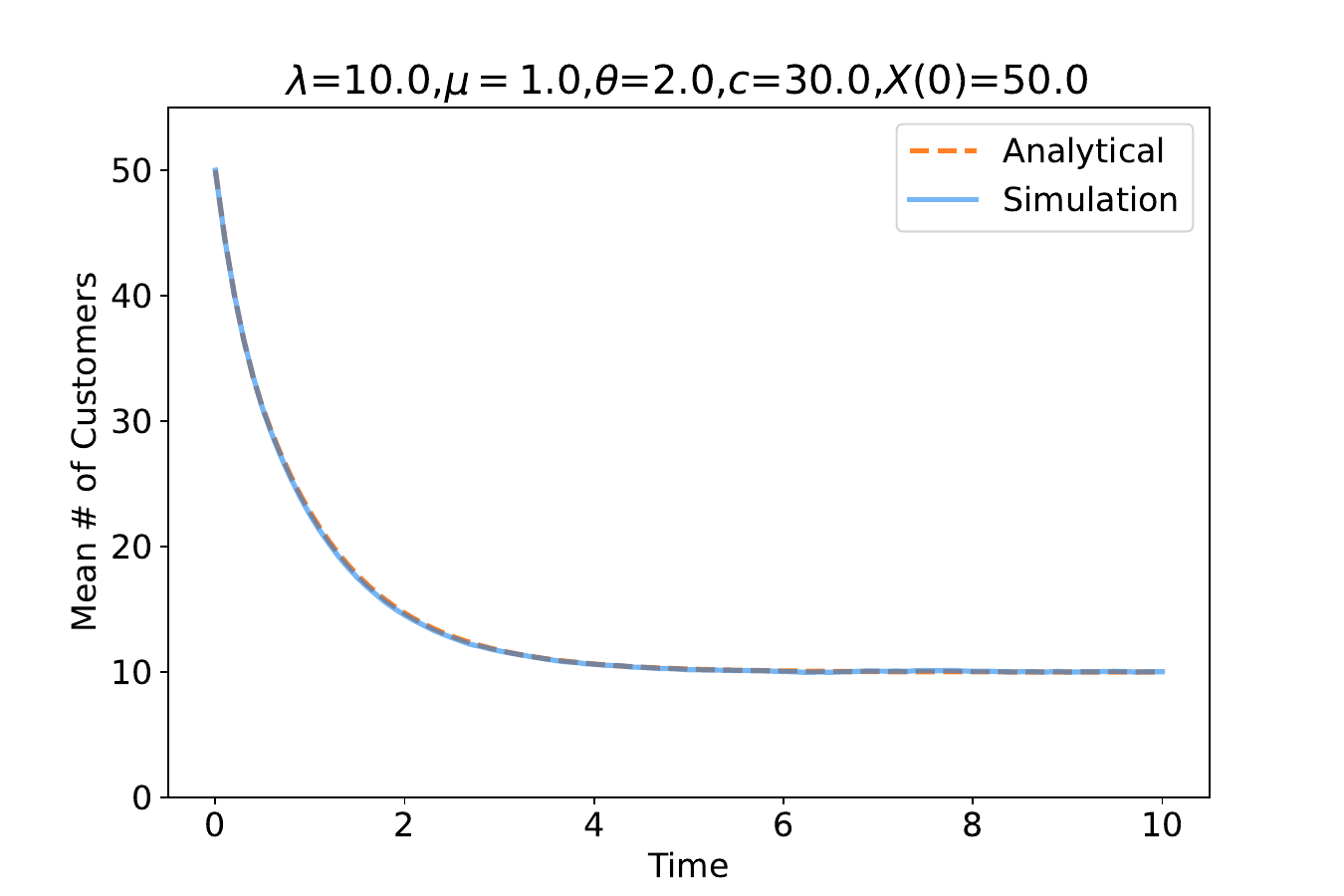}}
\\
\subfloat[]{\includegraphics[scale=.31]{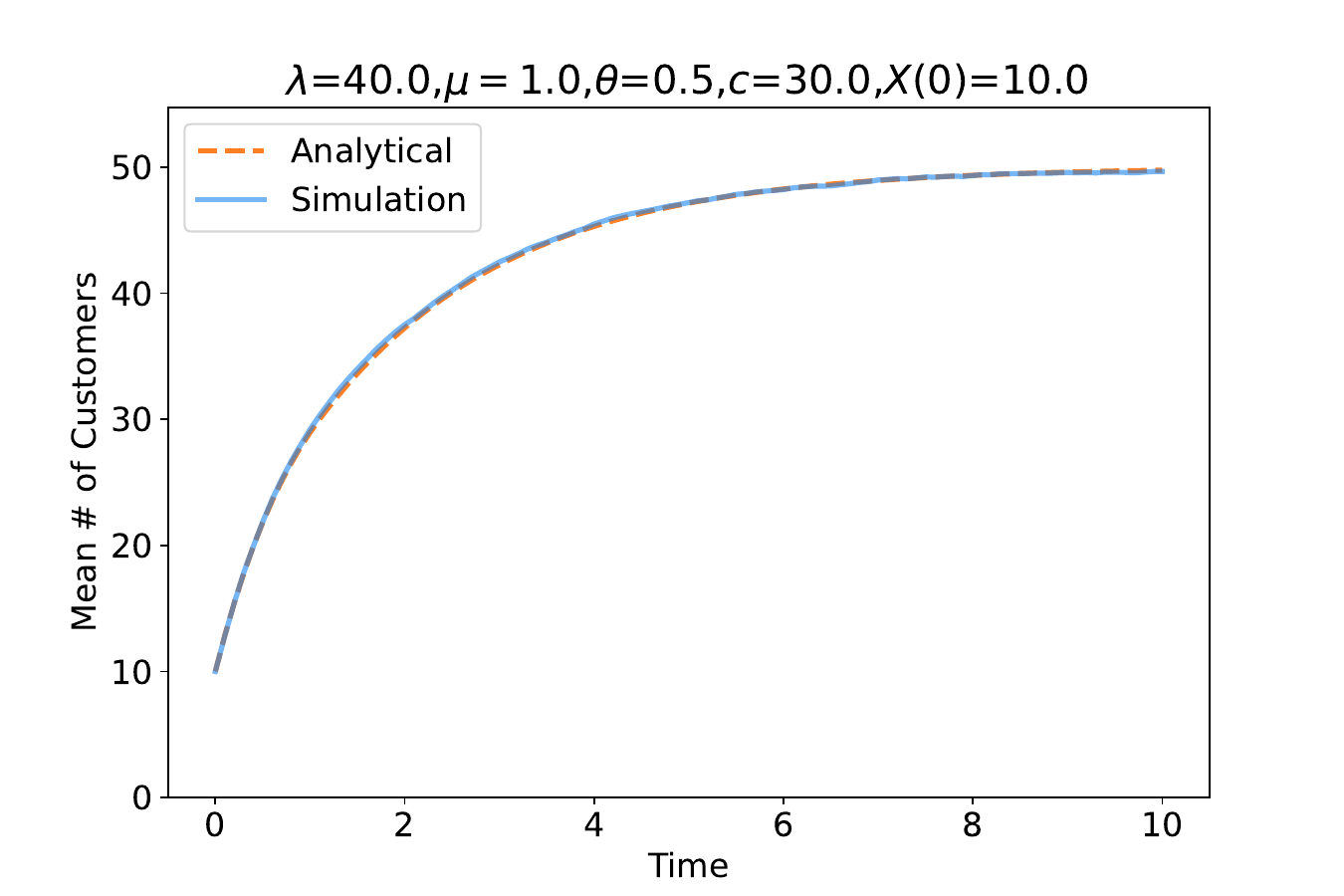}}
\subfloat[]{\includegraphics[scale=.31]{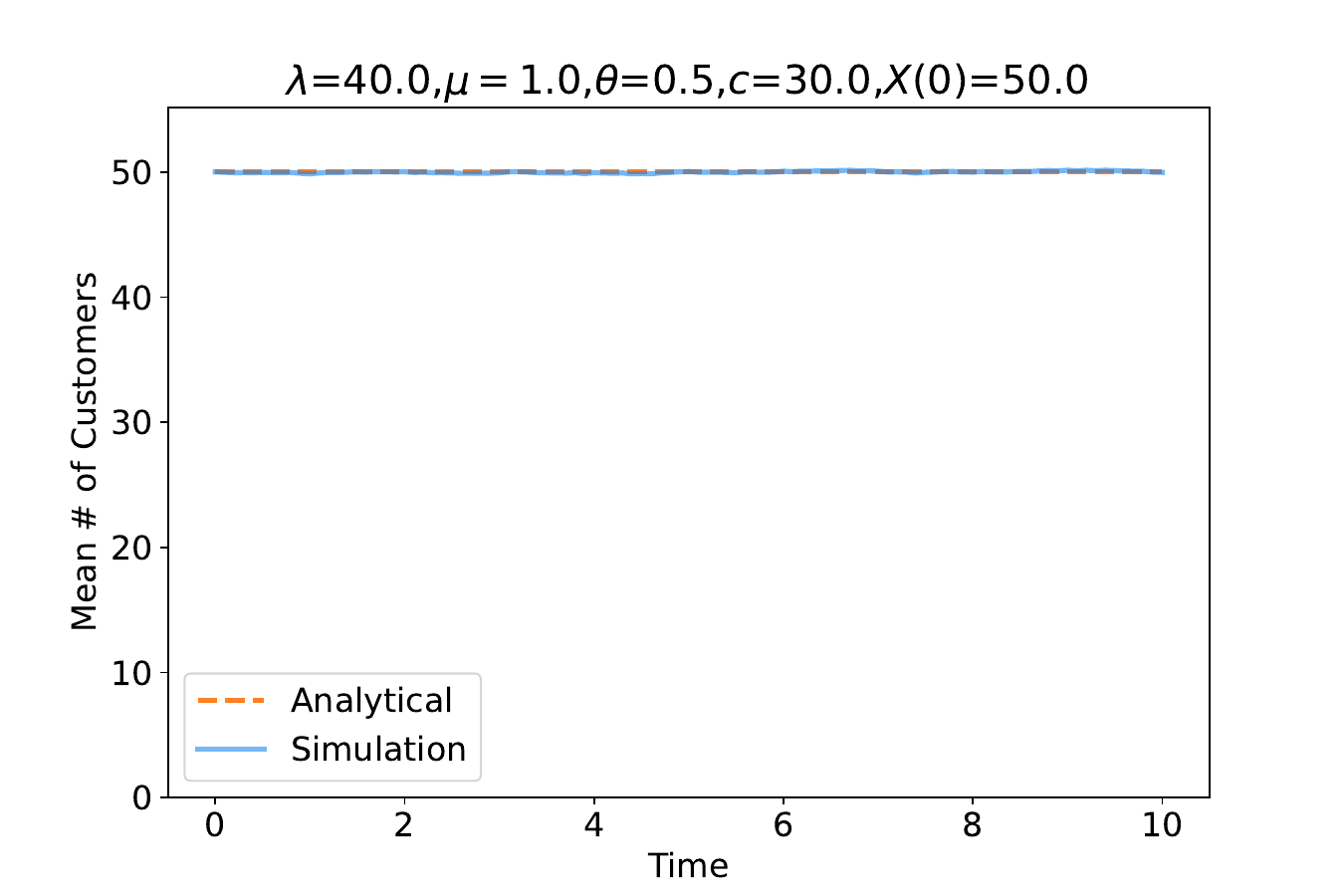}}
\\
\subfloat[]{\includegraphics[scale=.31]{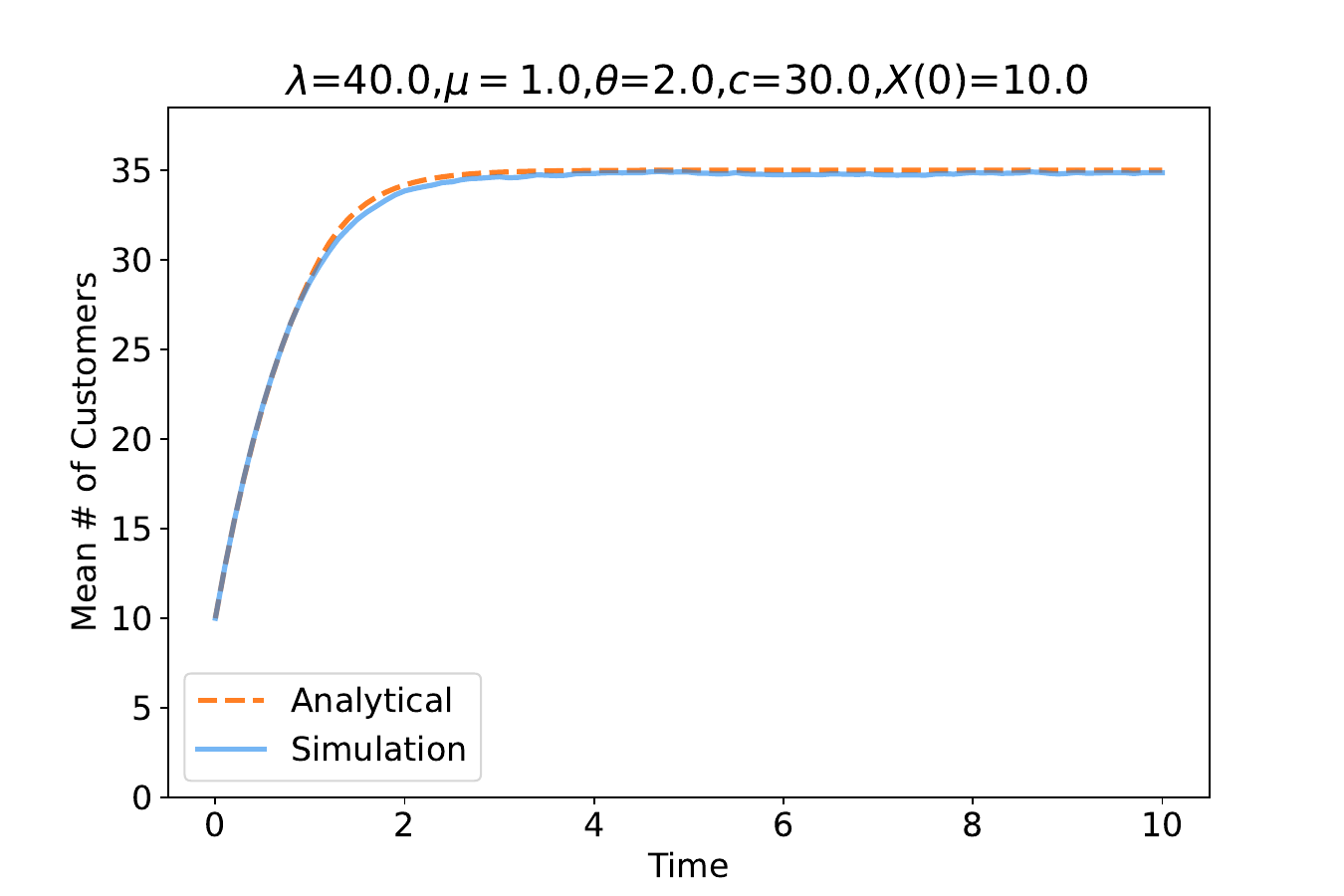}}
\subfloat[]{\includegraphics[scale=.31]{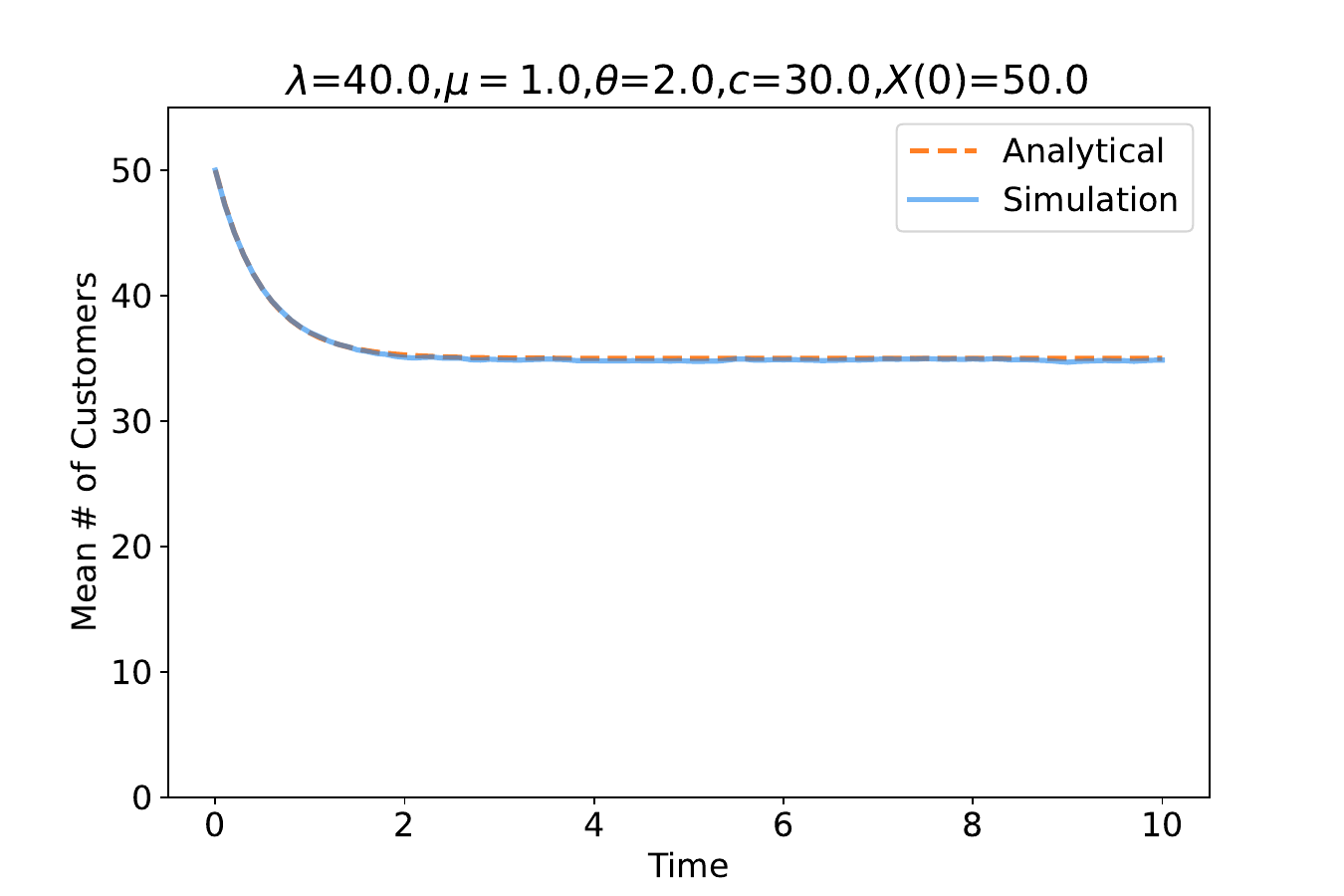}}
 \caption{Fluid Mean Number in System vs. Simulation.}  
\label{Figure_1}
\end{figure}

In Figure \ref{Figure_1}, we present a plot of the simulated mean queue length, which has been approximated by employing the fluid limit as expressed in Equation \eqref{fldmean}, denoted as ``analytical" in the figures.  Our findings indicate that across various parameter settings outlined in Table \ref{Table_1}, the fluid approximation consistently and accurately estimates the mean dynamics when juxtaposed with the simulated values. This observation underscores the reliability and effectiveness of the fluid approximation in accurately characterizing the behavior of the mean queue length. 

\begin{figure}[!htbp]
\centering
\subfloat[]{\includegraphics[scale=.31]{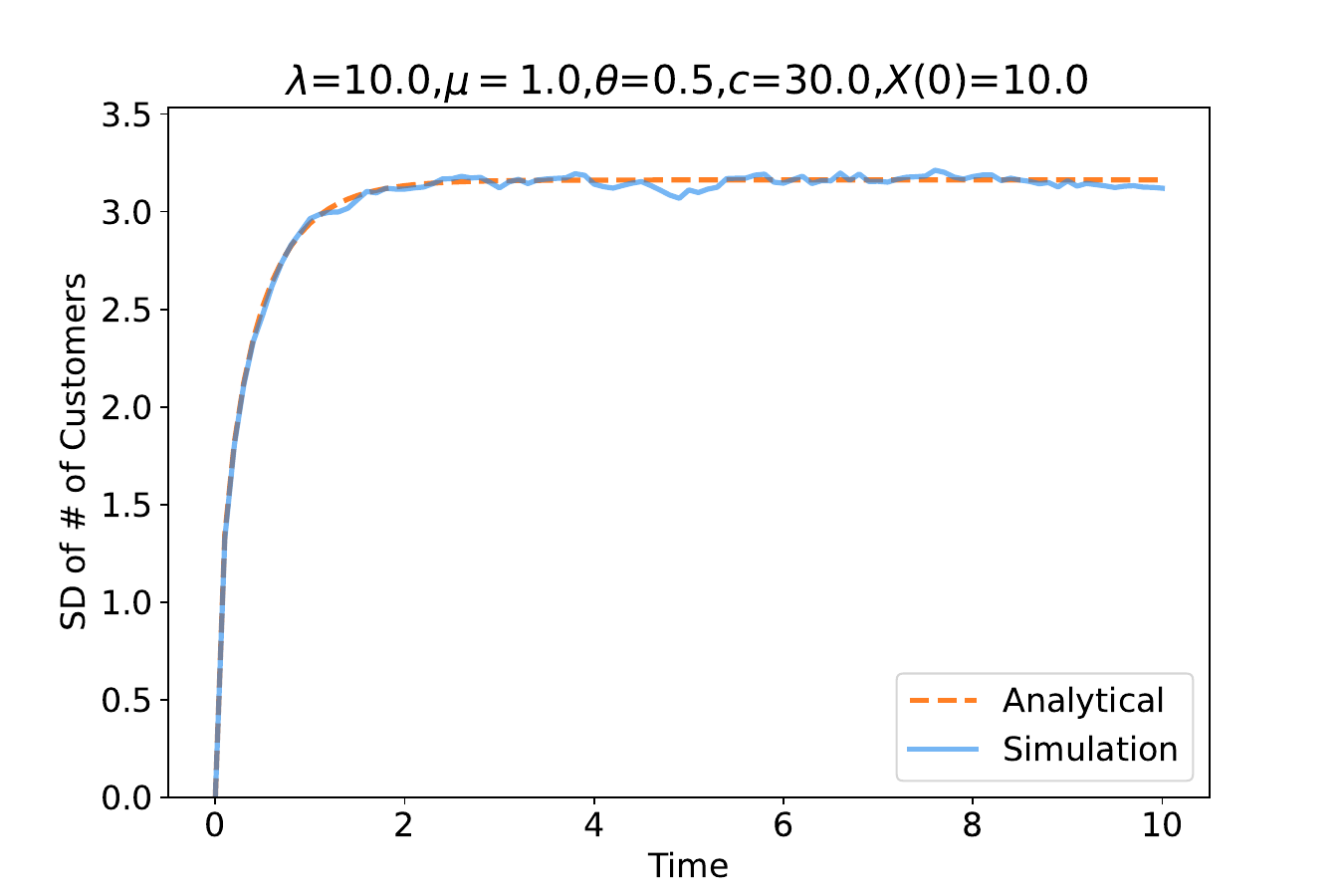}}
\subfloat[]{\includegraphics[scale=.31]{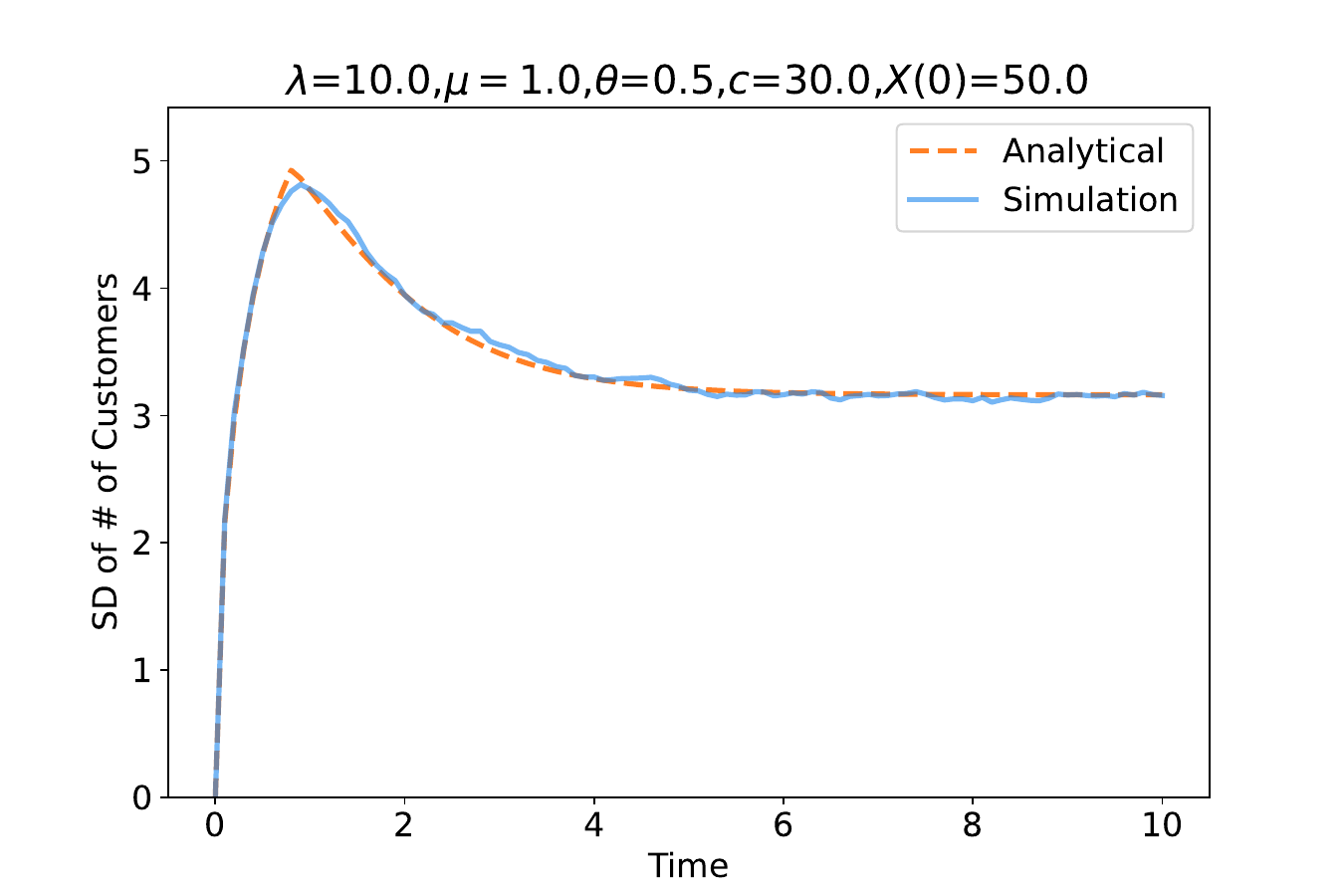}}
\\
\subfloat[]{\includegraphics[scale=.31]{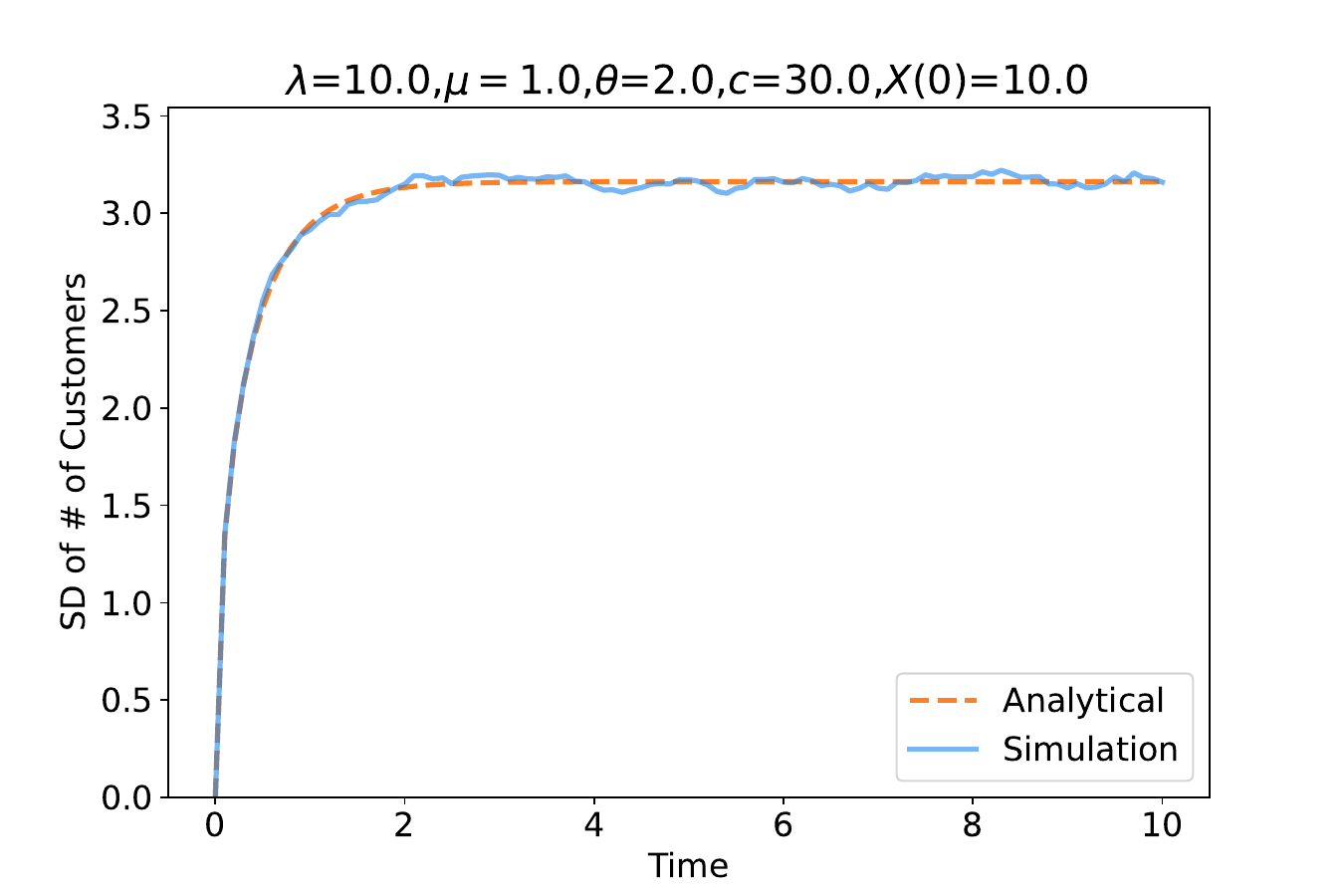}}
\subfloat[]{\includegraphics[scale=.31]{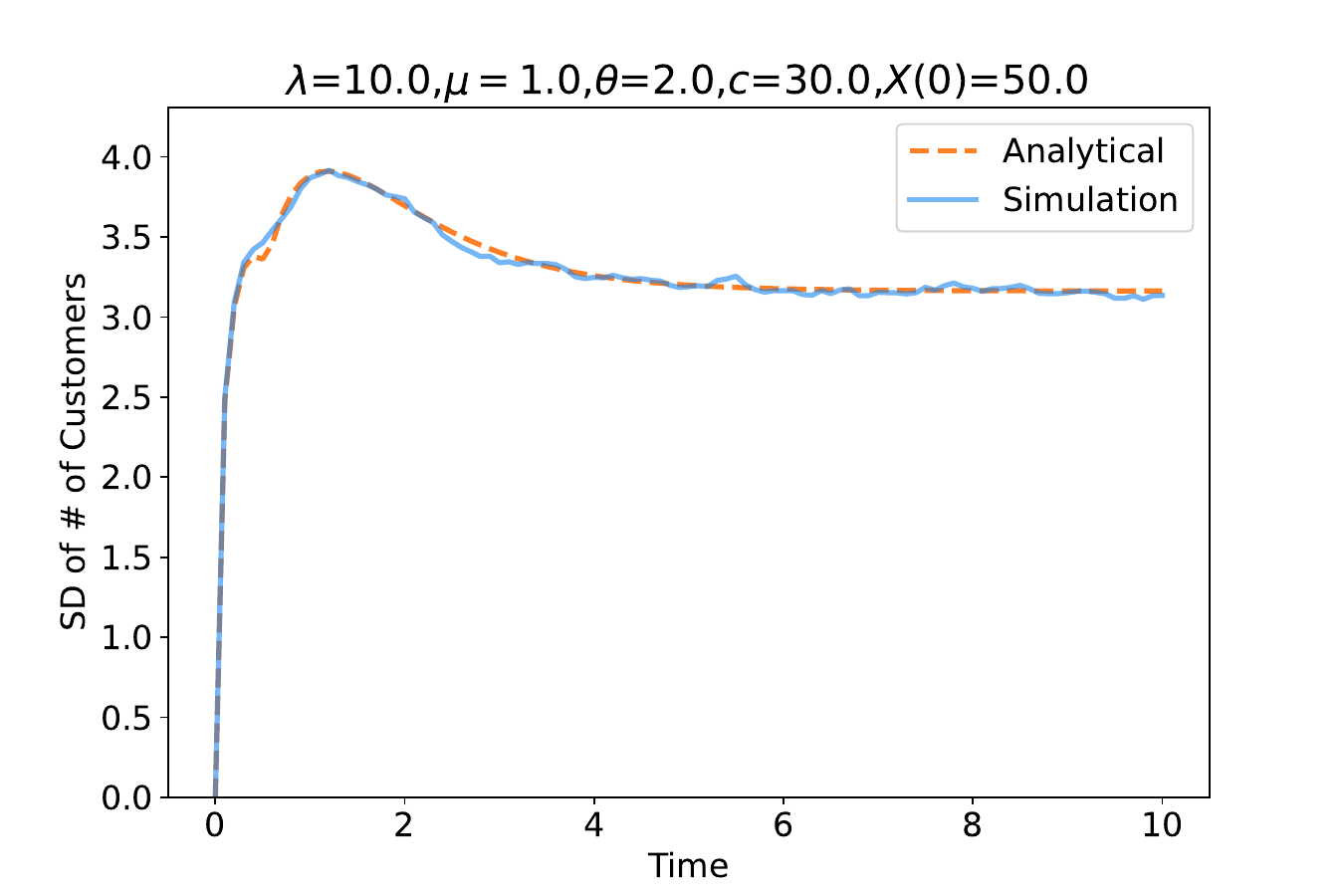}}
\\
\subfloat[]{\includegraphics[scale=.31]{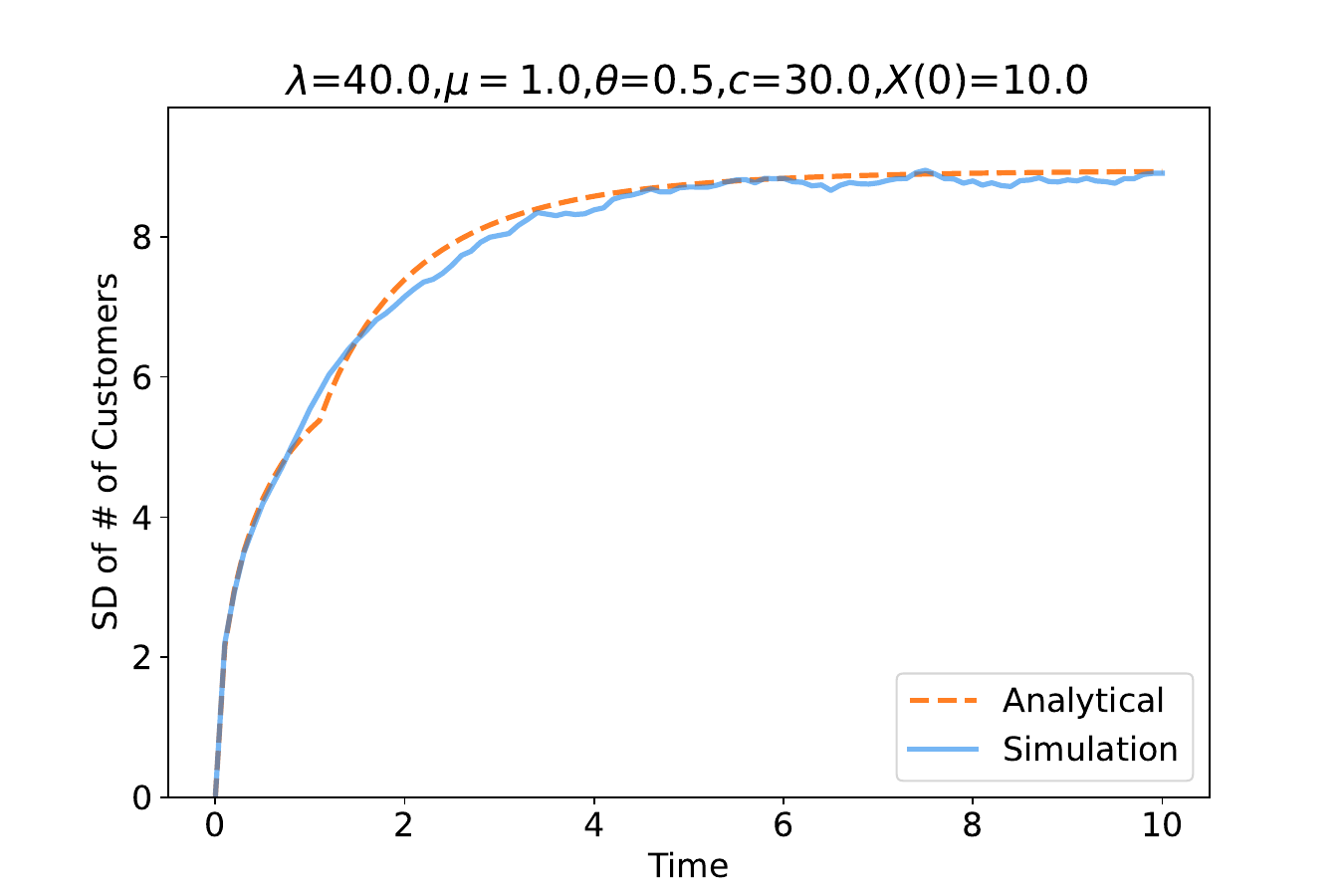}}
\subfloat[]{\includegraphics[scale=.31]{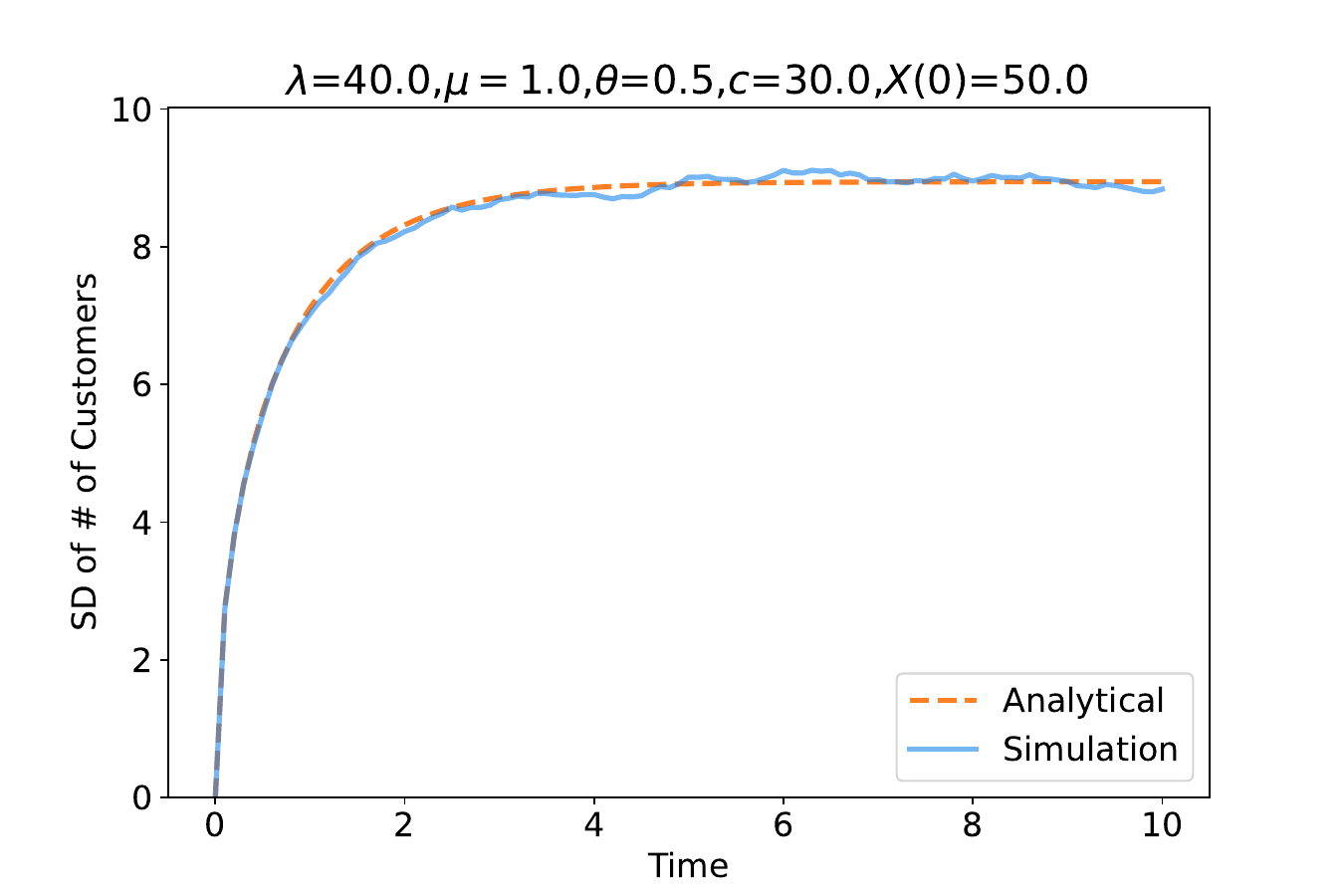}}
\\
\subfloat[]{\includegraphics[scale=.31]{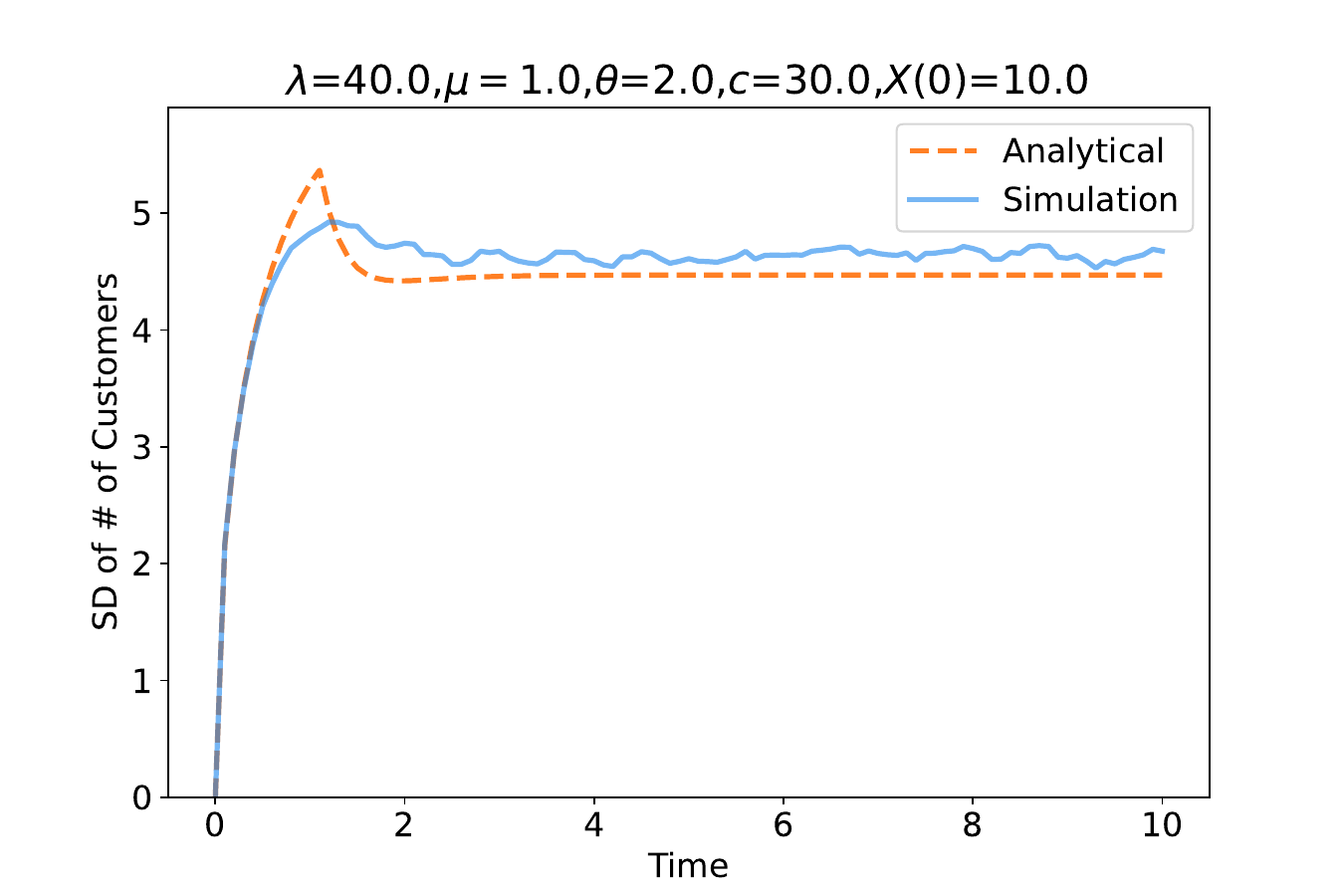}}
\subfloat[]{\includegraphics[scale=.31]{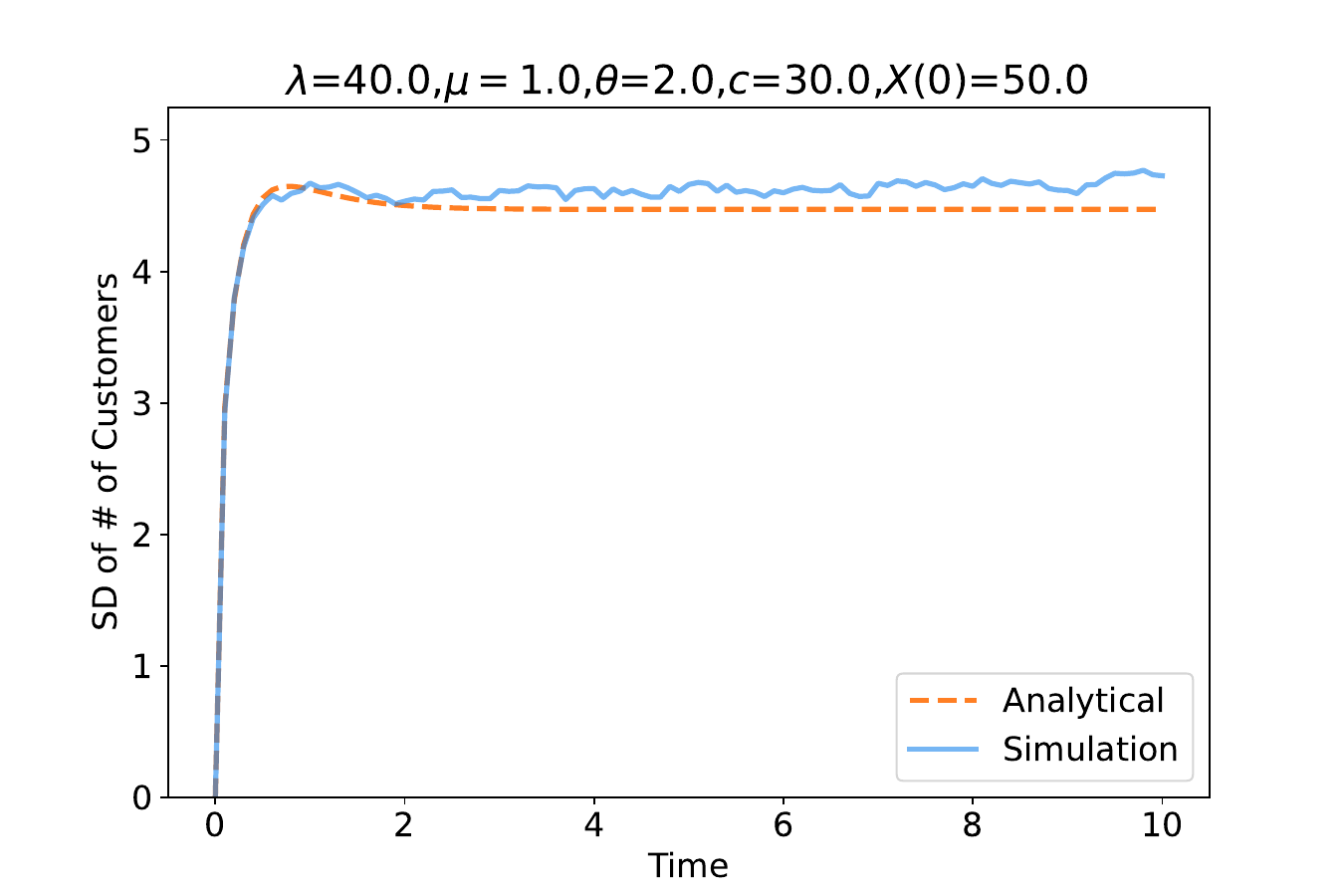}}
 \caption{Standard Deviation Number of Customers (Analytical vs. Simulation). }
\label{Figure_2}
\end{figure}

In Figure \ref{Figure_2}, we plot the standard deviation of the queue length for the parameters outlined in Table \ref{Table_1}. We observe that the approximation provided by the diffusion variance, as given in Equation \eqref{diffvar}, consistently performs well at approximating the corresponding simulated values, for all of the parameter values. Thus, this illustrates the  accuracy of the diffusion variance in effectively approximating the dynamic behavior of the queue length standard deviation.  Moreover, with good approximations for the standard deviation, it  becomes feasible to construct prediction intervals for the queue length. Such prediction intervals serve as valuable tools for assessing the range within which the actual queue length is likely to fall, providing a measure of confidence in our approximations.  

\begin{figure}[!htbp]
\centering
\subfloat[]{\includegraphics[scale=.31]{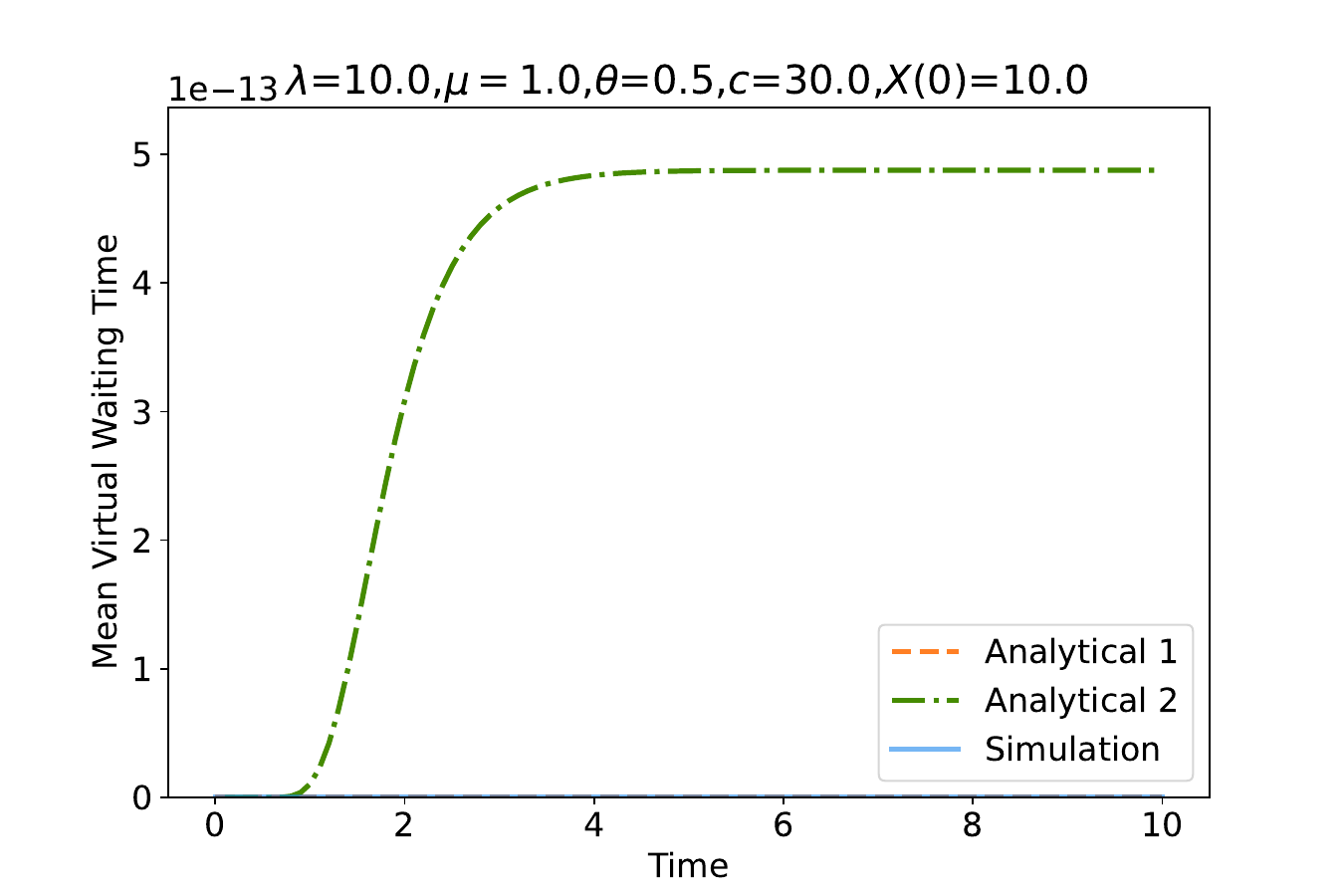}}
\subfloat[]{\includegraphics[scale=.31]{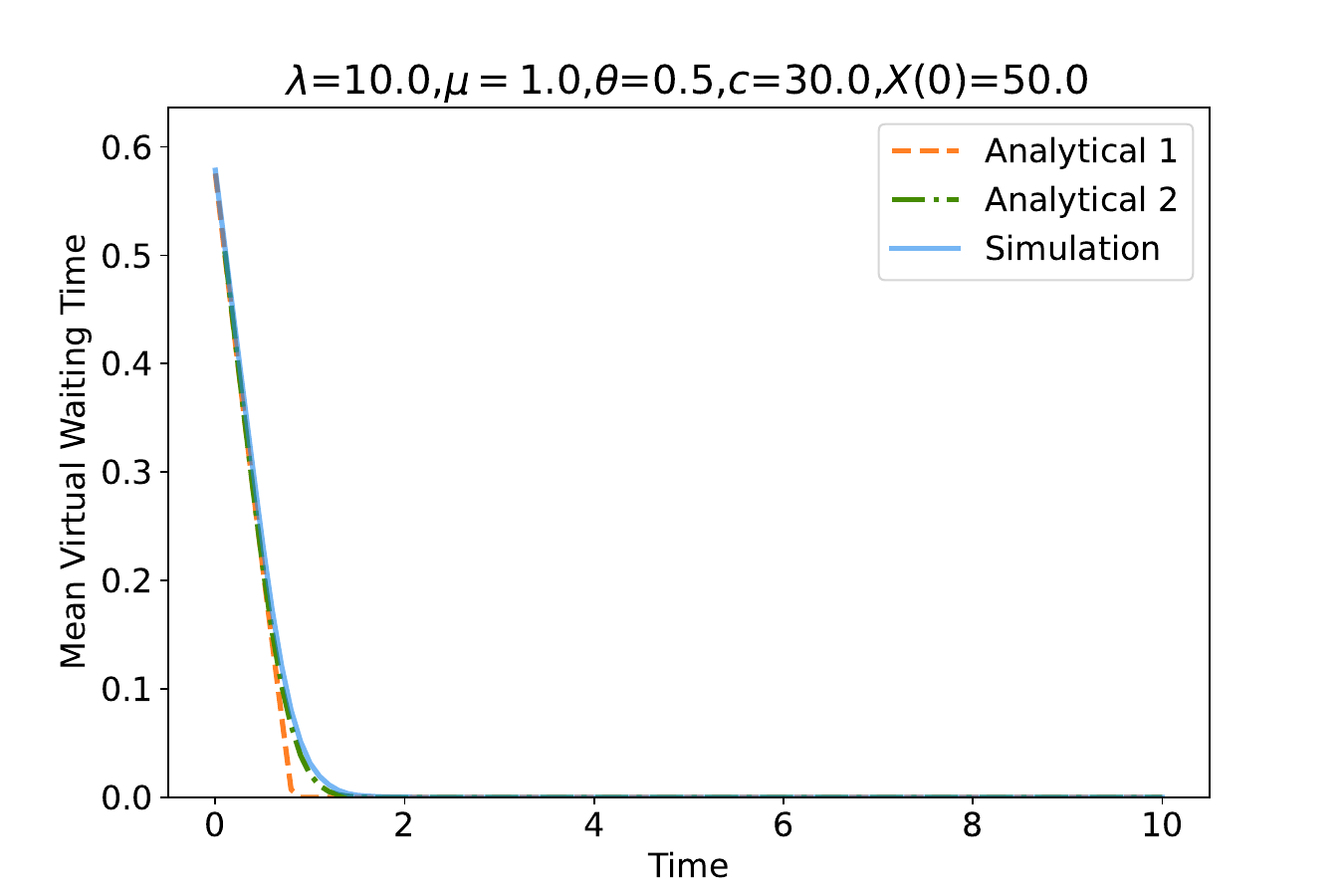}}
\\
\subfloat[]{\includegraphics[scale=.31]{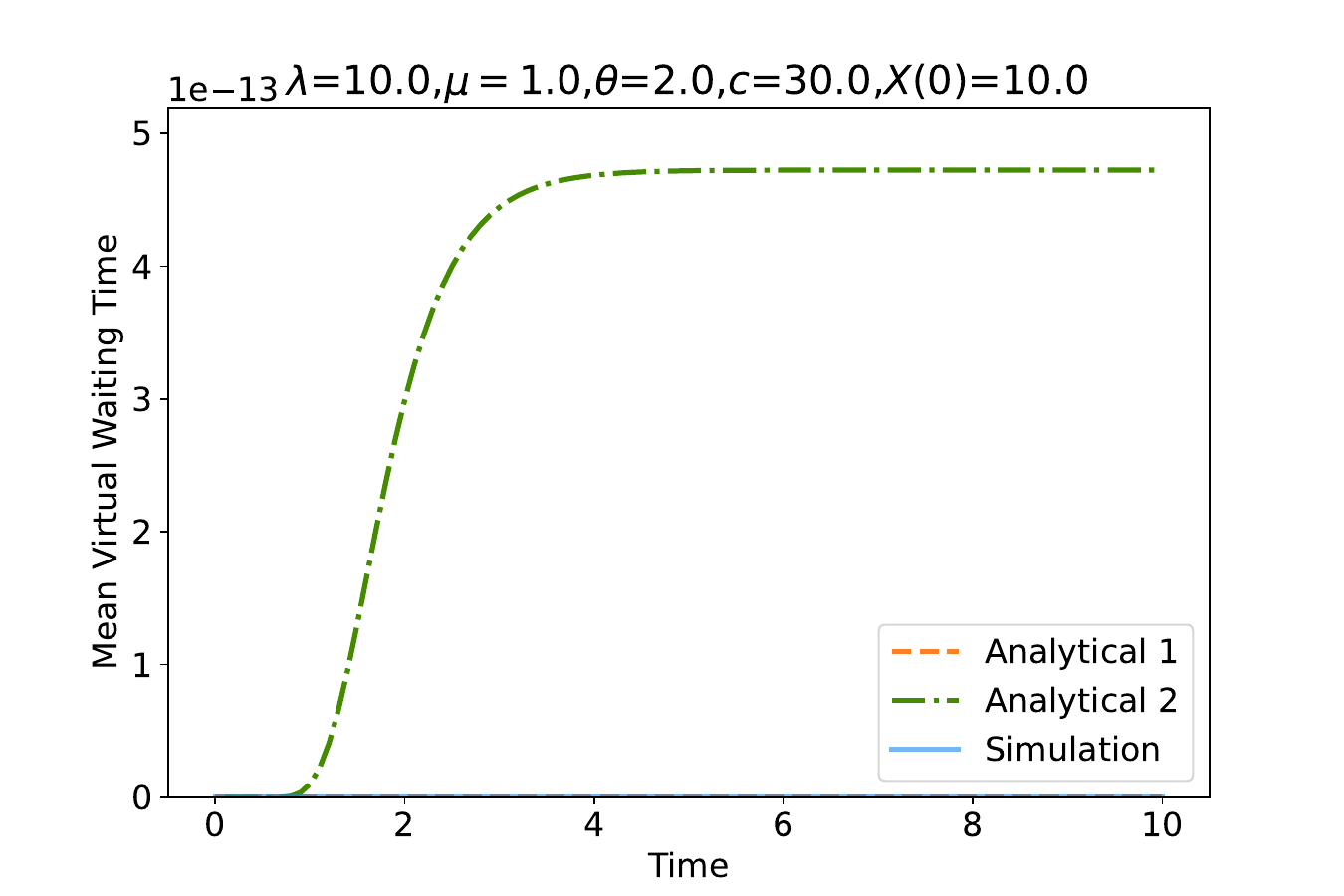}}
\subfloat[]{\includegraphics[scale=.31]{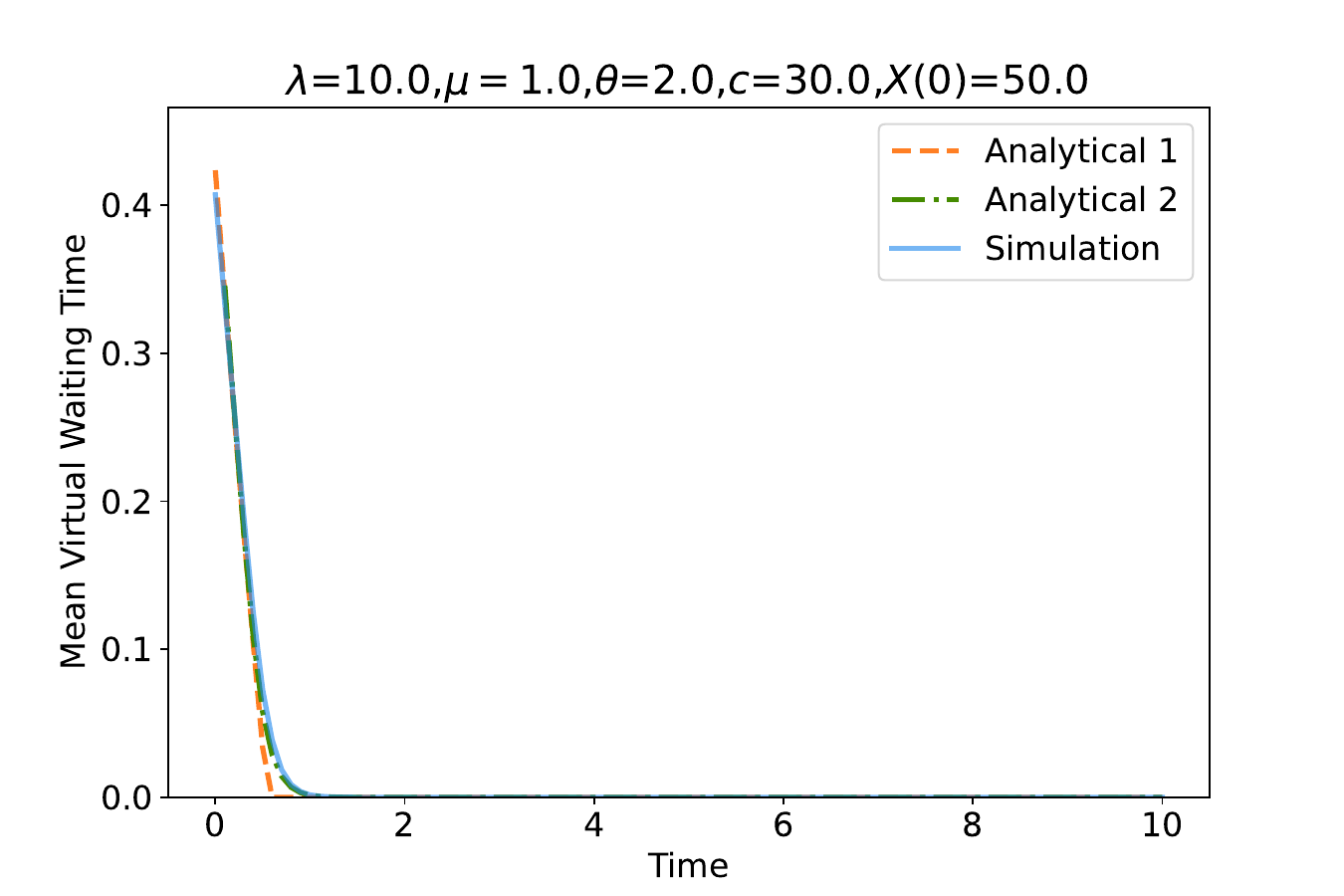}}
\\
\subfloat[]{\includegraphics[scale=.31]{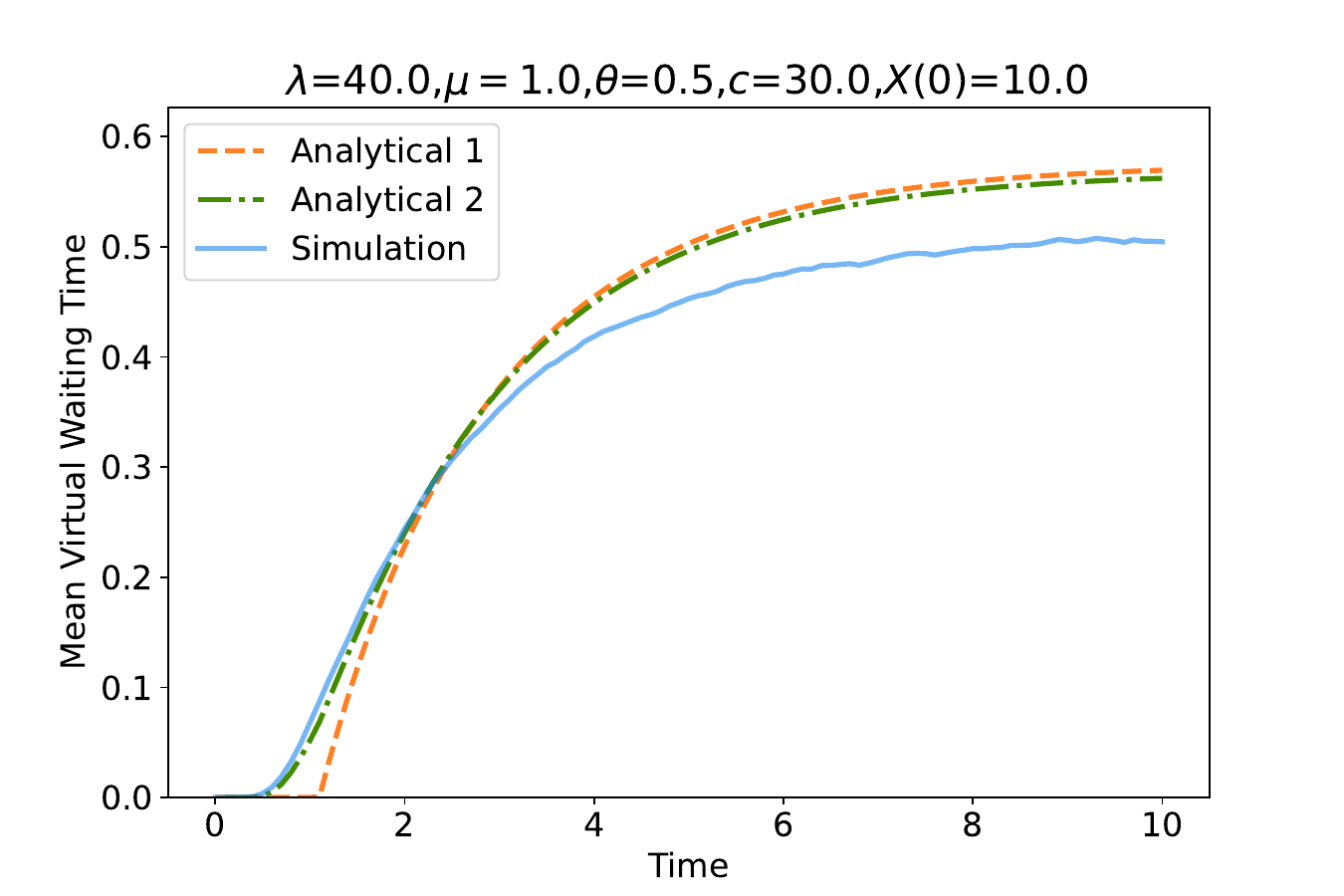}}
\subfloat[]{\includegraphics[scale=.31]{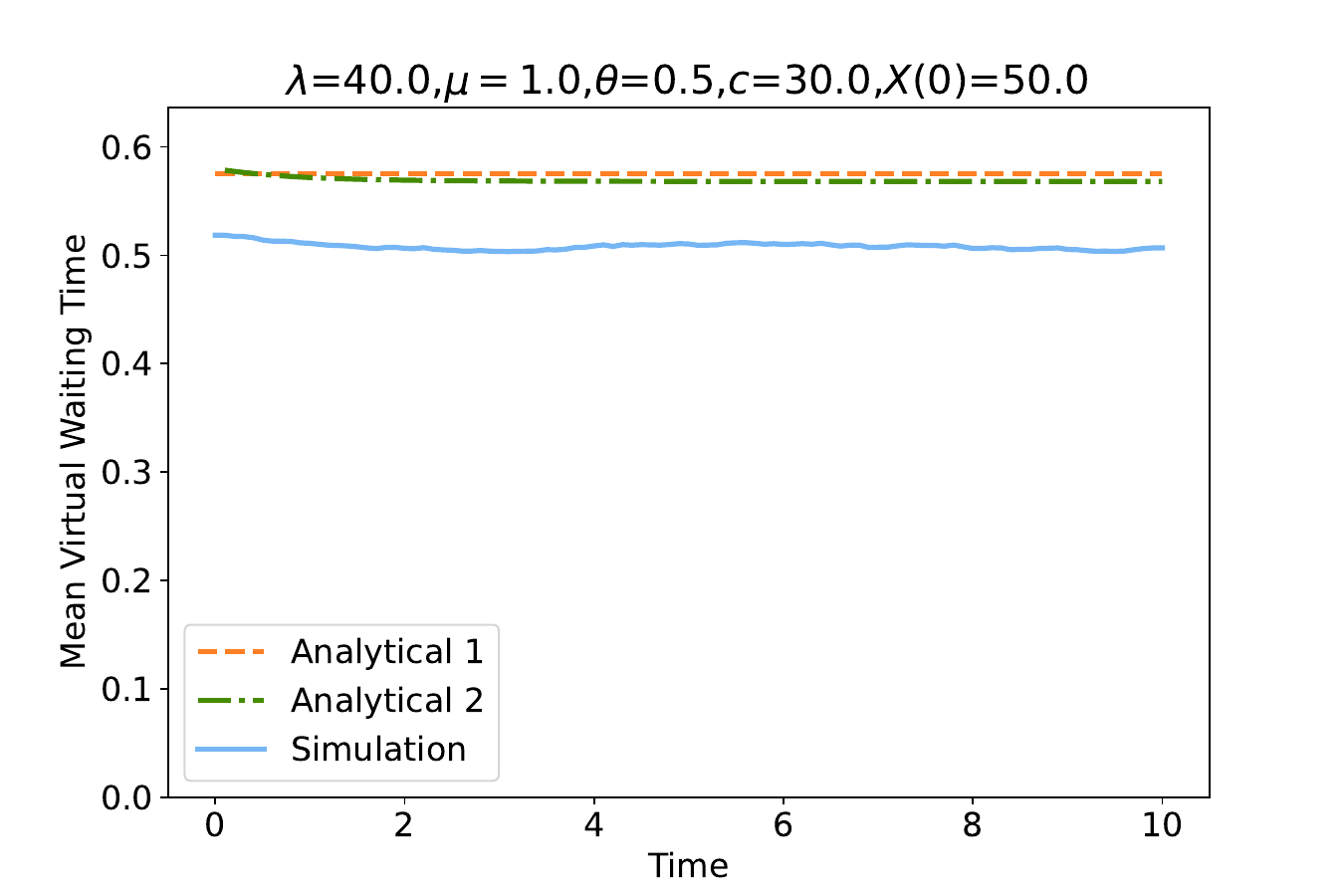}}
\\
\subfloat[]{\includegraphics[scale=.31]{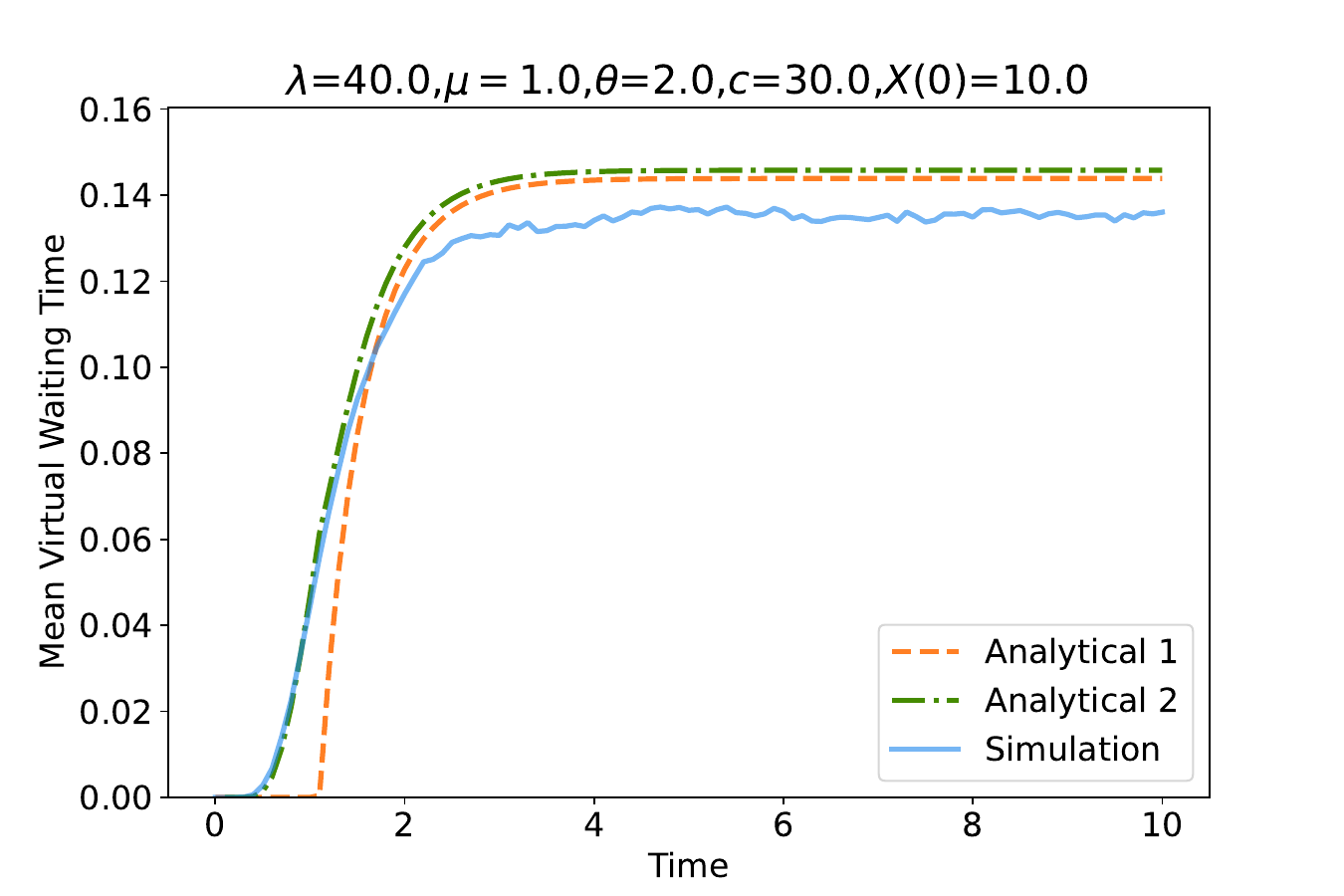}}
\subfloat[]{\includegraphics[scale=.31]{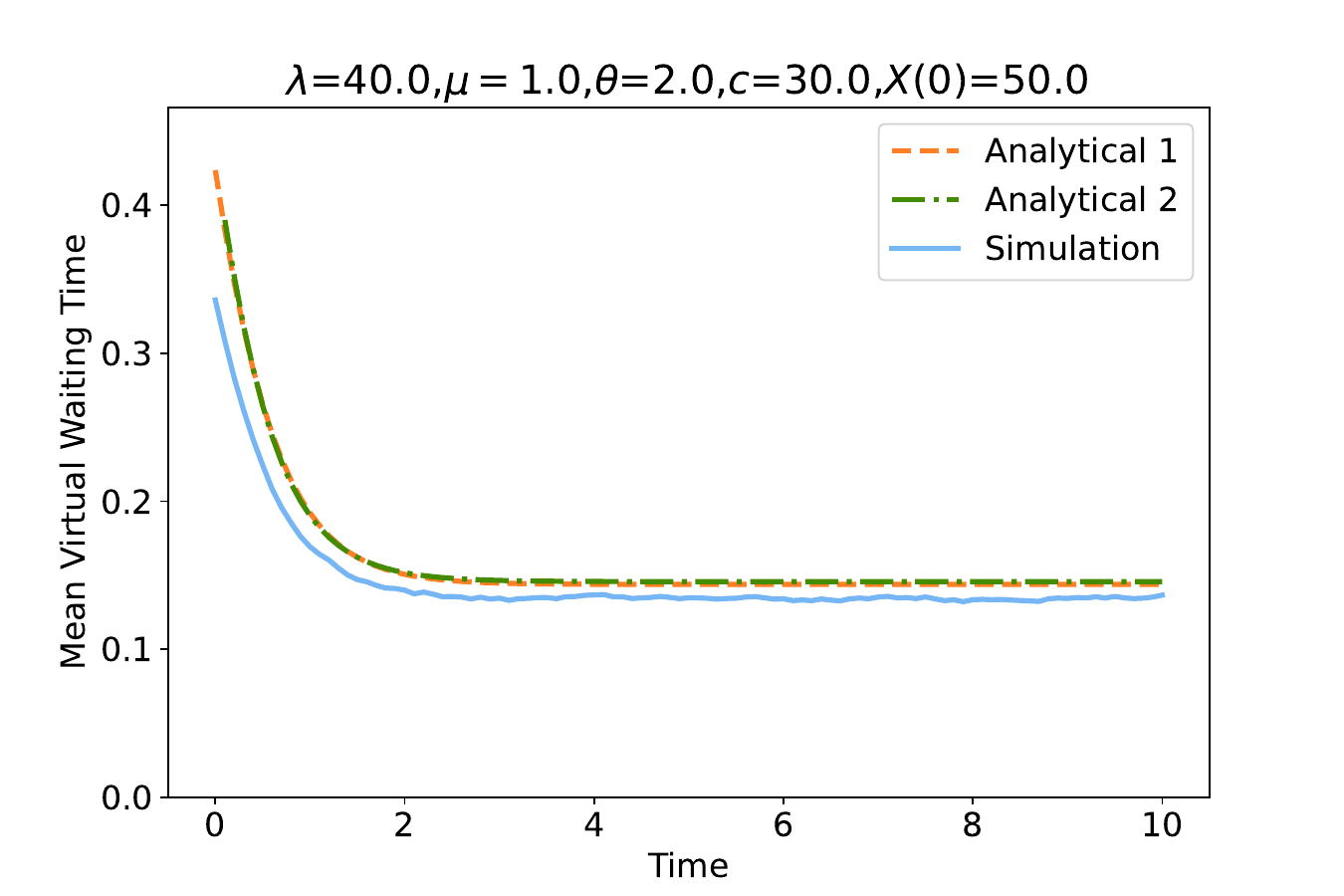}}
 \caption{Mean Virtual Waiting Time (Analytical vs Simulation). }
\label{Figure_3}
\end{figure}

In Figure \ref{Figure_3}, we plot the simulated mean virtual waiting time with two different approximations. The first approximation, denoted as ``analytical 1" in the figures, is derived from the fluid limit, as outlined in Theorem \ref{v_wait_time} and is given in Equation \eqref{fluid_wait}. The second approximation, denoted as ``analytical 2" in the figure, is given by Equation \eqref{mean_wait_digamma2}, which depends on the digamma function and the fluid queue length.  We observe that across all parameter settings provided in Table \ref{Table_1}, the approximation employing the digamma function consistently outperforms the fluid-based version in estimating the virtual waiting time. 

\begin{figure}[!htbp]
\centering
\subfloat[]{\includegraphics[scale=.31]{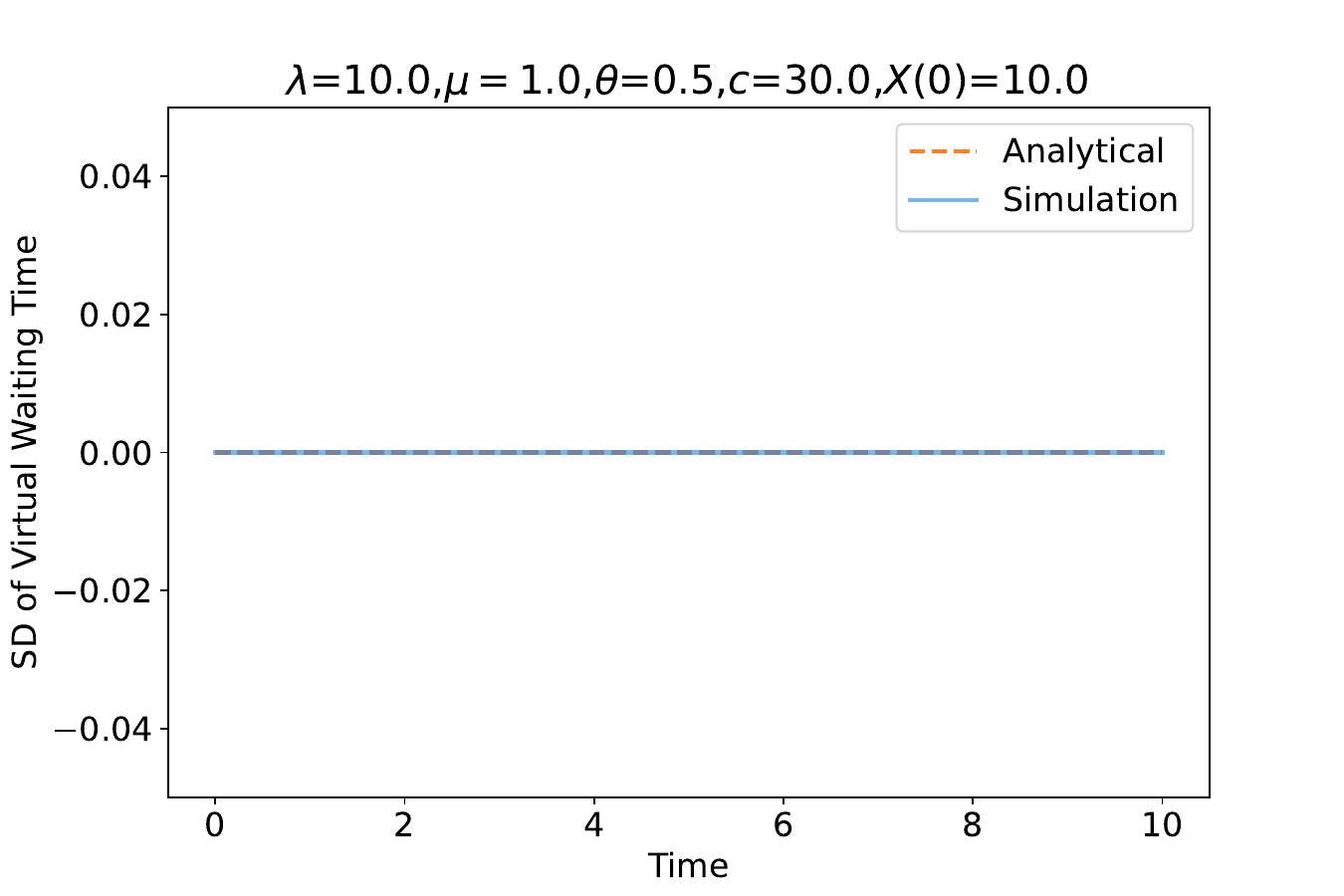}}
\subfloat[]{\includegraphics[scale=.31]{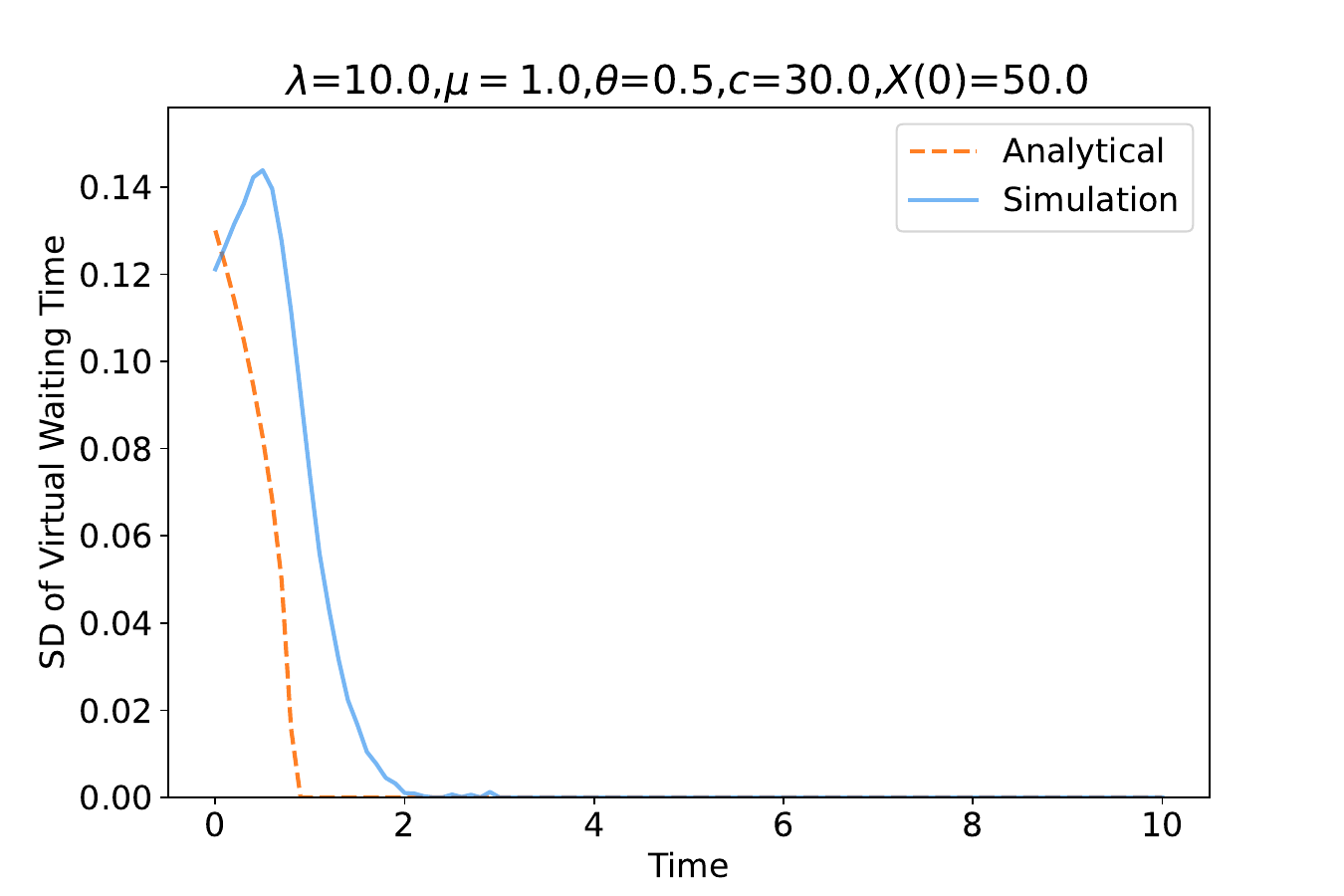}}
\\
\subfloat[]{\includegraphics[scale=.31]{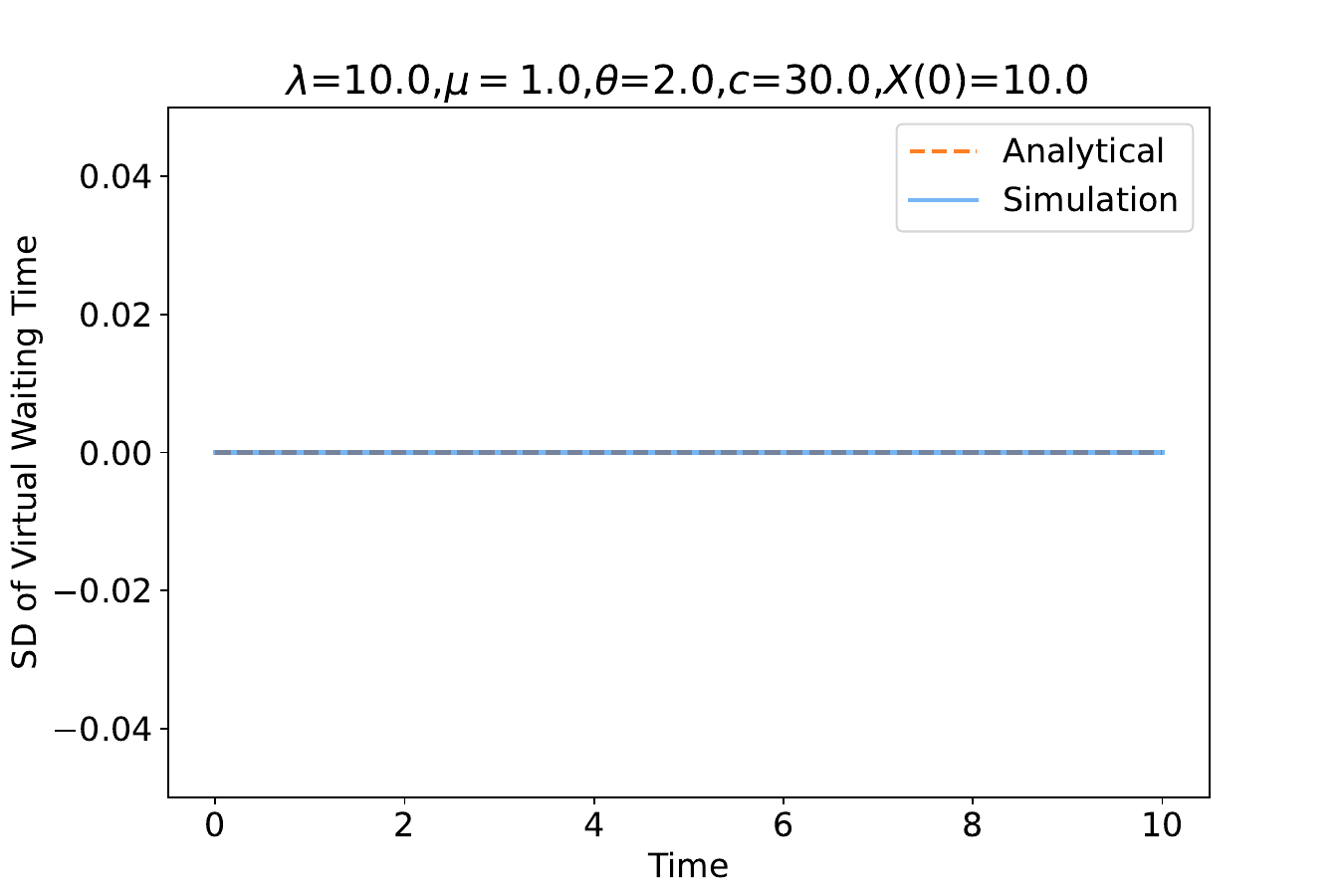}}
\subfloat[]{\includegraphics[scale=.31]{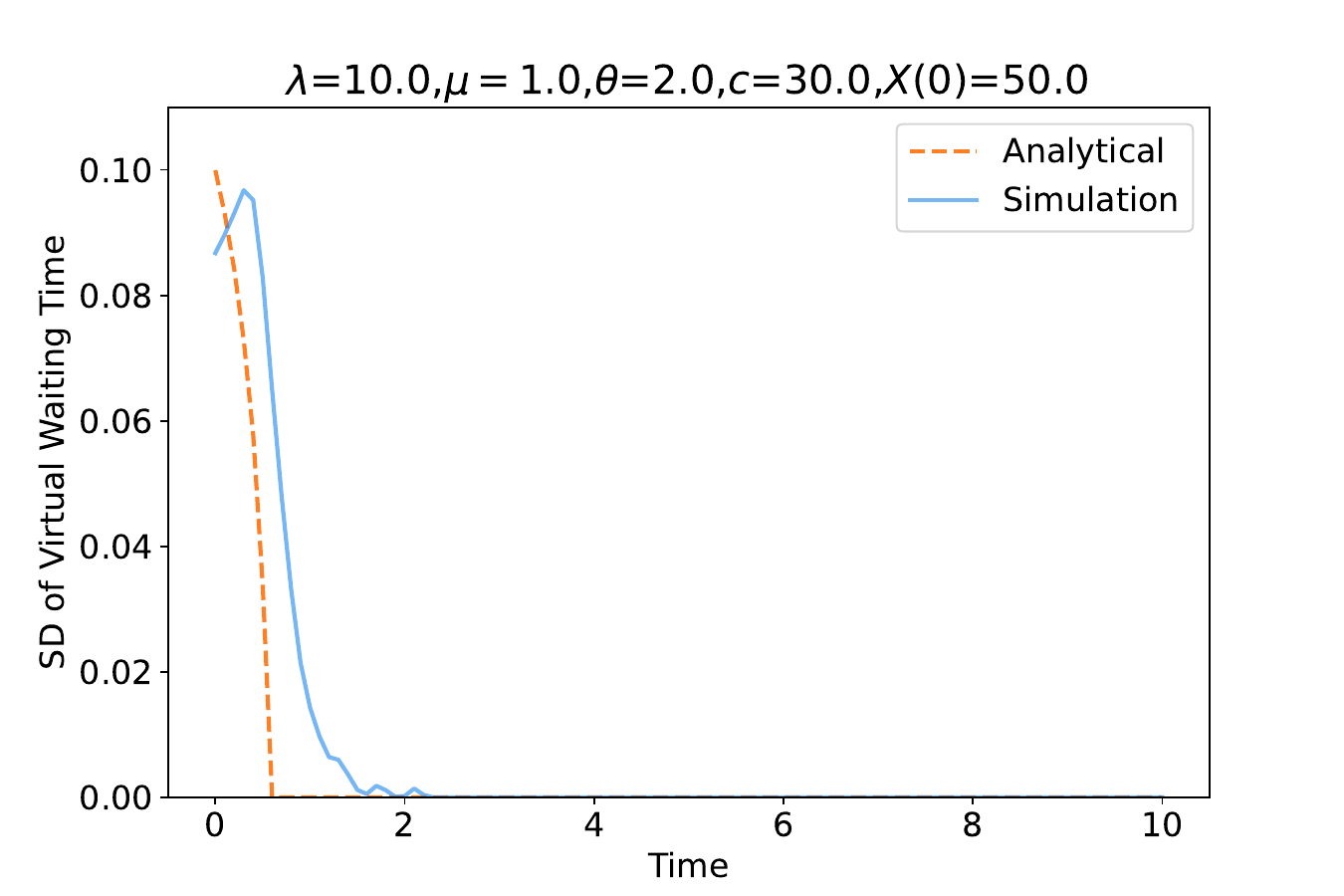}}
\\
\subfloat[]{\includegraphics[scale=.31]{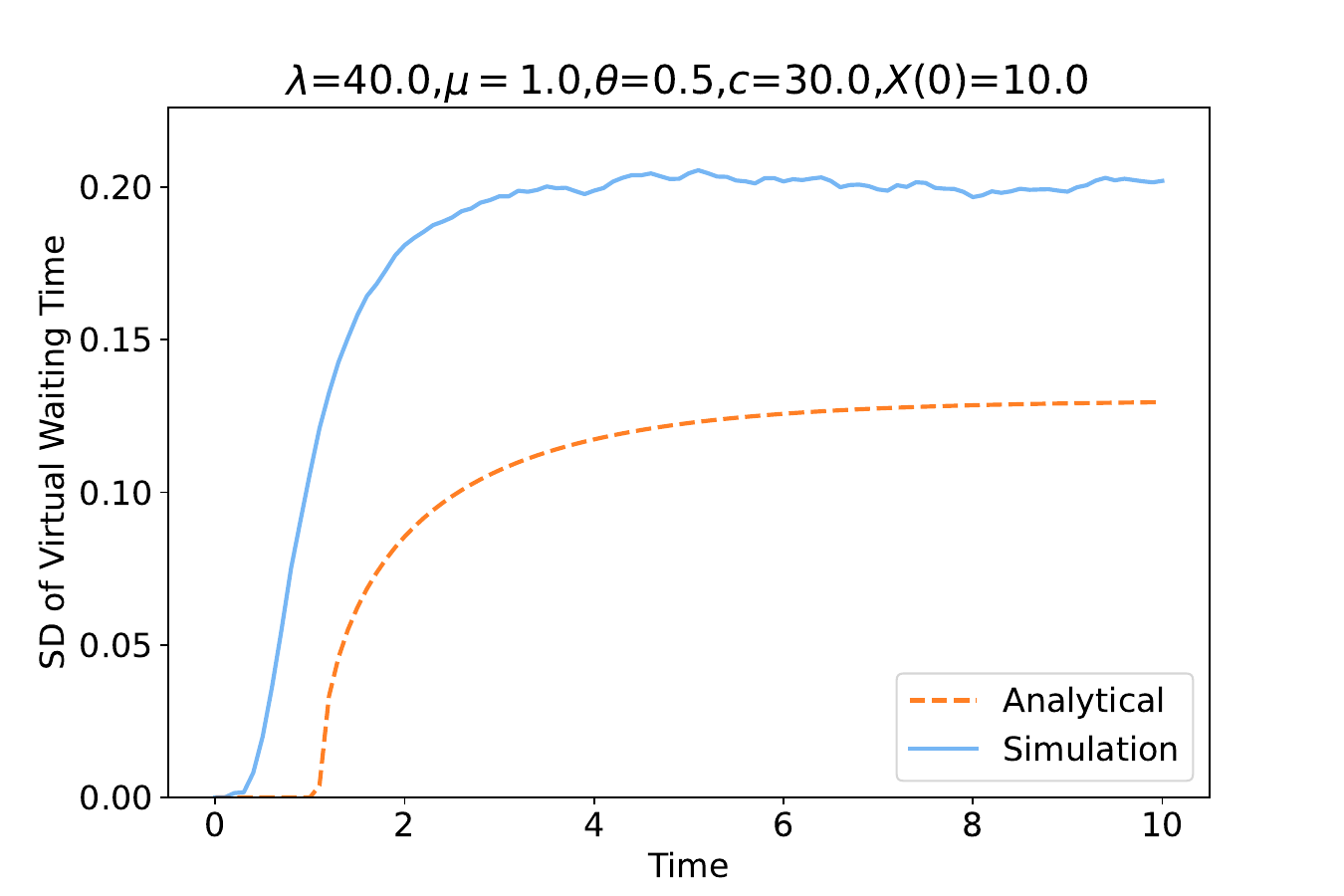}}
\subfloat[]{\includegraphics[scale=.31]{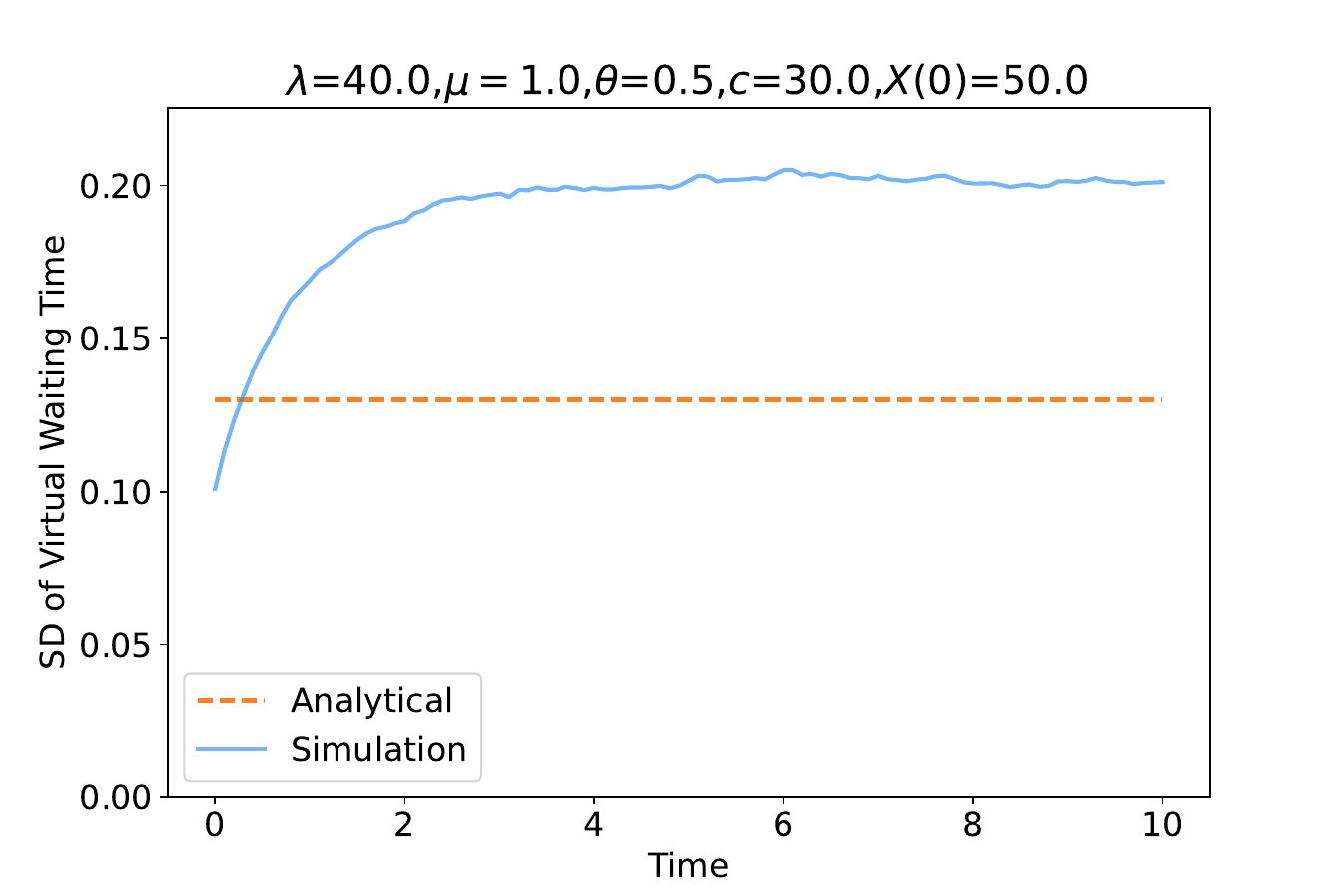}}
\\
\subfloat[]{\includegraphics[scale=.31]{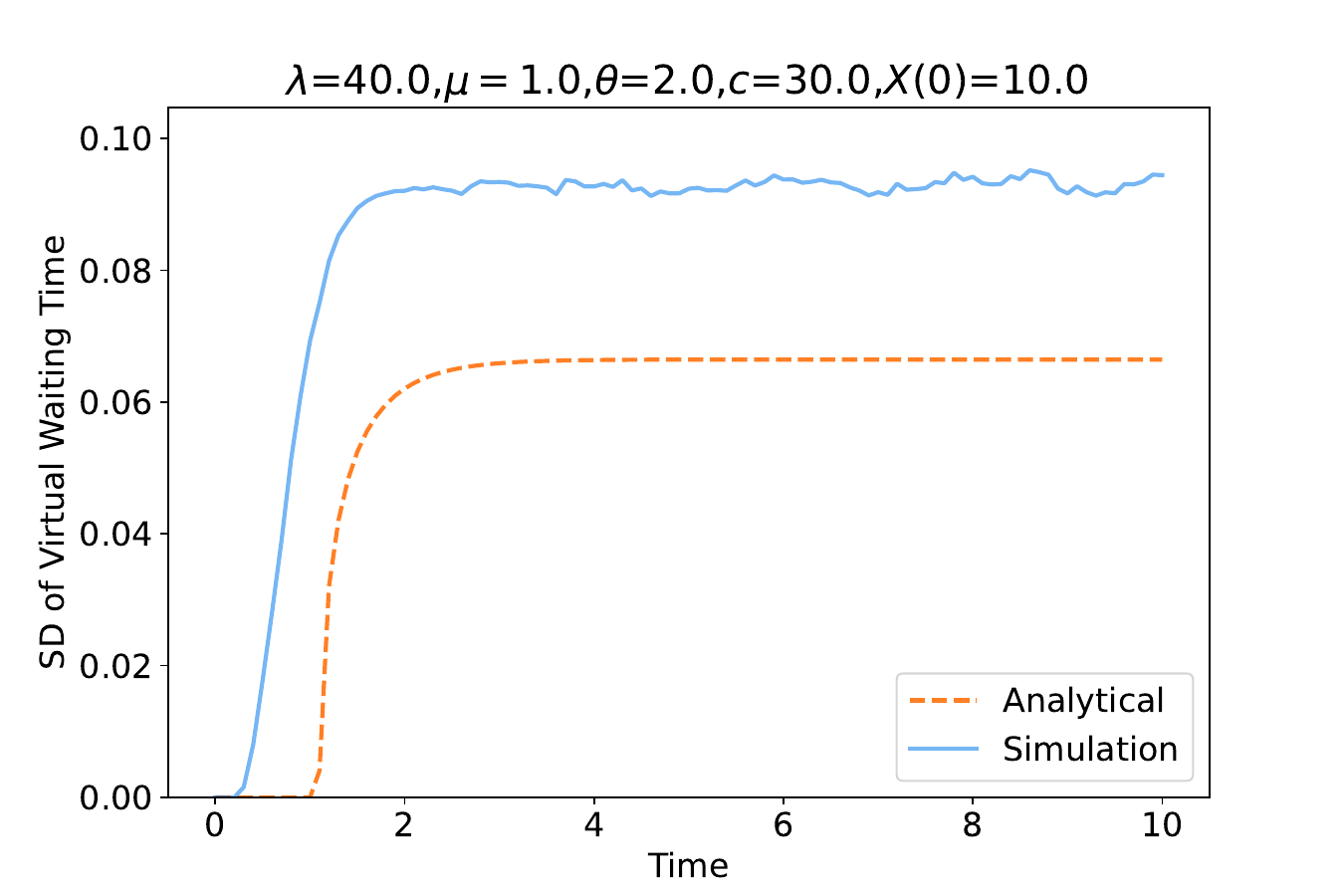}}
\subfloat[]{\includegraphics[scale=.31]{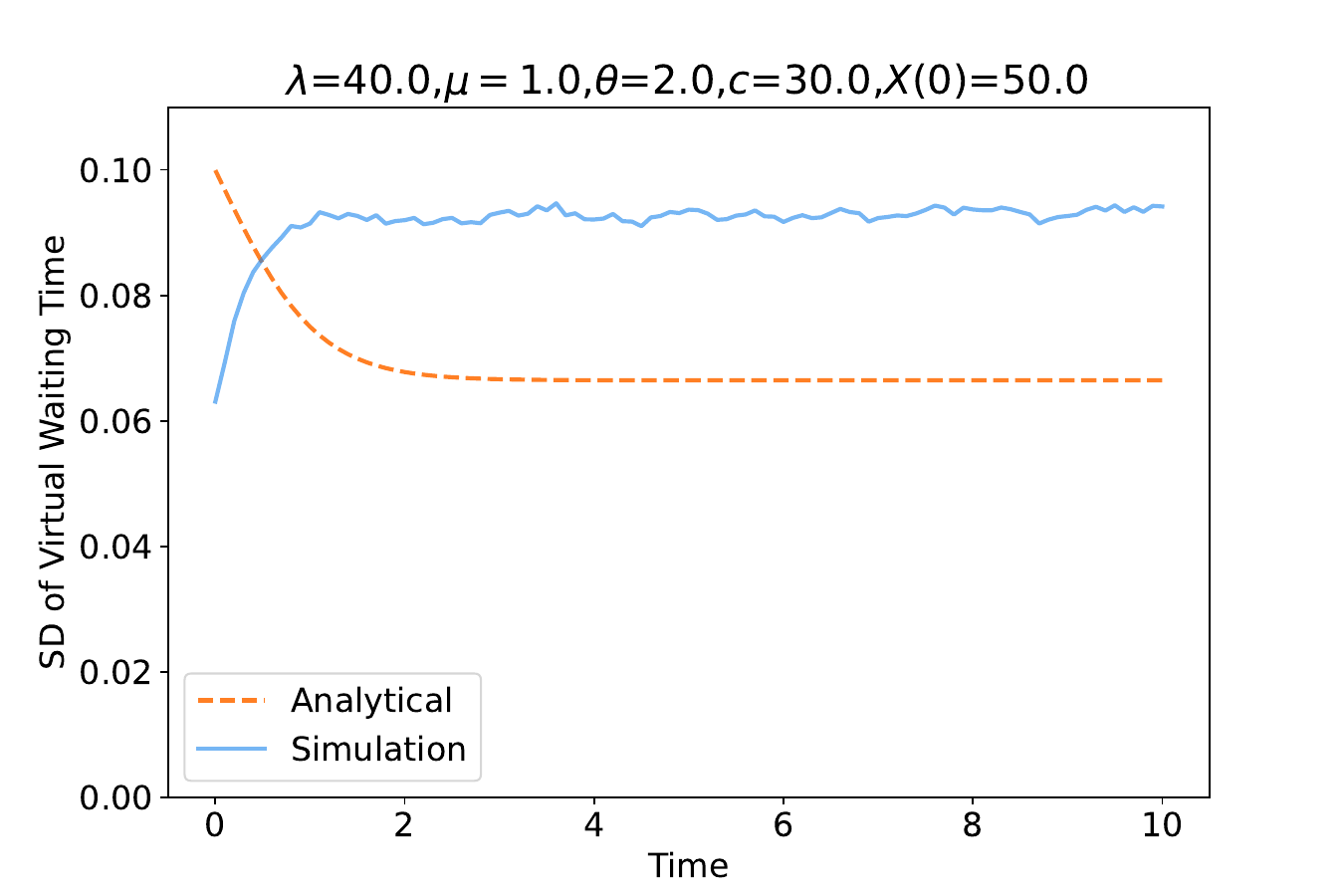}}
 \caption{Standard Deviation of Virtual Waiting Time (Analytical vs. Simulation). }
\label{Figure_4}
\end{figure}

In Figure \ref{Figure_4}, we plot the simulated standard deviation of the virtual waiting time, utilizing an approximation derived from the trigamma function as given in Equation \eqref{var_wait_trigamma2}. The derivation of this approximation involved a Taylor expansion of the variance and the utilization of conditional mean and variance formulas for the wait time, as outlined in Equations \eqref{cond_mean_wait} and \eqref{cond_var_wait}. Our observations reveal that regardless of the parameter settings examined in Table \ref{Table_1}, the trigamma function consistently demonstrates excellent performance in approximating the standard deviation of the virtual waiting time. This consistency and accuracy highlight the reliability and effectiveness of the trigamma-based approximation method.

\begin{figure}[!htbp]
\centering
\subfloat[]{\includegraphics[scale=.31]{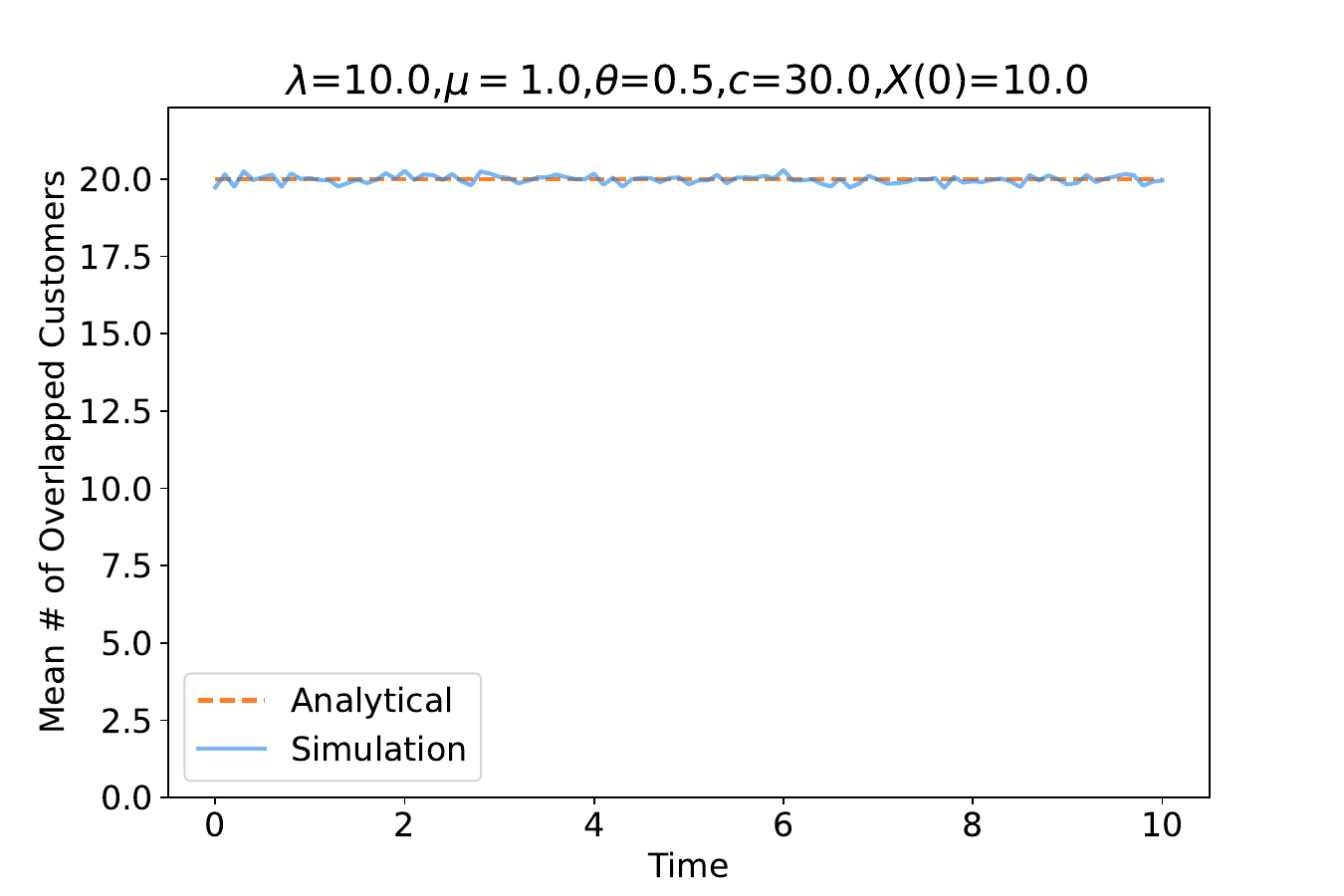}}
\subfloat[]{\includegraphics[scale=.31]{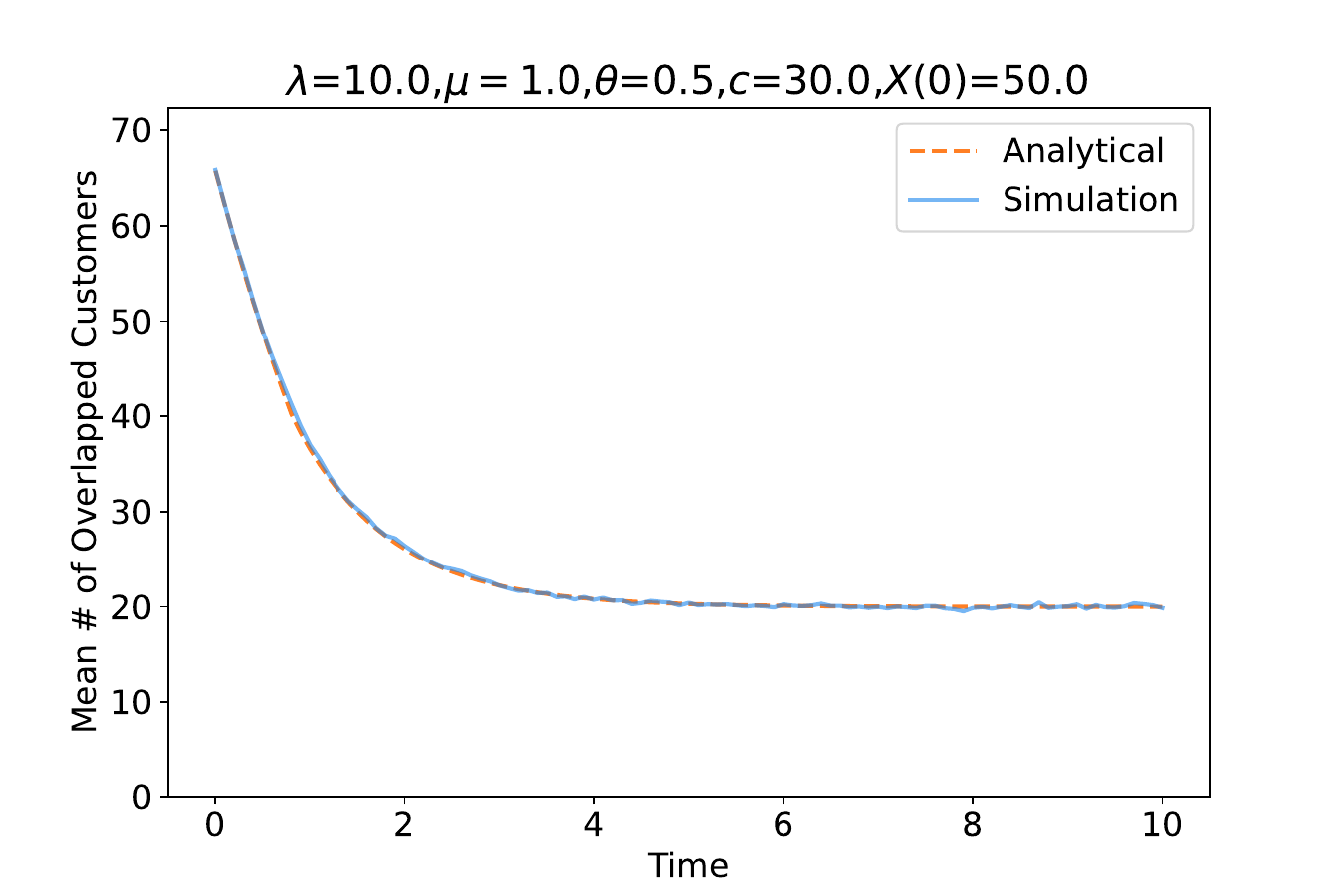}}
\\
\subfloat[]{\includegraphics[scale=.31]{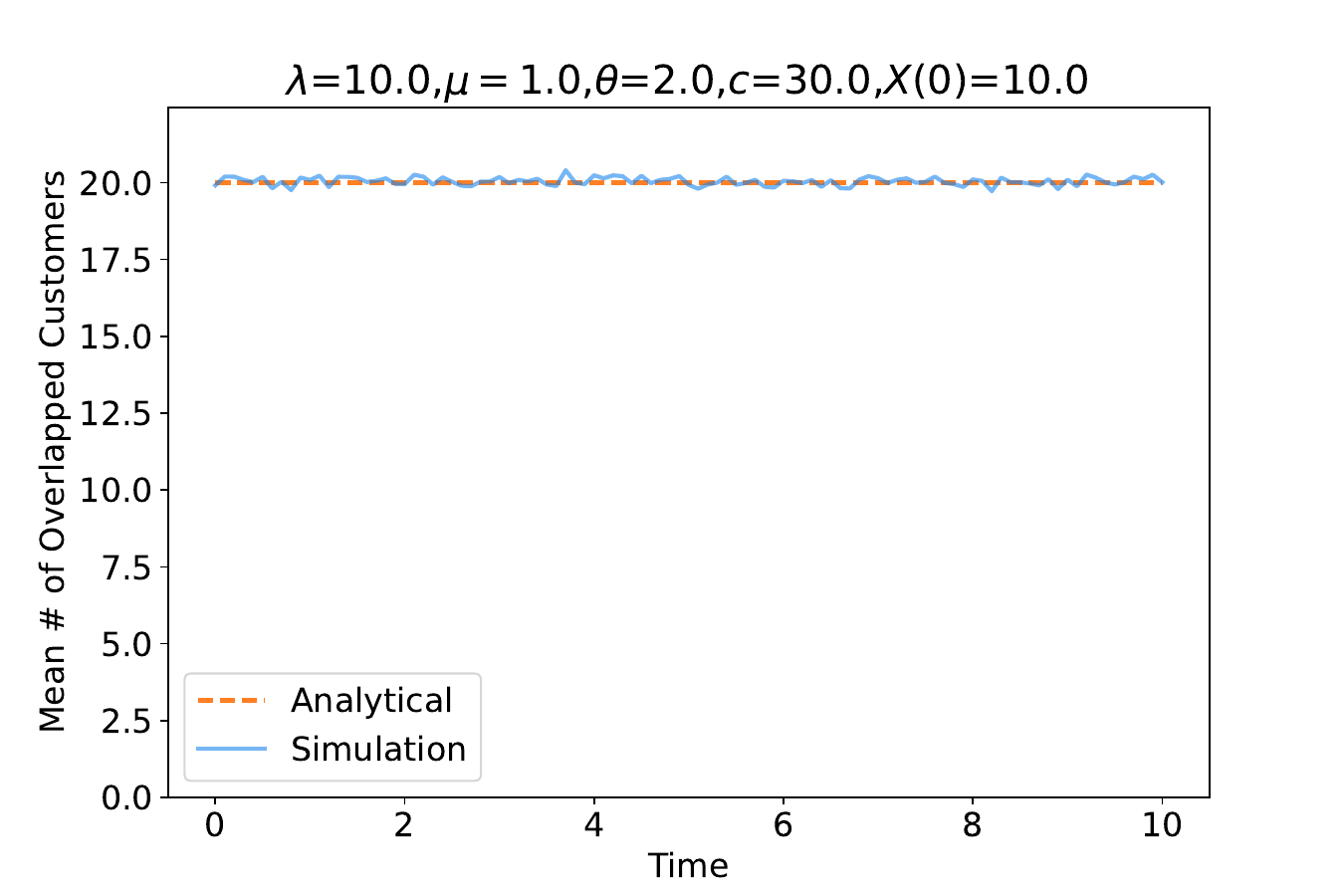}}
\subfloat[]{\includegraphics[scale=.31]{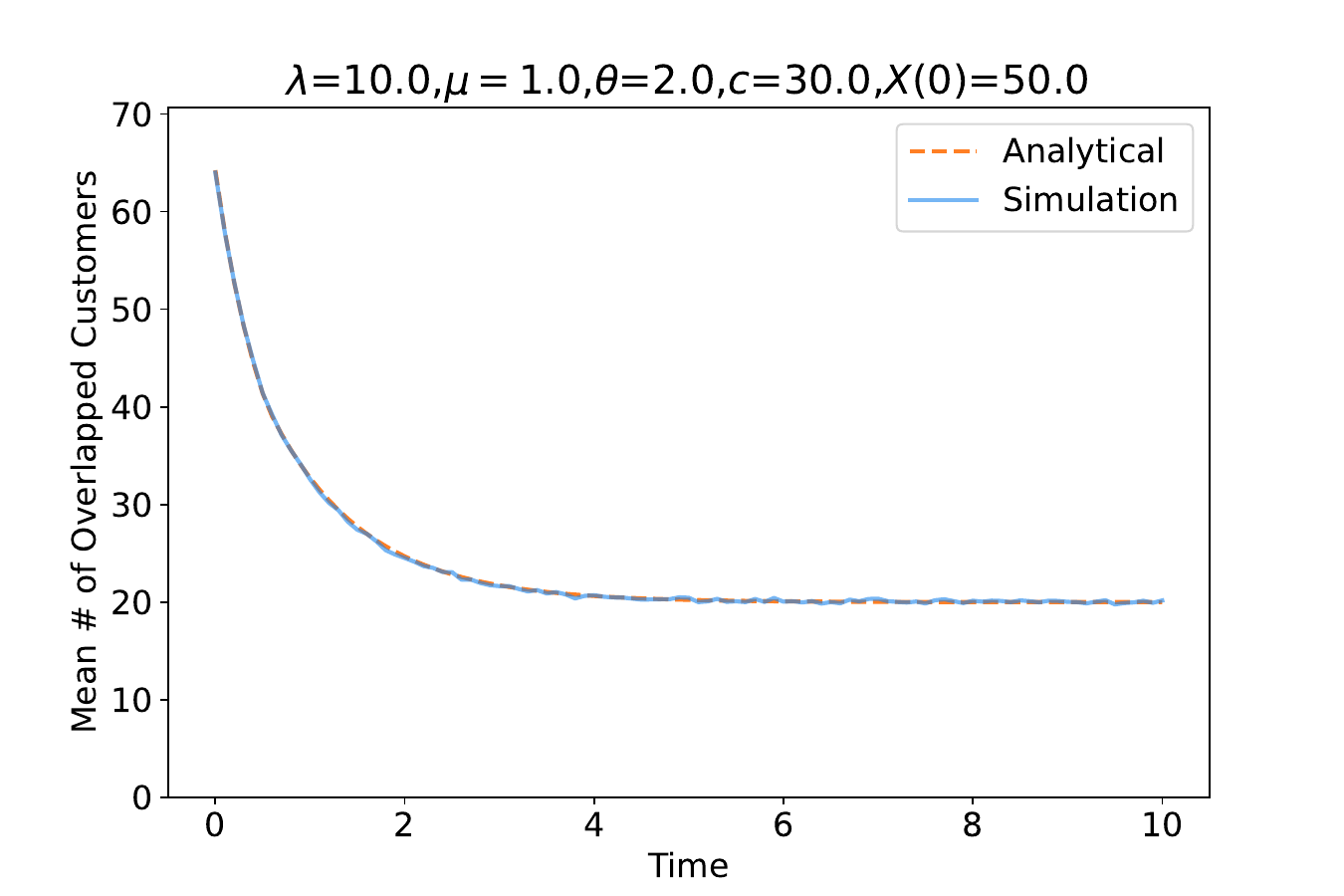}}
\\
\subfloat[]{\includegraphics[scale=.31]{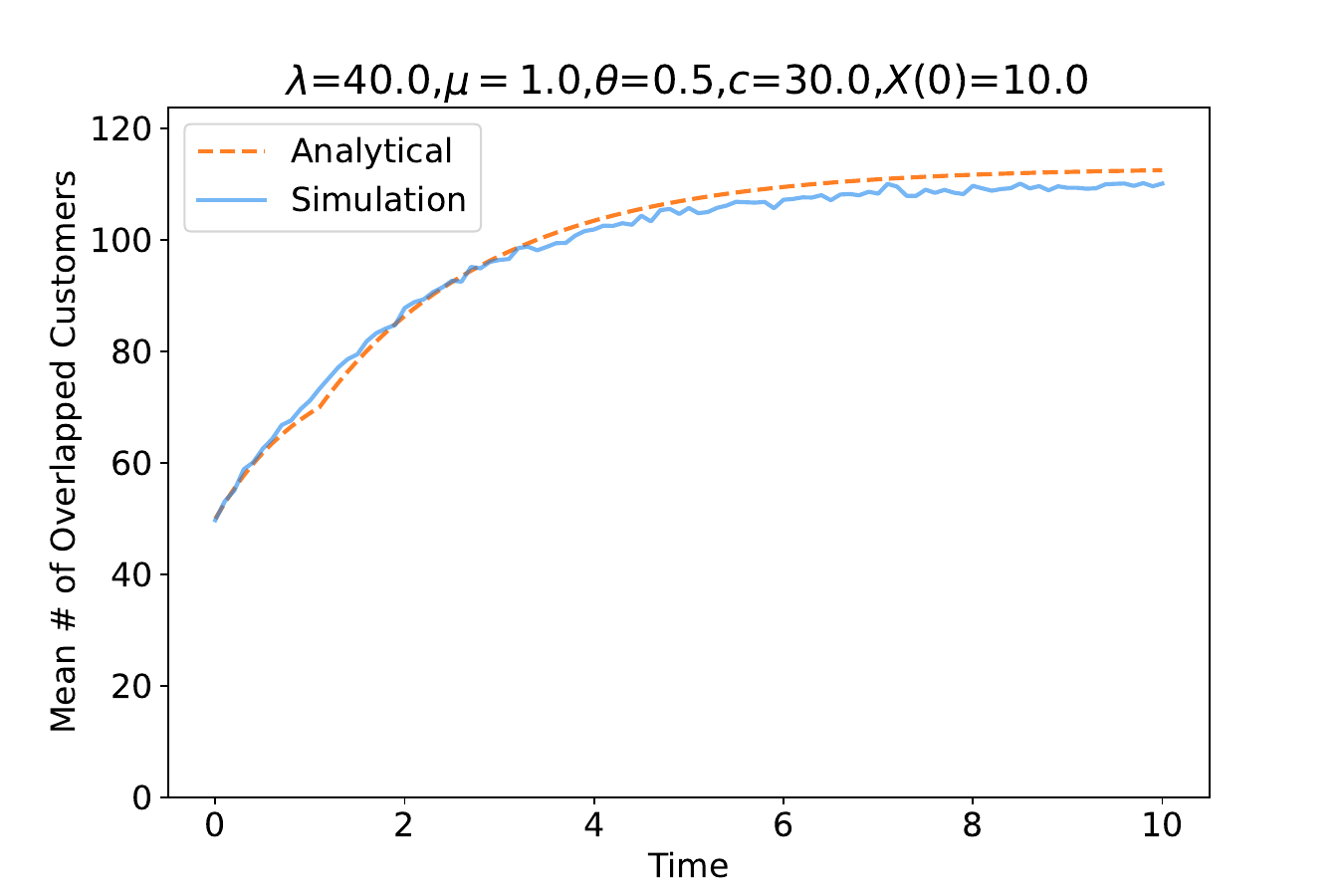}}
\subfloat[]{\includegraphics[scale=.31]{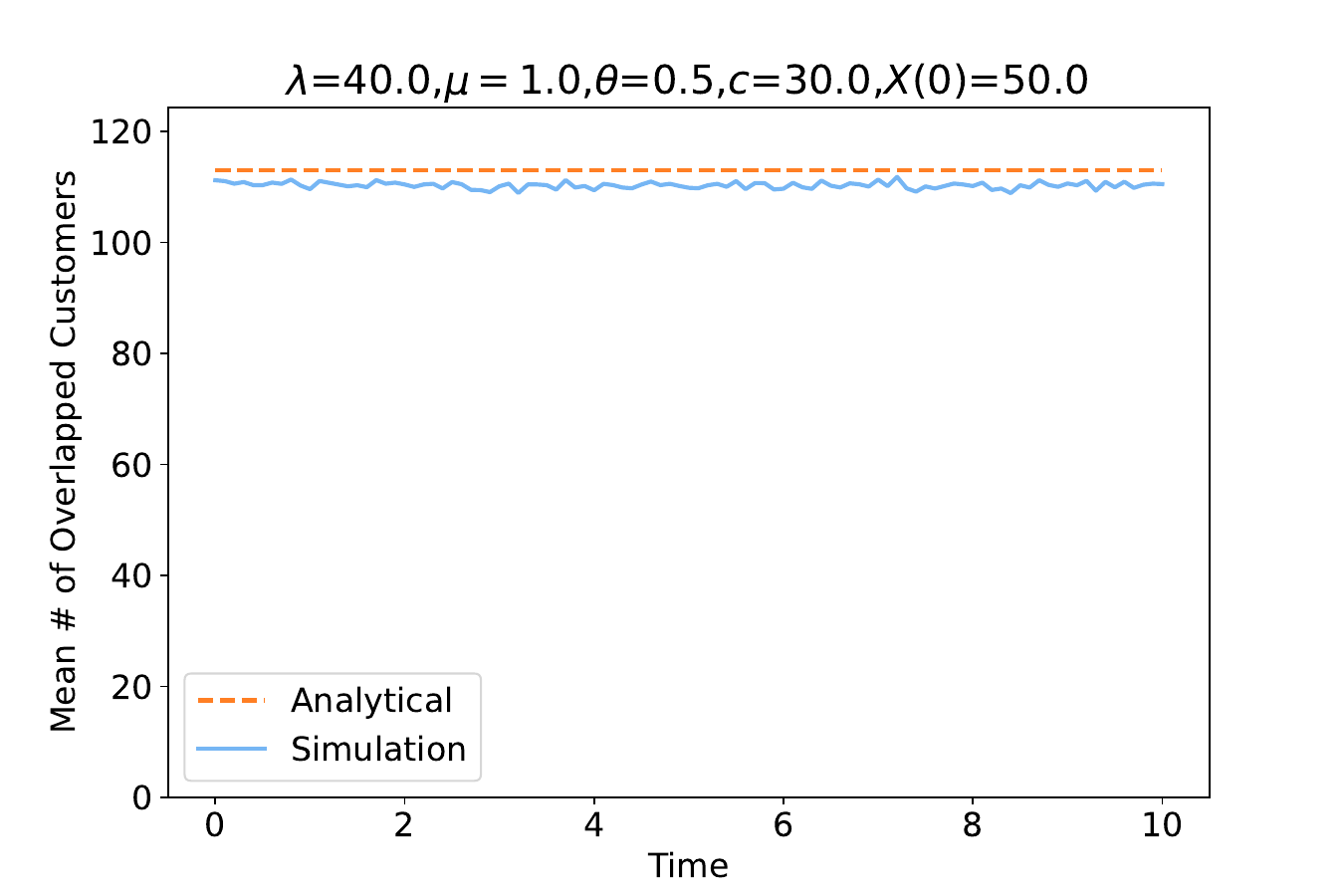}}
\\
\subfloat[]{\includegraphics[scale=.31]{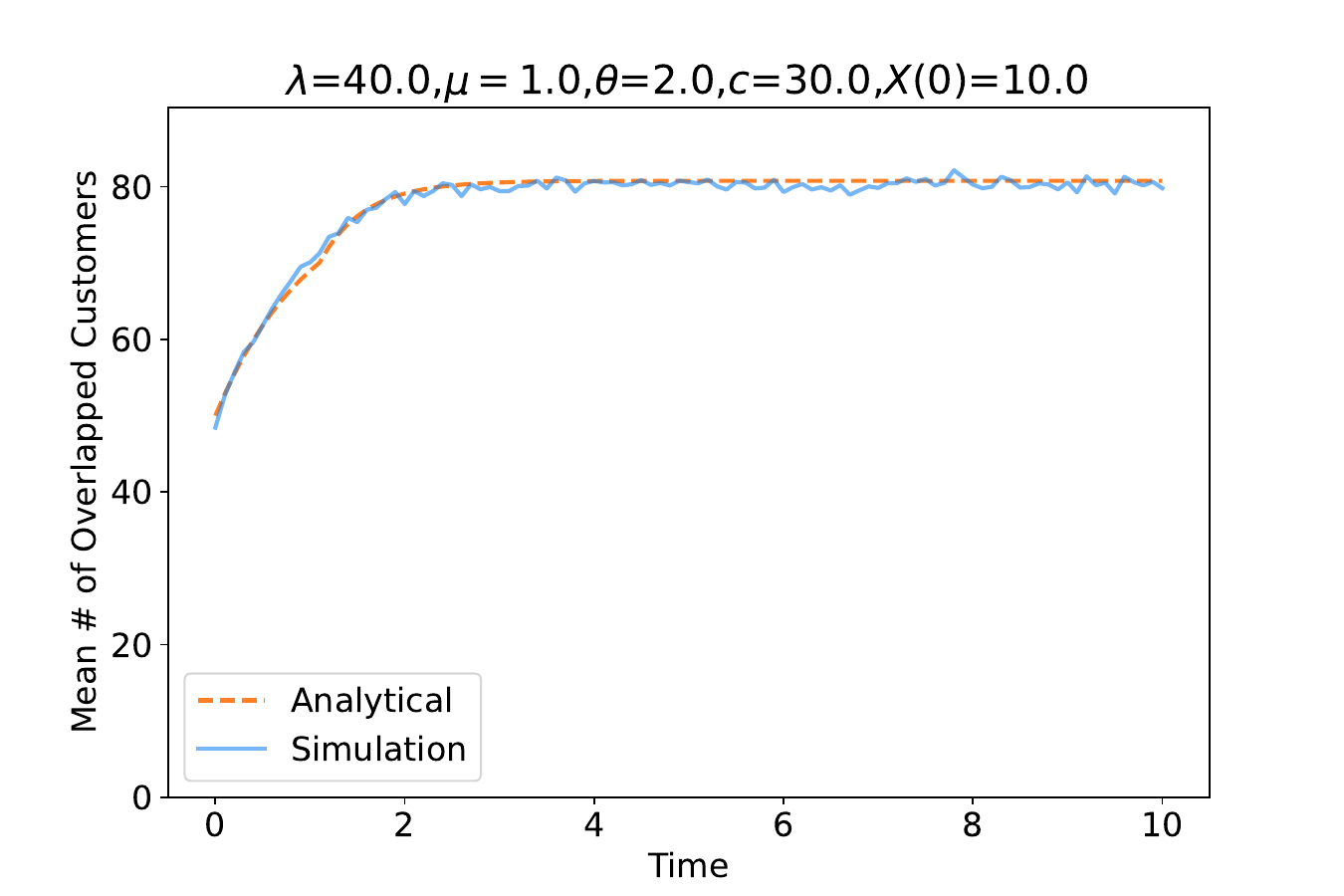}}
\subfloat[]{\includegraphics[scale=.31]{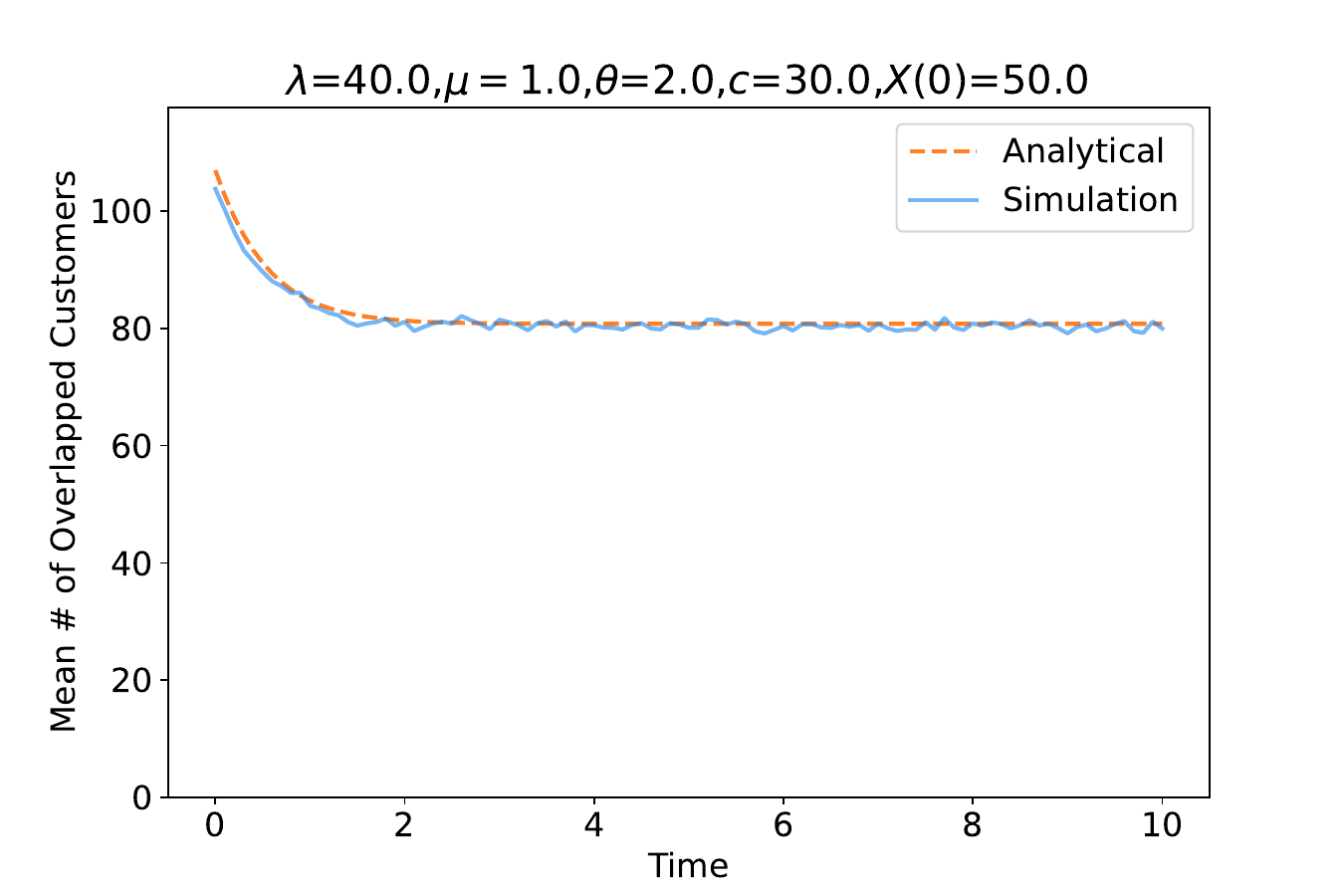}}
 \caption{Mean Number of Overlapping Customers (Analytical vs. Simulation). }
\label{Figure_5}
\end{figure}

Moving on to Figure \ref{Figure_5}, we explore the mean number of overlapping customers as a function of time. It is noteworthy that the approximations presented in Theorem \ref{overlap_fluid}, denoted as ``analytical" in the figures, show remarkable accuracy across all parameter values considered. The high level of precision maintained by these approximations reinforces their reliability and usefulness in practical settings. 

\begin{figure}[!htbp]
\centering
\subfloat[]{\includegraphics[scale=.31]{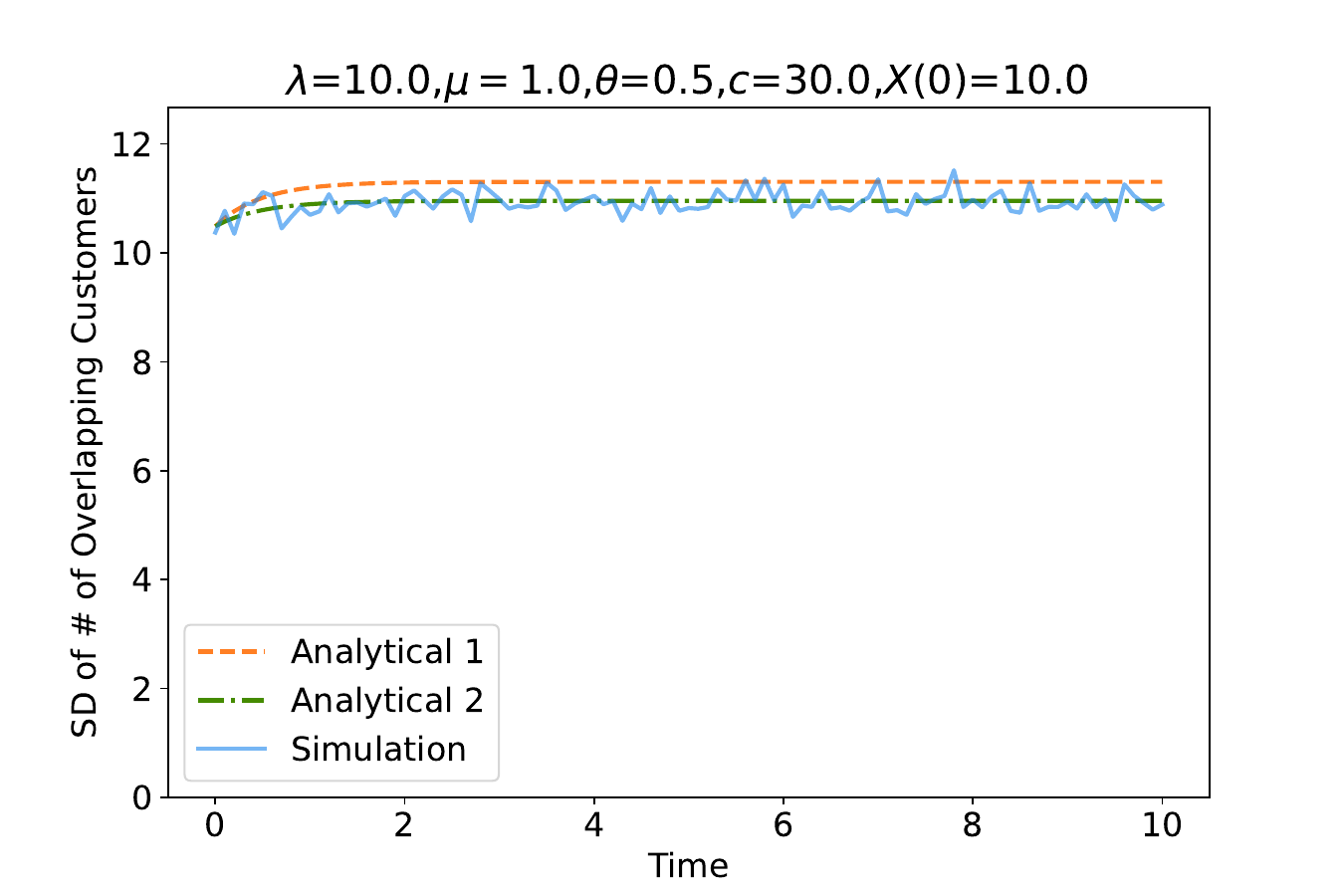}}
\subfloat[]{\includegraphics[scale=.31]{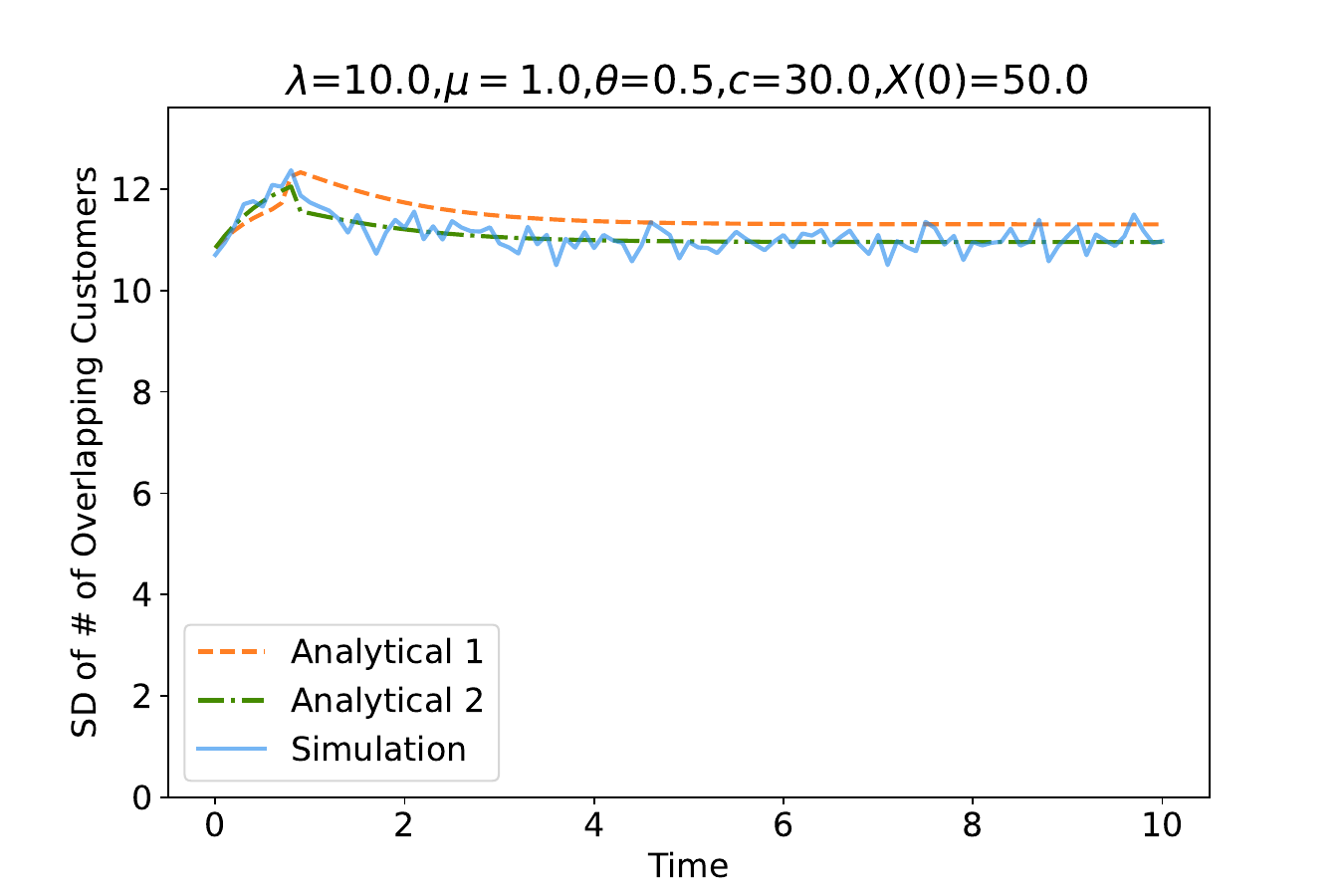}}
\\
\subfloat[]{\includegraphics[scale=.31]{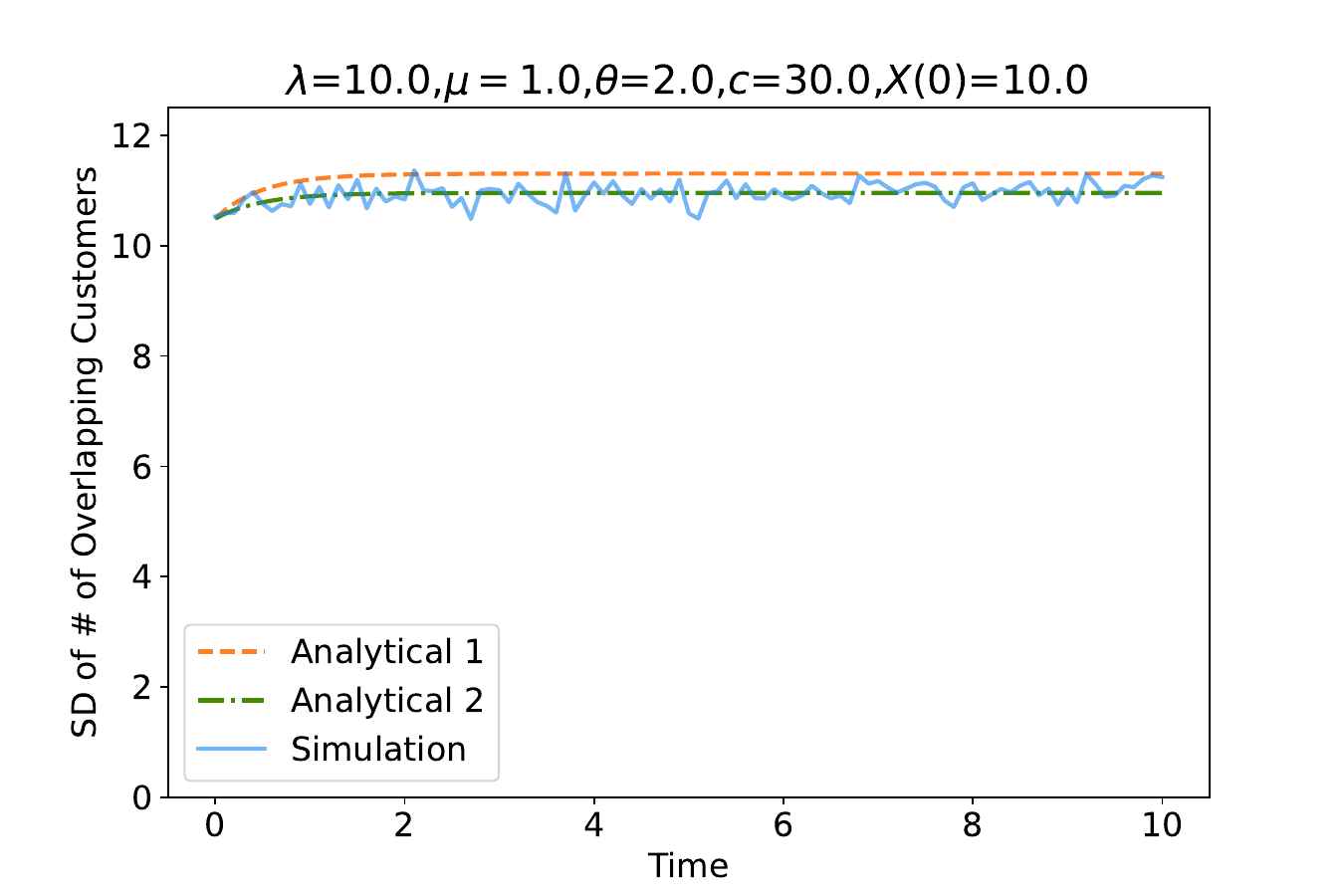}}
\subfloat[]{\includegraphics[scale=.31]{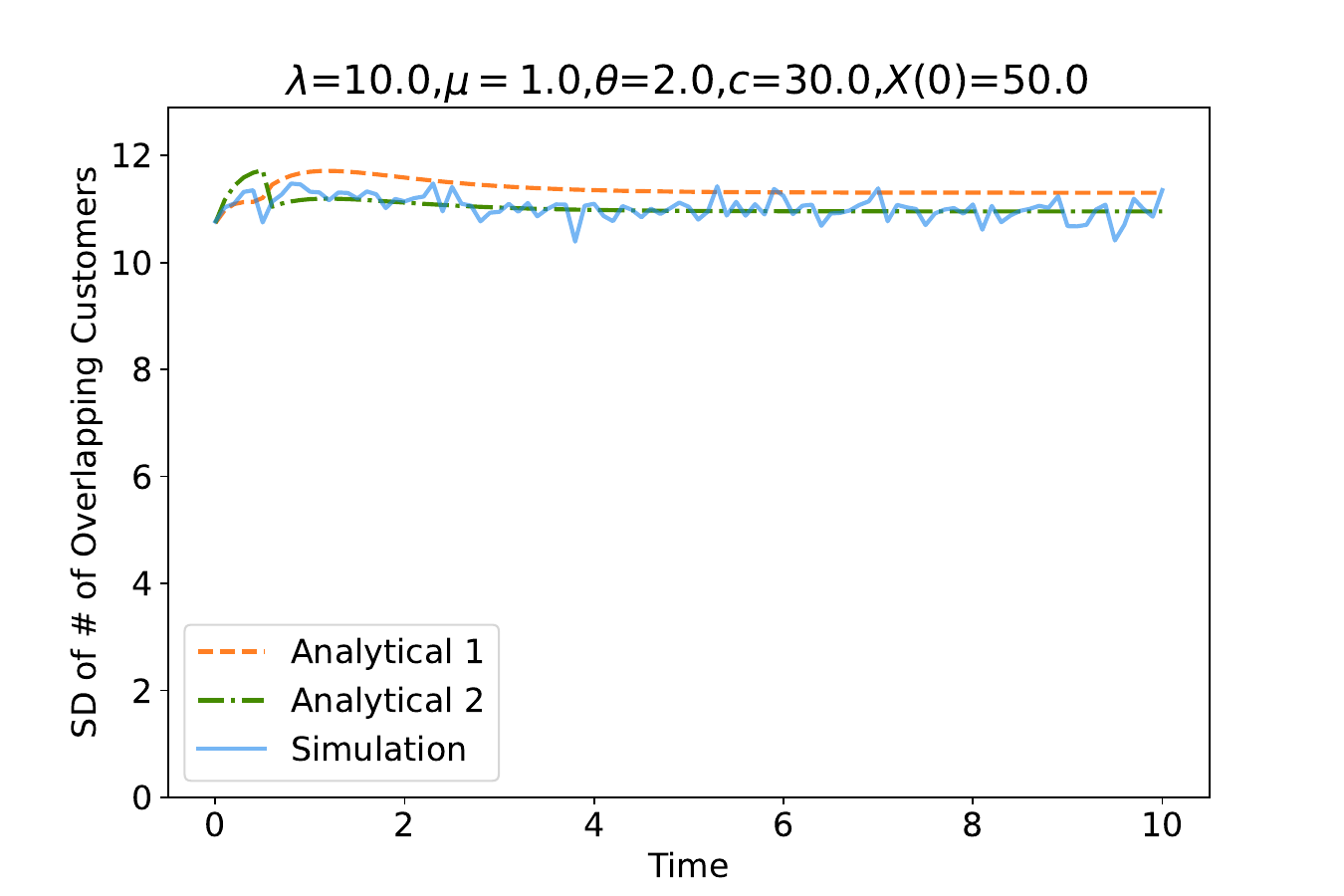}}
\\
\subfloat[]{\includegraphics[scale=.31]{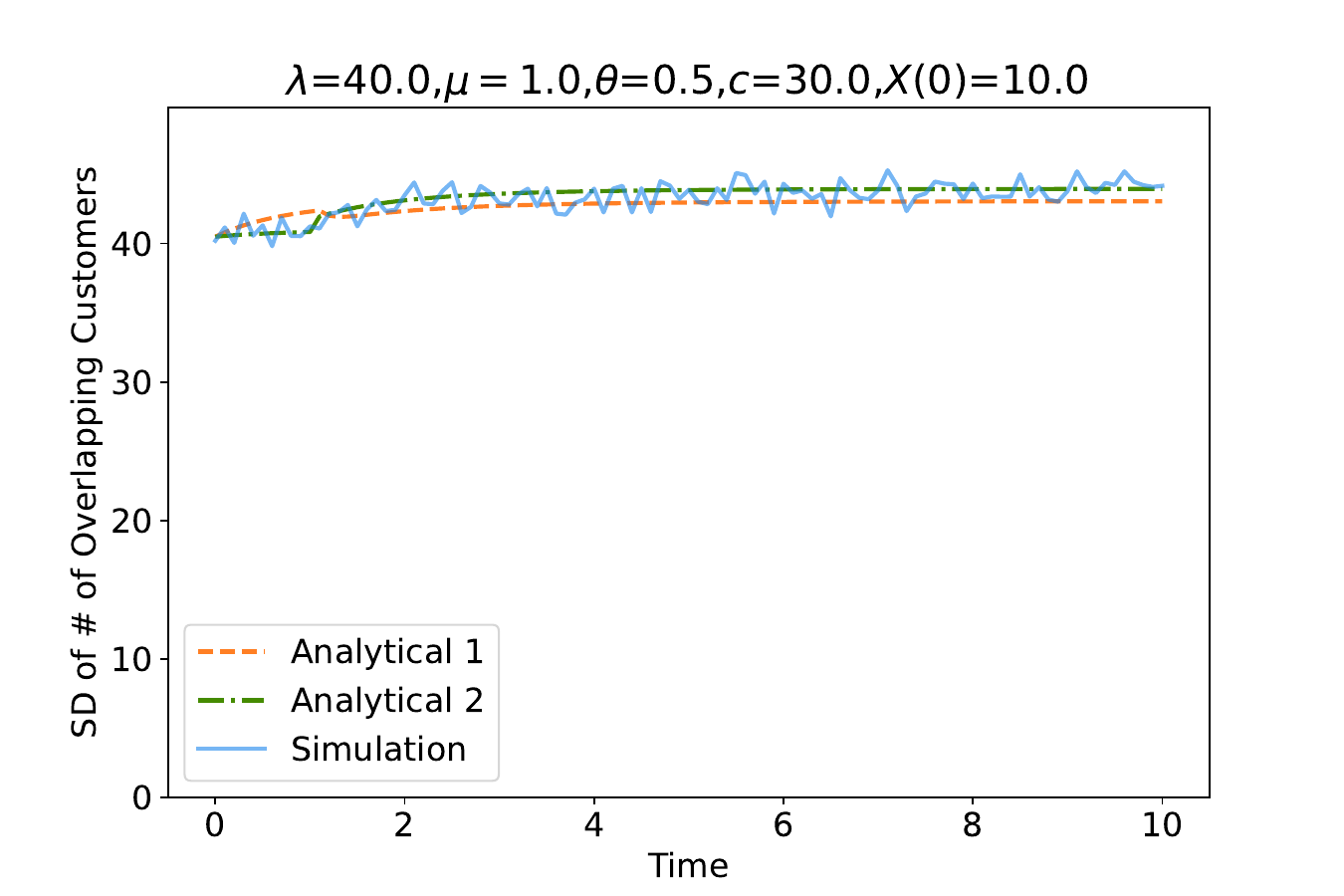}}
\subfloat[]{\includegraphics[scale=.31]{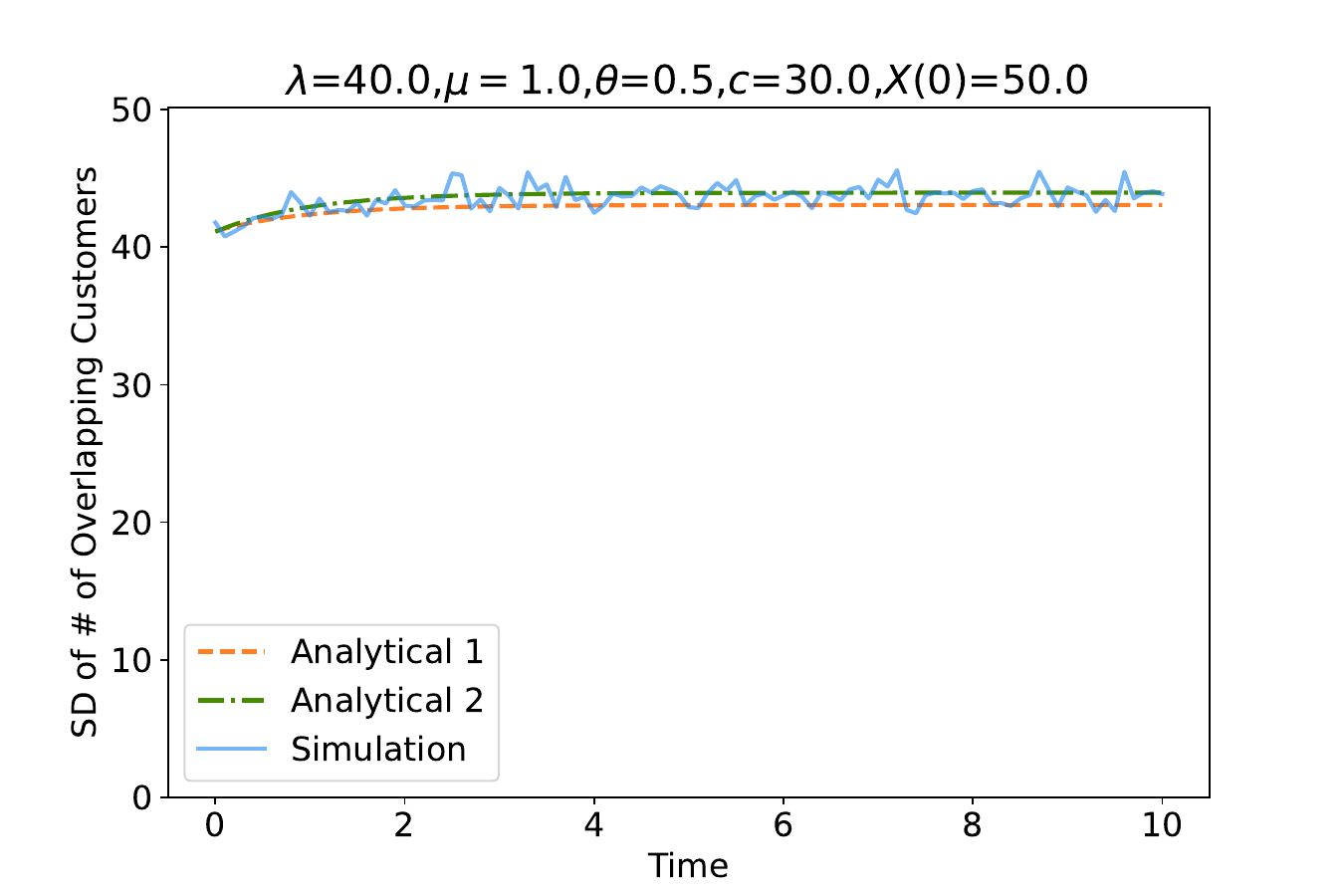}}
\\
\subfloat[]{\includegraphics[scale=.31]{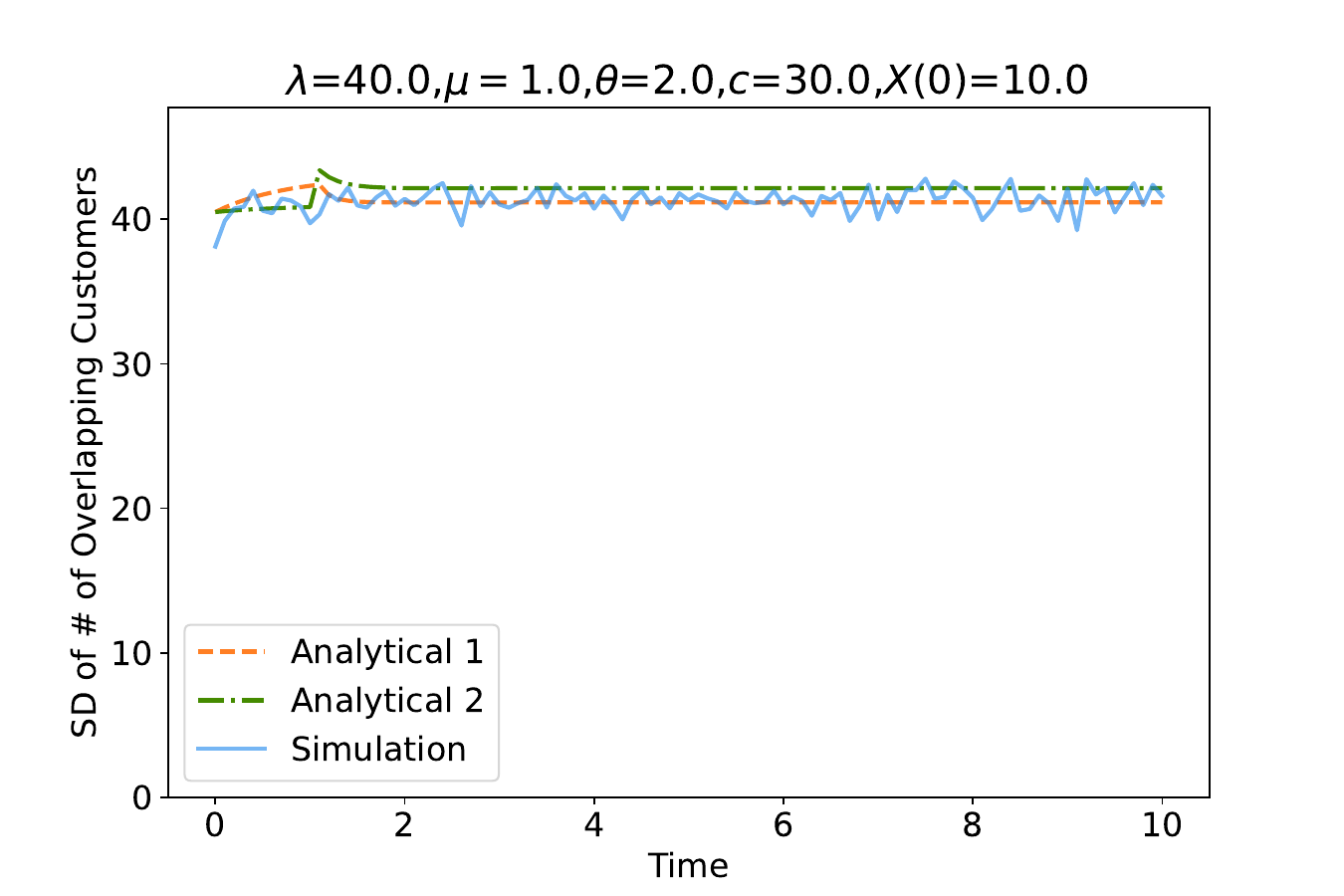}}
\subfloat[]{\includegraphics[scale=.31]{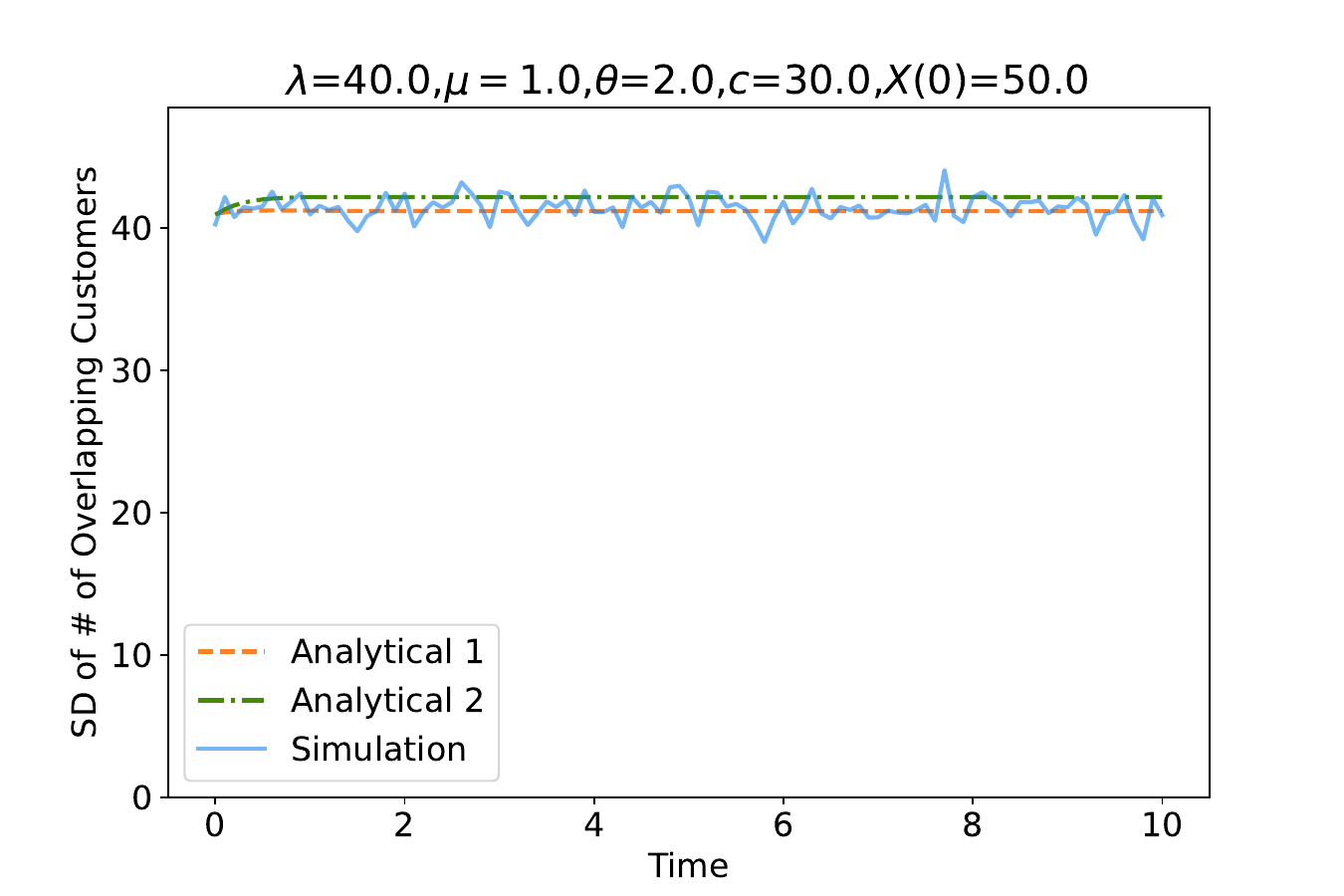}}
 \caption{Standard Deviation of Number of Overlapping Customers (Analytical vs. Simulation). }
\label{Figure_6}
\end{figure}

Figure \ref{Figure_6} plots the standard deviation of the number of overlapping customers over time for the aforementioned eight examples in Table \ref{Table_1}. Our investigation reveals that the two variance approximations presented in Equation \eqref{diffvar} (``analytical 1" in the figures) and Equation \eqref{var_wait} (``analytical 2" in the figures) both exhibit strong performance across all parameter values. The approximation based on the fluid and diffusion limits, as outlined in Equations \eqref{fldmean} and \eqref{diffvar}, consistently outperforms the approximation given by Equation \eqref{var_wait}. Consequently, with these robust approximations at our disposal, it becomes feasible to construct reliable prediction intervals for determining the customers that a virtual customer may overlap with.

By extending and refining our analysis using advanced approximation techniques derived from trigamma functions and fluid and diffusion limits, we have established a solid foundation for predicting and understanding various aspects of the virtual waiting time and customer overlap dynamics. These findings have significant implications for enhancing the efficiency and effectiveness of virtual customer service systems in a number of practical applications.


\section{Conclusion}\label{secConc}


 In this paper, we present a novel analysis of the mean and variance of the number of overlaps in the Erlang-A queue. Our contribution extends the current literature by considering both abandonment and a finite number of servers, thus providing a more realistic model. To achieve this, we employ a methodology based on the fluid and diffusion differential equations introduced by \citet{mandelbaum1998strong, mandelbaum1999waiting, massey2018dynamic}. Specifically, we derive exact expressions for these equations using the theory of linear differential equations. Moreover, we utilize these exact expressions to approximate the number of overlaps for a virtual customer that will not abandon. Our results show that our fluid and diffusion-based approximations offer reliable estimates of the mean and variance of the number of overlapping customers in the Erlang-A queue.

As a side result, we also present new approximations for the mean and variance of the waiting time in the Erlang-A queue. Notably, our approximations are functions of the digamma and trigamma functions, respectively. These approximations offer a significant improvement over existing results and can be used to enhance the performance of queueing systems in practical applications. Overall, our work contributes to a better understanding of the behavior of the Erlang-A queue and provides useful tools for its analysis in real-world scenarios.


 This work opens up several potential avenues for future research that could be valuable to pursue. Firstly, we could explore a more general queueing model with abandonment, such as the $G/G/C+G$ queue. Although some limit theorems exist for this model \citet{liu2012g, liu2014many}, the analysis of the virtual waiting time and its relationship with the queue length is currently unavailable. A thorough investigation of this relationship would provide insights into how the generality of the arrival, service, and abandonment processes could impact the number of overlaps.

Additionally, we could consider other types of queueing models, such as multidimensional network queueing models like those of \citet{liu2011network, liu2014stabilizing, pender2017approximating}, batch queueing models like those explored in \citet{pang2012infinite, daw2019distributions, daw2020non}, and even models with self-exciting arrivals \citet{koops2018infinite, daw2018queues, daw2020co}. These extensions would offer further opportunities to investigate the impact of various system parameters on the number of overlaps and the waiting time and they represent promising avenues for future research that we intend to pursue.


\section*{Acknowledgements}

Jamol Pender would like to acknowledge the gracious support of the National Science Foundation DMS Award \# 2206286. This work was also supported in part by the National Research Foundation of Korea (NRF) grants (No. 2021R1A2C1094699 and 2021R1A4A1031019) funded by the Korea government (Ministry of Science and ICT, MSIT).

\bibliographystyle{plainnat}
\bibliography{erlangA_refs}

\end{document}